\documentclass[10pt]{article}

\usepackage{amsfonts} 
\usepackage{amsmath,hhline,latexsym}
\usepackage{amssymb}

\usepackage[latin1]{inputenc}

\usepackage{color}
\usepackage{graphicx}
\usepackage{array}

\newtheorem{hyp}{Hypothesis}
\newtheorem{theorem}{\sc Theorem}[section]
\newtheorem{lemma}{\sc Lemma}[section]
\newtheorem{proposition}{\sc Proposition}[section]
\newtheorem{corollary}{\sc Corollary}[section]

\newtheorem{remark}{Remark}

\newcommand{\dis}{\displaystyle}
\newcommand{\eps}{\varepsilon}

\newcommand{\ph}{\varphi}

\newcommand{\Fin}{\hfill$\Box$}

\newcommand{\N}{\mbox{$I \kern -4pt N$}}
\newcommand{\Q}{\mbox{$Q \kern -8pt I$}}
\newcommand{\R}{\mbox{$I \kern -4pt R$}}
\newcommand{\C}{\mbox{$C \kern -8pt I$}}

\newcommand{\jnt}{\dis\int}
\newcommand{\jjntQT}{\jnt\!\!\!\!\jnt_{Q_{T}}}

\newcommand{\jjntGT}{\jnt\!\!\!\!\jnt_{\Gamma_{T}}}

\providecommand{\tabularnewline}{\\}

\def\R{{\bf R}}
\title{\textbf{Reconstruction of the solution and the source of hyperbolic equations from boundary measurements: mixed formulations}}

\author{ 
       \textsc{Nicolae C\^indea}\thanks{Laboratoire de Math\'ematiques, Universit\'e Blaise Pascal (Clermont-Ferrand 2), UMR CNRS 6620,  
	Campus Universitaire des C\'ezeaux, 3 place Vasarely, 63178, Aubi\`ere, France. E-mail: {\tt cindea@math.univ-bpclermont.fr.}}\quad 
	\and
	\textsc{Arnaud M\"unch}\thanks{Laboratoire de Math\'ematiques, Universit\'e Blaise Pascal (Clermont-Ferrand 2), UMR CNRS 6620, 
	Campus Universitaire des C\'ezeaux, 3 place Vasarely, 63178, Aubi\`ere, France.
	E-mail: {\tt arnaud.munch@math.univ-bpclermont.fr}.}}
	
\begin{document}


%



%



%
%
%

%

\maketitle

\begin{abstract}
We introduce a direct method allowing to solve numerically inverse type problems for linear hyperbolic equations. We first consider the reconstruction of the full solution of the wave equation posed in $\Omega\times (0,T)$ - $\Omega$ a bounded subset of $\mathbb{R}^N$ - from a partial boundary observation. We employ a least-squares technique and minimize the $L^2$-norm of the distance from the observation to any solution. Taking the hyperbolic equation as the main constraint of the problem, the optimality conditions are reduced to a mixed formulation involving both the state to reconstruct and a Lagrange multiplier. Under usual geometric optic conditions, we show the well-posedness of this mixed formulation (in particular the inf-sup condition) and then introduce a numerical approximation based on space-time finite elements discretization. We prove the strong convergence of the approximation and then discuss several examples for $N=1$ and $N=2$. The problem of the reconstruction of both the state and the source term is also addressed.  
\end{abstract}



\section{Introduction }
\label{sec:intro}

Let $\Omega$ be a bounded domain of $\mathbb{R}^N$ ($N\geq 1$) whose boundary $\partial\Omega$ is of class $C^2$ and let $T>0$. We define $Q_T:=\Omega\times (0,T)$, $\Sigma_T:=\partial\Omega\times (0,T)$ and denote by $\nu=\nu(x)$ the outward unit normal to $\Omega$ at any point $x\in \partial\Omega$. We are concerned with inverse type problems for the following linear hyperbolic equation 
\begin{equation}
 \label{eq:wave}
 \left\{
   \begin{array}{ll}
   y_{tt} - \nabla\cdot (c(x) \nabla y) + d(x, t) y= f, & \qquad (x, t) \in Q_T \\
   y = 0, & \qquad (x,t)\in \Sigma_T \\
   (y(\cdot, 0), y_t(\cdot,0)) = (y_0, y_1), & \qquad x \in \Omega.
   \end{array} 
 \right.
\end{equation}
We assume that $c\in C^1(\overline{\Omega}, \mathbb{R})$ with $c(x)\geq c_0>0$ in~$\overline{\Omega}$, $d \in L^\infty(Q_T)$, $(y_0,y_1) \in \boldsymbol{V}:=H_0^1(\Omega)\times L^2(\Omega)$ and $f\in L^2(Q_T)$.

   For any $(y_0,y_1)\in \boldsymbol{V}$ and any $f\in L^2(Q_T)$, there exists exactly one solution $y$ to (\ref{eq:wave}), with $y \in C^0([0, T]; H_0^1(\Omega)) \cap C^1([0, T]; L^2(\Omega))$ (see~\cite{JLL88,lionsmagenes72}). 
   In the sequel, for simplicity, we shall use the following notation:
\begin{equation}
\label{eq:L}
L\, y:=y_{tt}-\nabla\cdot (c(x) \nabla y) + d(x,t)y.
\end{equation}
Let now $\Gamma$ be any non empty open subset of $\partial\Omega$ and let $\Gamma_T:=\Gamma \times (0,T)\subset \Sigma_T$.
A typical inverse problem for (\ref{eq:wave}) is the following one :  from an \textit{observation} or \textit{measurement} of the normal derivative $y_{\nu,obs}$ in $L^2(\Gamma_T)$ on the sub-domain $\Gamma_T$, we want to recover a solution $y$ of the boundary value problem (\ref{eq:wave}) such that its normal derivative coincides with the observation on $\Gamma_T$. 

We denote by $\partial_{\nu}y=\nabla y\cdot \nu$ the normal derivative of $y$ on $\partial\Omega$. Introducing the operator $P:Z\to L^2(Q_T)\times L^2(\Gamma_T)$ defined by $P\,y:=(Ly,c(x)\partial_{\nu} y_{\vert \Gamma_T})$ and the Hilbert space $Z$ defined by (\ref{refZ})-(\ref{eq:pseta}), the problem is reformulated as : 
\begin{equation}
\label{IP}\tag{$IP$}
\text{\it find } y\in Z \text{ \it solution of } P\,y=(f,c(x) y_{\nu,obs}). 
\end{equation}

From the unique continuation property for (\ref{eq:wave}), if the set $\Gamma_T$ satisfies some geometric conditions and if $y_{\nu,obs}$ is a restriction to $\Gamma_T$ of a normal derivative of a solution of (\ref{eq:wave}), then the problem is well-posed in the sense that the state $y$ corresponding to the pair $(y_{\nu,obs},f)$ is unique. 



In view of the unavoidable uncertainties on the data $y_{\nu,obs}$ (coming from measurements, numerical approximations, etc), the problem (\ref{IP}) needs to be relaxed. In this respect, the most natural (and widely used in practice) approach consists in introducing the following extremal problem (of least-squares type)
\begin{equation}
\label{extremal_problem} \tag{\it LS}
\quad
\left\{
\begin{aligned}
&\textrm{minimize over } \boldsymbol{V} \quad J(y_0,y_1):=\frac{1}{2}  \Vert c(x)(\partial_{\nu} y-y_{\nu,obs}) \Vert^2_{L^2(\Gamma_T)}\\
& \textrm{where} \quad y \quad \textrm{solves} \quad (\ref{eq:wave}),
\end{aligned}
\right.
\end{equation}
since $y$ is uniquely and fully determined from $f$ and the data $(y_0,y_1)$. Here the constraint $c(x)(\partial_{\nu} y-y_{\nu,obs})=0$ in $L^2(\Gamma_T)$ is relaxed; however, if $y_{\nu,obs}$ is a restriction to $\Gamma_T$ of the normal derivative of a solution of (\ref{eq:wave}), then problems (\ref{extremal_problem}) and (\ref{IP}) coincide. A minimizing sequence for $J$ in $\boldsymbol{V}$ is easily defined in term of the solution of an auxiliary adjoint problem. Apart from a possible low decrease of the sequence near extrema, the main drawback, when one wants to prove the convergence of a discrete approximation is that, it is not in general possible to minimize over a discrete subspace of $\{y; Ly-f=0\}$ subject to the equality (in $L^2(Q_T)$) $Ly-f=0$. Therefore, a classical trick consists first in discretizing the functional $J$ and the system (\ref{eq:wave}); this raises the issue of uniform coercivity property (typically here some uniform discrete observability inequality for the adjoint solution) of the discrete functional with respect to the approximation parameter. As far as we know, this delicate issue has received answers only for specific and somehow academic situations (uniform Cartesian approximation of $\Omega$, constant coefficients in (\ref{eq:wave}), etc). We refer to \cite{NC-AM-mixedwave,Glo08,komornikloreti,munch05} and the references therein.

More recently, a different method to solve inverse type problems like (\ref{IP}) has emerged and use the so called Luenberger type observers:  roughly, this method consists in defining, from the observation on $\Gamma_T$, an auxiliary boundary value problem whose solution possesses the same asymptotic behavior in time than the solution of (\ref{eq:wave}): the use of the reversibility of the hyperbolic equation then allows to reconstruct the initial data $(y_0,y_1)$. We refer to \cite{cindea_moireau,ramdani2010} and the references therein. However, for the same reasons, from a numerically point of view, these methods require to prove uniform discrete observability properties. 

In a series of works, Klibanov and co-workers use different approaches to solve inverse problems (we refer to \cite{Klibanov-book} and the references therein): they advocate in particular the quasi-reversibility method which reads as follows : for any $\eps>0$, find $y_{\eps}\in \mathcal{A}$ the solution of 
\begin{equation}\label{FV}\tag{$QR_{\eps}$}
 \langle Py_\eps,P\overline{y} \rangle_{L^2(Q_T)\times L^2(\Gamma_T)} + \eps \langle y_{\eps},\overline{y} \rangle_{\mathcal{A}} = \left\langle(f,c(x)y_{\nu,obs}), P\overline{y}\right\rangle_{L^2(Q_T)\times L^2(\Gamma_T)}, \quad \forall \ \overline y \in \mathcal{A},
\end{equation}
where $\mathcal{A}$ denotes a Hilbert space subset of $L^2(Q_T)$ so that $Py\in L^2(Q_T)\times L^2(\Gamma_T)$ for all $y\in \mathcal{A}$ and $\eps>0$ is a Tikhonov like parameter which ensures the well-posedness. We refer for instance to \cite{clason_klibanov} where the lateral Cauchy problem for the wave equation with non constant diffusion is addressed within this method. Remark that (\ref{FV}) can be viewed as a least-squares problem since the solution $y_{\eps}$ minimizes over $\mathcal{A}$ the functional $y\to \Vert P y - (f,c(x)y_{\nu,obs})\Vert^2_{L^2(Q_T)\times L^2(\Gamma_T)}+ \eps \Vert y\Vert^2_{\mathcal{A}}$. Eventually, if $y_{\nu,obs}$ is the normal derivative of a restriction to $\Gamma_T$ of a solution of (\ref{eq:wave}), 
the corresponding $y_{\eps}$ converges in $L^2(Q_T)$ toward to the solution of (\ref{IP}) as $\eps\to 0$. There, unlike in Problem (\ref{extremal_problem}), the unknown is the state variable $y$ itself (as it is natural for elliptic equations) so that any standard numerical methods based on a conformal approximation of the space $\mathcal{A}$ together with appropriate observability inequalities allow to obtain a convergent approximation of the solution. In particular, there is no need to prove discrete observability inequalities. We refer to the book \cite{BeilinaKlibanov14}. We also mention \cite{bourgeoisbis,bourgeois2010} where a similar technique has been used recently to solve the inverse obstacle problem associated to the Laplace equation, which consists in finding an interior obstacle from boundary Cauchy data. 

In the spirit of the works \cite{bourgeois2010,clason_klibanov,Klibanov-book}, we explore the direct resolution of the optimality conditions associated to the extremal problem (\ref{extremal_problem}), without Tikhonov parameter while keeping $y$ as the unknown of the problem. This strategy, which avoids any iterative process, has been successfully applied in the close context of the exact controllability of (\ref{eq:wave}) in \cite{NC-AM-mixedwave} and \cite{CC-NC-AM,NC-EFC-AM}. The idea is to take into account the state constraint $Ly-f=0$ with a Lagrange multiplier. This allows to derive explicitly the optimality system associated to (\ref{extremal_problem}) in term of an elliptic mixed formulation and therefore reformulate the original problem. Well-posedness of such new formulation is related to a unique continuation property for the hyperbolic equation (\ref{eq:wave}).   

From the observation $y_{\nu,obs}$, we also address the simultaneous reconstruction problem of the source term $f$ and the solution $y$: 
\begin{equation}
\label{IPf}\tag{$IP_f$}
\text{\it find } (y,f)\in Z\times L^2(Q_T) \text{ \it solution of } P\,y=(f,c(x) y_{\nu,obs}), 
\end{equation}
where $(y,f)$ solves (\ref{eq:wave}). Without additional assumption on $f\in L^2(Q_T)$, the pair $(y,f)$
solution of (\ref{IPf}) is not unique: consider for instance a source term $f$ supported in a set which is near $\Omega\times \{T\}$: from the finite propagation of the solution, the source $f$ will not affect the solution $y$ on $\Gamma_T$. On the other hand, a result of Yamamoto and Zhang in \cite{yamamoto} asserts that the uniqueness holds true if the source takes the form $f(x,t)=\sigma(t)\mu(x)$, where the smooth time part $\sigma$ is given and the spatial part $\mu$ is a $H^{-1}(\Omega)$ function. 

We adapt in this work the arguments of \cite{NC-AM-IPboundary} where the observation is distributed in $Q_T$. The outline is as follow. In Section \ref{recovering_y}, we consider the least-squares problem (\ref{P}) and reconstruct the solution of the hyperbolic equation from a partial observation localized on a subset $\Gamma_T$ of $\Sigma_T$. For that, in Section \ref{sec2_direct}, we associate to (\ref{P}) the equivalent mixed formulation (\ref{eq:mf}) which relies on the optimality conditions of the problem. Assuming that $\Gamma_T$ satisfies the classical geometric optic condition (Hypothesis 1, see (\ref{iobs})), we then show the well-posedness of this mixed formulation, in particular, we check the Babuska-Brezzi inf-sup condition (see Theorem \ref{th:mf}). Interestingly, we also derive in Section \ref{sec2_dual} an equivalent dual extremal problem, which reduces the determination of the state $y$ to the minimization of a strongly elliptic functional with respect to the Lagrange multiplier. In Section \ref{recovering_y_mu}, we use the uniqueness result \cite[Theorem 2.1]{yamamoto} of Yamamoto and Zhang and apply the same procedure to recover from a partial observation both the state and the spatial part $\mu$ of the source term assumed in $H^{-1}(\Omega)$. Section \ref{sec_numer} is devoted to the numerical approximation, through a conformal space-time finite element discretization. The strong convergence of the approximation $\{y_h,\mu_h\}$ is proved as the discretization parameter $h$ tends to zero. In particular, we discuss the discrete  inf-sup property of the mixed formulation. We present numerical experiments in Section \ref{sec_experiment} for $\Omega=(0,1)$ and $\Omega \subset\mathbb{R}^2$, an check the robustness and convergence of the approximations in agreement with the theoretical part. Section \ref{sec_conclusion} concludes with some perspectives: in particular, we highlight that the parabolic case can be treated in a similar way.


\section[Recovering the solution from a partial observation]{Recovering the solution from a partial observation: a mixed re-formulation of the problem}\label{recovering_y}

In this section, assuming that the initial data $(y_0,y_1)\in \boldsymbol{V}$ are unknown, we address the inverse problem (\ref{IP}). Without loss of generality, in view of the linearity of the system (\ref{eq:wave}), we assume that the source term $f$ is zero.  

\par\noindent
We consider the non empty vector space $Z$ defined by 
\begin{equation}
Z:=\{y: y\in C([0,T], H_0^1(\Omega))\cap C^1([0,T], L^2(\Omega)), Ly\in L^2(Q_T)\}  \label{refZ}
\end{equation}
and then first recall that $\partial_{\nu}y \in L^2(\Sigma_T)$ for any $y\in Z$: precisely (see \cite[Theorem 4.1, Ch 1]{JLL88}), there exists a constant $C_T>0$ such that the following holds : 
\begin{equation}
\Vert c(x)\partial_{\nu} y\Vert^2_{L^2(\Gamma_T)} \leq C_T \biggl(\Vert (y(\cdot,0),y_t(\cdot,0))\Vert^2_{\boldsymbol{V}} + \Vert Ly\Vert^2_{L^2(Q_T)}  \biggr), \quad \forall y\in Z.  \label{dependancecontinu}
\end{equation}

We then introduce the following hypothesis :

\begin{hyp}
There exists a constant $C_{obs}=C(\omega,T,\Vert c\Vert_{C^1(\overline{\Omega})},\Vert d\Vert_{L^{\infty}(Q_T)})$ such that the following estimate holds : 
\begin{equation}
\Vert y(\cdot,0),y_t(\cdot,0)\Vert^2_{\boldsymbol{V}} \leq C_{obs} \biggl( \Vert c(x)\partial_{\nu} y\Vert^2_{L^2(\Gamma_T)}+ \Vert Ly\Vert^2_{L^2(Q_T)}  \biggr), \quad \forall y\in Z.  \label{iobs}\tag{$\mathcal{H}$}
\end{equation}
\end{hyp}

Condition (\ref{iobs}) is a generalized observability inequality for the solution of the hyperbolic equation (\ref{eq:wave}). For constant coefficients, this estimate is known to hold if the triplet ($\Gamma,T,\Omega$) satisfies a geometric optic condition. We refer to \cite{BLR} for the case of constant velocity $c$. In particular, $T$ must be large enough. In the one-dimensional case, for non constant velocity $c$ and potential $d$, we refer to \cite{NC-EFC-AM} and the references therein.

Then, within this hypothesis, for any $\eta>0$, we define on $Z$ the bilinear form 
\begin{equation}
\label{eq:pseta}
\begin{aligned}
\langle y,\overline{y}\rangle_Z:= & \jjntGT c^2(x)\,\partial_{\nu} y\,\partial_{\nu}\overline{y}\,d\sigma dt + \eta \jjntQT Ly\, L\overline{y}  \,dxdt  \quad \forall y,\overline{y}\in Z
\end{aligned}
\end{equation}
and we denote the corresponding semi-norm $\Vert y\Vert_Z:=\sqrt{\langle y,y \rangle_Z}$.
We have the following result : 
\begin{lemma}\label{le:ZHilbert}
Under the hypothesis $(\mathcal{H})$, the space $Z$ is a Hilbert space with the scalar product $\langle \cdot,\ \cdot \rangle_Z$ defined by (\ref{eq:pseta}).
\end{lemma}
\textsc{Proof-}The two main properties we need to verify are that the semi-norm associated to this inner product $\| \cdot \|_Z$ is indeed a norm, and that $Z$ is closed with respect to this norm. The first property is a direct consequence of the inequality $(\mathcal{H})$. 

To check the second property, let us consider a convergent sequence $\{ z_k\}_{k\geq 1} \subset Z$ such that $z_k \to z$ in the norm $\| \cdot \|_Z$. We have to see that $z \in Z$. 
As a consequence of  $(\mathcal{H})$, there exist $(z_0,z_1)\in \boldsymbol{V}$ and $f\in L^2(Q_T)$ such that $(z_k(\cdot,0),z_{k,t}(\cdot, 0))\to (z_0,z_1)$ in $\boldsymbol{V}$ and $Lz_k \to f$ in $L^2(Q_T)$. Therefore, $z_k$ can be considered as a sequence of solutions of the hyperbolic equation with convergent initial data and second hand term $Lz_k \to f$. 

By the continuous dependence of the solutions of the wave equation on the data, $z_k \to z$ in $C([0,T];H_0^1(\Omega))\cap C^1([0,T];L^2(\Omega))$, where $z$ is the solution of the hyperbolic equation with initial data $(z_0,z_1)\in \boldsymbol{V}$ and second hand term $Lz=f \in L^2(Q_T)$. Therefore, in view of (\ref{dependancecontinu}), $z\in Z$. \Fin

\
\par\noindent
We consider the following extremal problem : 
\begin{equation}
\label{P}
\tag{$\mathcal{P}$}
\left\{
\begin{aligned}
& \inf  J(y):=  \frac{1}{2}\Vert c(x)(\partial_{\nu} y-y_{\nu,obs})\Vert^2_{L^2(\Gamma_T)}, \\
& \textrm{subject to}\quad y\in W
\end{aligned}
\right.
\end{equation}
where $W$ is the closed subspace of $Z$ defined by 
$$
W:=\{y\in Z; \,  Ly=0 \,\, \textrm{in}\,\, L^2(Q_T)\}
$$
and endowed with the norm of $Z$. 

The extremal problem (\ref{P}) is well posed : the functional $J$ is continuous over $W$, is strictly convex and is such that $J(y)\to +\infty$ as $\Vert y\Vert_W\to \infty$. Note also that the solution of (\ref{P}) in $W$ does not depend on $\eta$. 

We recall that from the definition of $Z$, $Ly$ belongs to $L^2(Q_T)$. Furthermore, the uniqueness of the solution is lost if the hypothesis (\ref{iobs}) is not fulfilled, for instance if $T$ is not large enough.  Eventually, from (\ref{iobs}), the solution $y$ in $Z$ of (\ref{P}) satisfies $(y(\cdot,0),y_t(\cdot,0))\in \boldsymbol{V}$, so that Problem (\ref{P}) is equivalent to the minimization of $J$ with respect to $(y_0,y_1)\in \boldsymbol{V}$ as in problem (\ref{IP}), Section 1. 


\subsection{Direct approach} \label{sec2_direct}

In order to solve (\ref{P}), we have to deal with the constraint equality which appears in the space $W$. Proceeding as in \cite{NC-AM-mixedwave}, we introduce a Lagrange multiplier $\lambda\in L^2(Q_T)$ and the following mixed formulation: find $(y, \lambda)\in Z\times L^2(Q_T)$ solution of 
\begin{equation} \label{eq:mf}
\left\{
\begin{array}{rcll}
\noalign{\smallskip} a(y, \overline{y}) + b(\overline{y}, \lambda) & = & l(\overline{y}), & \qquad \forall \overline{y} \in Z \\
\noalign{\smallskip} b(y, \overline{\lambda}) & = & 0, & \qquad \forall \overline{\lambda} \in L^2(Q_T),
\end{array}
\right.
\end{equation}
where
\begin{align}
\label{eq:a} & a : Z \times Z \to \mathbb{R},  \quad a(y,\overline{y}) := \jjntGT c^2(x)\,\partial_{\nu}y\,\partial_{\nu}\overline{y}\,d\sigma dt,
\\
\label{eq:b} & b: Z \times L^2(Q_T)  \to \mathbb{R},  \quad b(y,\lambda) := \jjntQT \lambda\, Ly  \,dxdt,\\
\label{eq:l} & l: Z \to \mathbb{R},  \quad l(y) := \jjntGT c^2(x)\,y_{\nu,obs} \,\partial_{\nu}y\, d\sigma dt. 
\end{align}
System (\ref{eq:mf}) is nothing else than the optimality system corresponding to the extremal problem (\ref{P}). Precisely, the following result holds: 
\begin{theorem}\label{th:mf} Under the hypothesis (\ref{iobs}), 
\begin{enumerate}
\item The mixed formulation (\ref{eq:mf}) is well-posed.
\item The unique solution $(y, \lambda) \in Z\times L^2(Q_T)$ to (\ref{eq:mf}) is the unique saddle-point of the Lagrangian $\mathcal{L}:Z\times L^2(Q_T)\to \mathbb{R}$ defined by
\[
\begin{aligned}
\mathcal{L}(y, \lambda):=& \frac{1}{2}a(y,y) + b(y,\lambda)- l(y).
\end{aligned}
\]
\item We have the estimate 
\begin{align}
& \Vert y\Vert_Z= \left( \Vert c(x)\,\partial_{\nu} y\Vert_{L^2(\Gamma_T)}^2 + \eta \Vert Ly \Vert_{L^2(Q_T)}^2\right)^\frac12 \leq \Vert c(x)\,y_{\nu,obs}\Vert_{L^2(\Gamma_T)},  \nonumber\\
& \Vert \lambda\Vert_{L^2(Q_T)} \leq 2\sqrt{C_{\Omega,T}+\eta}  \Vert c(x)\,y_{\nu,obs}\Vert_{L^2(\Gamma_T)}.  \label{estimate_lambda} 
\end{align}
\end{enumerate}
\end{theorem}

\textsc{Proof-} The proof is based on classical results for saddle point problems  (see \cite{brezzi_new}, chapter 4). 

We easily obtain the continuity of the bilinear form $a$ over $Z\times Z$, the continuity of bilinear $b$ over $Z\times L^2(Q_T)$ and the continuity of the linear form $l$ over $Z$. In particular, we get 
\begin{equation}
\Vert l\Vert_{Z^{\prime}}\leq \Vert c(x)\,y_{\nu,obs}\Vert_{L^2(\Gamma_T)}, \qquad \Vert a\Vert_{(Z\times Z)^{\prime}}\leq 1, \quad \Vert b\Vert_{(Z\times L^2(Q_T))^{\prime}} \leq  \eta^{-1/2}. \label{simple_ineq}
\end{equation}
Moreover, the kernel $\mathcal{N}(b)=\{y\in Z;\  b(y,\lambda)=0 \quad \forall \lambda\in L^2(Q_T)\}$ coincides with $W$: we easily get 
\begin{equation}
a(y,y)=\Vert y \Vert^2_Z, \quad \forall y\in \mathcal{N}(b)=W.  \nonumber
\end{equation}
Therefore, in view of \cite[Theorem 4.2.2]{brezzi_new}, it remains to check the inf-sup constant property : $\exists \delta>0$ such that
\begin{equation}
\inf_{\lambda\in L^2(Q_T)} \sup_{y \in Z}  \frac{b(y,\lambda)}{\Vert y \Vert_Z \Vert \lambda\Vert_{L^2(Q_T)}} \geq \delta. \nonumber
\end{equation}
We proceed as follows. For any fixed $\lambda\in L^2(Q_T)$, we define $y_0\in Z$ as the unique solution of 
\begin{equation}
Ly_0=\lambda \,\,\,\textrm{in} \,\,\, Q_T,  \quad (y_0(\cdot,0),y_{0,t}(\cdot,0))=(0,0)\,\,\, \textrm{in}\,\,\, \Omega ,\quad y_0=0\,\,\, \textrm{on}\,\,\, \Sigma_T.   \label{infsup_y0}
\end{equation}  
We get $b(y_0,\lambda)=\Vert \lambda\Vert^2_{L^2(Q_T)}$
and
\begin{equation}
\Vert y_0\Vert^2_Z= \Vert c(x)\,\partial_{\nu} y_0\Vert^2_{L^2(\Gamma_T)} + \eta \Vert \lambda \Vert_{L^2(Q_T)}^2.   \nonumber
\end{equation}
Using (\ref{dependancecontinu}), the estimate $\Vert c(x)\partial_{\nu} y_0\Vert_{L^2(\Gamma_T)} \leq \sqrt{C_{\Omega,T}} \Vert \lambda \Vert_{L^2(Q_T)}$  
 implies that 
\[
\sup_{y\in Z}  \frac{b(y,\lambda)}{\Vert y \Vert_Z \Vert \lambda \Vert_{L^2(Q_T)}}\geq  \frac{b(y_0,\lambda)}{\Vert y_0 \Vert_Z \Vert \lambda \Vert_{L^2(Q_T)}}\geq \frac{1}{\sqrt{C_{\Omega,T}+\eta}}>0
\]
leading to the result with $\delta=(C_{\Omega,T}+\eta)^{-1/2}$. 

The third point is the consequence of classical estimates  (see \cite{brezzi_new}, Theorem 4.2.3.) : 
\begin{equation}
\Vert y\Vert_{Z}\leq \frac{1}{\alpha_0} \Vert l\Vert_{Z^{\prime}}, \quad \Vert \lambda \Vert_{L^2(Q_T)} \leq \frac{1}{\delta}\biggl(1+\frac{\Vert a \Vert}{\alpha_0} \biggr) \Vert l\Vert_{Z^{\prime}} \nonumber
\end{equation}
where 
\begin{equation}
\alpha_0:=\inf_{y\in \mathcal{N}(b)}  \frac{a(y,y)}{\Vert y\Vert^2_Z}.  \label{eq:alpha0}
\end{equation}
Estimates (\ref{simple_ineq}) and the equality $\alpha_0=1$ lead to the results.  Eventually, from (\ref{simple_ineq}), we obtain that 
\begin{equation}
\Vert \lambda\Vert_{L^2(Q_T)} \leq \frac{2}{\delta} \Vert c(x)\,y_{\nu,obs}\Vert_{L^2(\Gamma_T)}   \nonumber
\end{equation}
and that $\delta \geq (C_{\Omega,T}+\eta)^{-1/2}$ to get (\ref{estimate_lambda}). \Fin

In practice, it is very convenient to "augment" the Lagrangian (see \cite{fortinglowinski}) and consider instead the Lagrangian $\mathcal{L}_r$ defined for any $r>0$ by 
\begin{equation}
\begin{aligned}
& \mathcal{L}_r(y,\lambda):=\frac{1}{2}a_r(y,y)+b(y,\lambda)-l(y), \\
& a_r(y,y):=a(y,y)+r\Vert Ly\Vert^2_{L^2(Q_T)}. 
\end{aligned}
 \nonumber
\end{equation}
Since $a_r(y,y)=a(y,y)$ on $W$, the Lagrangian $\mathcal{L}$ and $\mathcal{L}_r$ share the same saddle-point. The positive real $r$ is an augmentation parameter. 

\begin{remark}\label{rk_lambda_sys}
Assuming additional hypotheses on the regularity of the solution $\lambda$, precisely $L\lambda\in L^2(0,T, H^{-1}(\Omega))$ and $(\lambda,\lambda_t)_{\vert t=0,T} \in L^2(\Omega)\times H^{-1}(\Omega)$,
we easily check, by writing the optimality conditions for $\mathcal{L}$, that the multiplier $\lambda$ satisfies the following relations :
%
%
\begin{equation}
\label{system_lambda}
\left\{
\begin{aligned}
& L\lambda=0 &&\quad \textrm{in}\quad Q_T , \quad \\
& \lambda= c(x)\left(\partial_{\nu}y-y_{\nu,obs}\right) &&\quad \textrm{on}\quad \Gamma_T, \\
&  \lambda=0 &&\quad\textrm{on}\quad \Sigma_T\setminus\Gamma_T, \\
& \lambda=\lambda_{t}=0 &&\quad\textrm{on}\quad \Omega \times \{0,T\}.
\end{aligned}
\right.
\end{equation}
Therefore, $\lambda$ (defined in the weak sense) is a exact null controlled solution of the hyperbolic equation (\ref{eq:wave}) through the boundary control $(\partial_{\nu} y-y_{\nu,obs})\, 1_{\Gamma_T}\in L^2(\Gamma_T)$. 
\begin{itemize}
\item  If $y_{\nu,obs}$ is the normal derivative of a solution of (\ref{eq:wave})  restricted to $\Gamma_T$, then the unique multiplier $\lambda$ must vanish almost everywhere. In that case, 
we have  
\[
\sup_{\lambda\in L^2(Q_T)}\inf_{y\in Z} \mathcal{L}_r(y,\lambda) = \inf_{y\in Z} \mathcal{L}_r(y,0)=\inf_{y\in Z} J_r(y)
\]
with 
\begin{equation}
\label{def_Jyr}
J_r(y):=\frac{1}{2}\Vert c(x)(\partial_{\nu}y-y_{\nu,obs})\Vert^2_{L^2(\Gamma_T)} + \frac{r}{2}\Vert Ly\Vert^2_{L^2(Q_T)}.
\end{equation}
The corresponding variational formulation is then : find $y\in Z$ such that 
\begin{equation}
a_r(y,\overline{y})=\jjntGT c^2(x) \partial_{\nu}y\,\partial_{\nu}\overline{y} \, d\sigma dt + r\jjntQT  Ly \, L\overline{y} \, dxdt = l(\overline{y}), \quad \forall \overline{y}\in Z.  \nonumber
\end{equation}
\item In the general case, the mixed formulation can be rewritten as follows: find $(y, \lambda)\in Z\times L^2(Q_T)$ solution of 
\begin{equation*}
\left\{
\begin{aligned}
\langle P_r y, P_r \overline{y}\rangle_{L^2(Q_T)\times L^2(\Gamma_T)} + \langle L\overline{y},\lambda \rangle_{L^2(Q_T)}   & = \langle (0,c(x)y_{\nu,obs}),P_r \overline{y} \rangle_{L^2(Q_T)\times L^2(\Gamma_T)}, \quad \forall \overline{y}\in Z, \\
\langle Ly,\overline{\lambda} \rangle_{L^2(Q_T)}   & = 0, \quad \forall \overline{\lambda}\in L^2(Q_T)
\end{aligned}
\right.
\end{equation*}
with $P_r:Z\times L^2(Q_T)\times L^2(Q_T)$ defined by $P_r y:=(\sqrt{r} L\,y,c(x)\,\partial_{\nu} y_{\vert \Gamma_T})$. This formulation may be seen as a generalization of the quasi-reversibility formulation (\ref{FV}), where the variable $\lambda$ is adjusted automatically (while the choice of the parameter $\eps$ in (\ref{FV}) is in general a delicate issue).
\end{itemize}
\end{remark}

System (\ref{system_lambda}) can be used to define a equivalent saddle-point formulation, very suitable at the numerical level.  Precisely, we define - in view of (\ref{system_lambda}) - the space $\Lambda$ by 
$$
\nonumber
\begin{aligned}
\Lambda:=\{\lambda: \lambda\in C([0,T]; & L^2(\Omega))\cap C^1([0,T]; H^{-1}(\Omega)), \\
&L\lambda\in L^2([0, T]; H^{-1}(\Omega)), \lambda(\cdot,0)=\lambda_t(\cdot,0)=0, \lambda_{|\Gamma_T}\in L^2(\Gamma_T)\}.
\end{aligned}
$$
Similarly to Lemma \ref{le:ZHilbert}, we can prove that $\Lambda$ is a Hilbert space endowed with the following inner product
\[
\langle \lambda, \ \overline\lambda \rangle_\Lambda := \int_0^T \langle L\lambda(t), L \overline \lambda(t) \rangle_{H^{-1}(\Omega)} dt +\jjntGT c^2(x) \lambda \overline \lambda d\sigma dt, \qquad \forall \ \lambda, \ \overline\lambda \in \Lambda
\]
using notably, that the elements of the non-empty vector space $\Lambda$ satisfy the inequality 
\begin{equation}
\Vert \lambda\Vert_{L^2(Q_T)}\leq C_{\Omega,T} \sqrt{\langle \lambda,\lambda \rangle_{\Lambda}}  \label{global_lambda}
\end{equation}
 for some positive constant $C_{\Omega,T}$ which depend on $\Vert c\Vert_{C^1(\overline{\Omega})}$ and $\Vert d\Vert_{L^{\infty}(Q_T)}$. In what follows we denote $\Vert \lambda\Vert_{\Lambda}:=\sqrt{\langle \lambda,\lambda\rangle_{\Lambda}}$.

Then, for every parameter $\alpha \in (0,1)$, we consider the following mixed formulation:
\begin{equation}\label{mf_stab}
\left\{
\begin{array}{rcll} 
a_{r,\alpha}(y, \overline{y}) + b_{\alpha}(\overline{y}, \lambda) & = & l_{1, \alpha}(\overline{y}), & \qquad \forall\ \overline{y} \in Z \\
b_\alpha(y, \overline{\lambda}) - c_\alpha(\lambda, \overline\lambda) & = &  l_{2, \alpha}(\overline \lambda), &\qquad \forall\ \overline \lambda \in \Lambda,
\end{array}
\right.
\end{equation}
where 
\begin{align}
\label{aralpha}
a_{r,\alpha} : Z \times Z \to \mathbb{R}, & \quad a_{r, \alpha}(y, \overline y) = (1 - \alpha) \jjntGT c^2(x) \partial_{\nu}y\,\partial_{\nu}\overline{y} \, d\sigma dt + r\jjntQT  Ly \, L\overline{y} dxdt \\
\label{balpha}
b_\alpha : Z \times \Lambda \to \mathbb{R}, & \quad b_\alpha(y,\lambda) = \jjntQT Ly \lambda dx dt - \alpha \jjntGT c^2(x) \partial_\nu y \lambda d\sigma dt \\
\label{calpha}
c_\alpha : \Lambda \times \Lambda \to \mathbb{R}, &\quad c_\alpha(\lambda, \overline \lambda) = \alpha \int_0^T \langle L\lambda(t), L \overline \lambda(t) \rangle_{H^{-1}(\Omega)} dt + \alpha \jjntGT c^2(x) \lambda \overline \lambda d\sigma dt \\
\label{l1alpha}
l_{1, \alpha} : Y \to \mathbb{R}, & \quad l_{1, \alpha}(y) = (1 - \alpha) \jjntGT c^2(x)\partial_\nu y y_{\nu, obs} d\sigma dt \\
\label{L2alpha}
l_{2, \alpha} : \Lambda \to \mathbb{R}, & \quad l_{2, \alpha}(\lambda) = -\alpha \jjntGT c^2(x) y_{\nu, obs} \lambda d\sigma dt.
\end{align}

From the symmetry of $a_{r,\alpha}$ and $c_{\alpha}$, we easily check that this formulation corresponds to the saddle point problem : 
\begin{equation}
\label{Lralpha}
\left\{
\begin{aligned}
& \sup_{\lambda\in \Lambda}\inf_{y\in Z} \mathcal{L}_{r,\alpha}(y,\lambda), \\
&\mathcal{L}_{r,\alpha}(y,\lambda):=\mathcal{L}_r(y,\lambda)-\frac{\alpha}{2}\Vert L\lambda\Vert^2_{L^2(H^{-1}(\Omega))} -\frac{\alpha}{2}\Vert c(x)\left(\lambda-(\partial_{\nu} y-y_{\nu,obs})\right)\Vert_{L^2(\Gamma_T)}^2.
\end{aligned}
\right.
\end{equation}
Precisely, the following holds true. 
\begin{proposition}\label{pr:mfs}
Under the hypothesis (\ref{iobs}), for every $\alpha \in (0,1)$ the stabilized mixed formulation (\ref{mf_stab}) is well-posed. Moreover, the unique pair $(y, \lambda) \in Z \times \Lambda$ satisfies
\begin{equation}\label{eq:estheta}
\theta \|y\|_{Z}^2 + \alpha \|\lambda\|^2_\Lambda \le \frac{(1 - \alpha)^2 + \alpha \theta}{\theta}  \|y_{\nu, obs}\|^2_{L^2(\Gamma_T)}
\end{equation}
with $\theta := \min \left( 1 - \alpha, r/\eta \right)$.
\end{proposition}
\textsc{Proof-} 
We easily get the continuity of the bilinear form $a_{r,\alpha}, b_{\alpha}$ and $c_{\alpha}$:
\begin{equation}
\nonumber
\begin{aligned}
&\vert a_{r, \alpha}(y,\overline y)| \le \max \left( 1-\alpha, \frac{r}{\eta}\right) \|y\|_Z \|\overline y\|_Z, \quad \forall y, \ \overline y \in Z,\\
& | b_\alpha(y, \alpha) |  \leq \frac{C_{\Omega, T} + \alpha}{\eta} \|y\|_Z  \|\lambda\|_\Lambda, \quad \forall y\in Z, \quad \forall \lambda\in \Lambda, \\
& c_\alpha(\lambda, \overline\lambda) \le \alpha \|\lambda\|_\Lambda \|\overline\lambda\|_\Lambda, \quad \forall \lambda, \overline \lambda \in \Lambda
\end{aligned}
\end{equation}
and of the linear form $l_{1,\alpha}$ and $l_{2,\alpha}$:  $\|l_{1,\alpha}\|_{Z'} \leq (1 - \alpha) \|y_{\nu, obs}\|_{L^2(\Gamma_T)}$ and $\|l_{2,\alpha}\|_{\Lambda'} \leq \alpha \|y_{\nu, obs}\|_{L^2(\Gamma_T)}$. 

Moreover, with  $\alpha\in (0,1)$, we also obtain the coercivity of $a_{r,\alpha}$ and of $c_{\alpha}$: precisely, we check that 
$a_{r,\alpha}(y,y) \geq \theta \Vert y\Vert^2_Z$ for all $y\in Z$, while $c_\alpha(\lambda, \lambda) \geq \alpha \|\lambda\|_{\Lambda}^2$ for all $\lambda \in \Lambda$.

The result \cite[Prop 4.3.1]{brezzi_new} then implies the well-posedness of the mixed formulation (\ref{mf_stab}) and the estimate (\ref{eq:estheta}).
\Fin

The $\alpha$-term in $\mathcal{L}_{r,\alpha}$ is a stabilization term: it ensures a coercivity property of $\mathcal{L}_{r,\alpha}$ with respect to the variable $\lambda$ and automatically the well-posedness. In particular, there is no need to prove any inf-sup property for the application $b_{\alpha}$.

\begin{proposition}
If the solution $(y,\lambda)\in Z\times L^2(Q_T)$ of (\ref{eq:mf}) verifies $\lambda\in \Lambda$, then the solutions of (\ref{eq:mf}) and (\ref{mf_stab}) coincide. 
\end{proposition}
\textsc{Proof-} The hypothesis of regularity and the relation (\ref{system_lambda}) imply that the solution $(y,\lambda)\in Z\times L^2(Q_T)$ of (\ref{eq:mf}) is also a solution of (\ref{mf_stab}). The result then follows from the uniqueness of the two formulations. \Fin

\begin{remark}
Remark that the following unique continuation type property for (\ref{eq:wave})
$$
\biggl( Ly=0\,\, \textrm{in}\,\, Q_T, \,\,\, \partial_{\nu}y=0\,\,\textrm{on}\,\, \Sigma_T  \biggr)  \Longrightarrow y=0 \,\, \textrm{in}\,\, Q_T,
$$
which is weaker than (\ref{iobs}), suffices to prove that the bilinear form $ \langle \cdot,\cdot\rangle_Z$ in (\ref{eq:pseta}) is a scalar product and therefore the well-posedness of (\ref{eq:mf}). We shall use specifically the observability inequality (\ref{iobs}) in the numerical Section \ref{sec_numer} to get an estimate of $\Vert y-y_h\Vert_{L^2(Q_T)}$ from estimate of $\Vert L(y-y_h)\Vert_{L^2(Q_T)}$  and $\Vert \partial_{\nu}(y-y_h)\Vert_{L^2(\Sigma_T)}$.
\end{remark}

\subsection{Dual formulation of the extremal problem (\ref{eq:mf}) and remarks} \label{sec2_dual}

As discussed at length in \cite{NC-AM-mixedwave} for which we refer to detail, we may also associate to the extremal problem (\ref{P}) an equivalent problem involving only the variable $\lambda$. Again, this is particularly interesting at the numerical level. This requires a strictly positive augmentation parameter $r$.

For any $r>0$, let us define the linear operator $\mathcal{P}_r$ from $L^2(Q_T)$ into $L^2(Q_T)$ by 
\[
\mathcal{P}_r\lambda:= Ly, \quad \forall \lambda\in L^2(Q_T)
\]
where $y \in Z$ is the unique solution to
\begin{equation}\label{eq:imageA}
a_r(y, \overline y) = b(\overline y, \lambda), \quad \forall \overline y \in Z.   
\end{equation}
The assumption $r>0$ is necessary here in order to guarantee the well-posedness of (\ref{eq:imageA}). Precisely, for any $r>0$, the form $a_r$ defines a norm equivalent to the norm on $Z$.

We then have the following two results, proved in \cite[Section 2.2]{NC-AM-mixedwave} in the very close context of the exact controllability for (\ref{eq:wave}) (see also Remark \ref{analogie_control} below).
\begin{lemma}\label{propA}
For any $r>0$, the operator $\mathcal{P}_r$ is a strongly elliptic, symmetric isomorphism from $L^2(Q_T)$ into $L^2(Q_T)$.
\end{lemma}
\begin{proposition}\label{prop_equiv_dual}
For any $r>0$, let $y_0\in Z$ be the unique solution of 
\[
a_r(y_0,\overline{y})= l(\overline{y}), \quad \forall \overline{y}\in Z
\]
and let $J_r^{\star\star}:L^2(Q_T)\to \mathbb{R}$ be the functional defined by 
\[
J_r^{\star\star}(\lambda) = \frac{1}{2} \jjntQT \mathcal{P}_r \lambda \, \lambda dx\,dt - b(y_0, \lambda).
\]
The following equality holds : 
\begin{equation}\nonumber
\sup_{\lambda\in L^2(Q_T)}\inf_{y\in Z} \mathcal{L}_r(y,\lambda) = - \inf_{\lambda\in L^2(Q_T)} J_r^{\star\star}(\lambda)\quad + \mathcal{L}_r(y_0,0).
\end{equation}
\end{proposition}
This proposition reduces the search of $y$, solution of problem (\ref{P}), to the minimization of $J_r^{\star\star}$. The well-posedness is a consequence of the ellipticity of the operator $\mathcal{P}_r$ stated in Lemma \ref{propA}.

\begin{remark}\label{analogie_control}
We emphasize that the mixed formulation (\ref{eq:mf}) has a structure very close to the one we get when we address - using the same approach - the null controllability of (\ref{eq:wave}): more precisely, the control of minimal $L^2(\Gamma_T)$-norm which drives to rest the initial data $(y_0,y_1)\in L^2(\Omega)\times H^{-1}(\Omega)$ is given by $v=\partial_{\nu}\ph\, 1_{\Gamma_T}$ where $(\ph,\lambda)\in \Phi\times L^2(Q_T)$ solves the mixed formulation
\begin{equation} 
\nonumber 
\left\{
\begin{array}{rcll}
\noalign{\smallskip} a(\ph, \overline{\ph}) + b(\overline{\ph}, \lambda) & = & \hat{l}(\overline{\ph}), & \qquad \forall \overline{\ph} \in \Phi \\
\noalign{\smallskip} b(\ph, \overline{\lambda}) & = & 0, & \qquad \forall \overline{\lambda} \in L^2(Q_T),
\end{array}
\right.
\end{equation}
where $\Phi:=Z$, $a$ and $b$ are given by (\ref{eq:a}) and (\ref{eq:b}) respectively, while $\hat{l}$ is given by   
\begin{align}
\hat{l}: \Phi \to \mathbb{R},  \quad l(\ph) = -\int_\Omega \ph_t(x, 0) y_0(x) dx + \langle  \ph(\cdot,0),y_1\rangle_{H^1_0(\Omega),H^{-1}(\Omega)}, \nonumber
\end{align}
and depends on the initial data $(y_0,y_1)$ to be controlled. We refer to \cite{NC-AM-mixedwave} for the one dimensional case. 
\end{remark}

\begin{remark}
Reversing the order of "priority" between the constraint  $\partial_{\nu}y-y_{\nu,obs}=0$ in $L^2(\Gamma_T)$ and $Ly-f=0$ in $L^2(Q_T)$, a possibility could be to minimize the functional 
$y\mapsto \Vert Ly-f\Vert^2_{L^2(Q_T)}$ over $y\in Z$ subject to the constraint $\partial_{\nu}y-y_{\nu,obs}=0$ in $L^2(\Gamma_T)$ via the introduction of a Lagrange multiplier in $L^2(\Gamma_T)$. 
However, the fact that the following inf-sup property: there exists $\delta>0$ such that 
$$
\inf_{\lambda\in L^2(\Gamma_T)} \sup_{y\in Z} \frac{\int\!\!\!\int_{\Gamma_T} \lambda y \,d\sigma dt}{\Vert \lambda\Vert_{L^2(\Gamma_T)} \Vert y\Vert_Z}   \geq \delta
$$
associated to the corresponding mixed formulation holds true is an open issue.
On the other hand, if a $\eps$-term is added as in (\ref{FV}), this property is satisfied (we refer again to the book \cite{Klibanov-book}).
\end{remark}


\section[Recovering the solution and the source]{Recovering the source and the solution from a partial observation: a mixed re-formulation of the problem}\label{recovering_y_mu}

Given a partial observation $y_{\nu,obs}$ of the solution on the subset $\Gamma_T$, we now consider the reconstruction of the full solution as well as the source term $f$. 
The situation is different with respect to the previous section, since without additional assumptions on $f$, the couple $(y,f)$ is not unique. Consider the case of a source $f$ supported in a set which is close to $\Omega\times \{T\}$: from the finite propagation of the solution, the source $f$ will not affect the solution $y$ on $\Gamma_T$.

Consequently, we assume that $f(x,t)=\sigma(t)\mu(x)$ with $\sigma\in C^1([0,T])$ such that $\sigma(0)\neq 0$ and $\mu\in H^{-1}(\Omega)$ and then recall the following result of Yamamoto-Zhang (see \cite[Theorem 2.1]{yamamoto}): 

\begin{theorem}\cite{yamamoto} \label{yamamoto_th}
Let us assume that the triplet $(\Gamma_T,T,Q_T)$ satisfies the geometric optic condition. Let $y=y(\mu) \in C([0,T];H_0^1(\Omega))\cap C^1([0,T]; L^2(\Omega))$ be the weak solution of (\ref{eq:wave}) with $c:=1$ and $(y_0,y_1)=(0,0)$. Then, there exists a positive constant $C$ such that 
\begin{equation}
\nonumber
C^{-1}\Vert \mu\Vert_{H^{-1}(\Omega)} \leq  \Vert c(x)\,\partial_{\nu} y\Vert_{L^2(\Gamma_T)} \leq C \Vert \mu\Vert_{H^{-1}(\Omega)}, \quad \forall\mu\in H^{-1}(\Omega).
 \end{equation}
\end{theorem} 
Therefore, assuming that the initial condition vanishes, this stability result implies the uniqueness of the source. 

\
We then consider the following extremal problem : 
\begin{equation}
\label{Pmu}
\tag{$\mathcal{P}_{y,\mu}$}
\left\{
\begin{aligned}
& \inf  J(y,\mu):=  \frac{1}{2}\Vert c(x)(\partial_{\nu} y-y_{\nu,obs})\Vert^2_{L^2(\Gamma_T)}, \\
& \textrm{subject to}\quad (y,\mu)\in W
\end{aligned}
\right.
\end{equation}
where $W$ is the space defined by 
\begin{equation}
\label{defW}
\begin{aligned}
W:=\biggl\{(y,\mu); \ & y \in C([0,T];H_0^1(\Omega))\cap C^1([0,T]; L^2(\Omega)), \  \mu\in H^{-1}(\Omega), \\
& Ly-\sigma \mu=0 \, \textrm{in } Q_T, \ y(\cdot,0)=y_t(\cdot,0)=0\biggr\}.
\end{aligned}
\end{equation}
Note that, from Theorem \ref{yamamoto_th},  $\partial_{\nu} y\in L^2(\Gamma_T)$ for any $(y,\mu)\in W$. We easily check that $W$ is a Hilbert space endowed with the norm $\Vert (y,\mu)\Vert_W:=\Vert c(x) \partial_{\nu}y\Vert_{L^2(\Gamma_T)}$. 
The extremal problem $(\mathcal{P}_{y,\mu})$ is well posed : the functional $J$ is continuous over $W$, is strictly convex and is such that $J(y,\mu)\to +\infty$ as $\Vert (y,\mu)\Vert_W\to \infty$. 
Moreover, in view of Theorem \ref{yamamoto_th}, the solution $\mu$ is uniformly bounded in $H^{-1}$-norm.

As in the previous section, in order to avoid the resolution of the extremal problem (\ref{Pmu}) by an iterative process, we introduce a mixed formulation taking the equation $Ly-\sigma \mu=0$
as a constraint equality in $L^2(Q_T)$. In this respect, we define the following non empty vector space: 

\begin{equation}
\nonumber
\begin{aligned}
Y:=\biggl\{(y,\mu); \ & y \in C([0,T];H_0^1(\Omega))\cap C^1([0,T]; L^2(\Omega)),\  \mu\in H^{-1}(\Omega), \\
& Ly-\sigma \mu\in L^2(Q_T),\ y(\cdot,0)=y_t(\cdot,0)=0\biggr\}.
\end{aligned}
\end{equation}

Then, as in Section \ref{recovering_y}, we introduce the following hypothesis: 

\begin{hyp}
There exists a positive constant $C_{obs}=C(\Gamma_T,T,\Vert c\Vert_{C^1(\overline{\Omega})},\Vert d\Vert_{L^{\infty}(\Omega)})$ such that the following estimate holds: 
\begin{equation}
\Vert \mu\Vert^2_{H^{-1}(\Omega)} \leq C_{obs} \biggl( \Vert c(x)\partial_{\nu} y\Vert^2_{L^2(\Gamma_T)}+ \Vert Ly-\sigma \mu\Vert^2_{L^2(Q_T)}  \biggr), \quad \forall (y,\mu)\in Y.  \label{iobs2}\tag{$\mathcal{H}_2$}
\end{equation}
\end{hyp}

Again, in the case where the velocity $c$ is constant, this inequality is a consequence of Theorem \ref{yamamoto_th}: it suffices to decompose any $y$ with $(y,\mu)\in Y$ as $y=y_1+y_2$
where  $y_1$ and $y_2$ solve  
\begin{equation}
\nonumber
Ly_1=\sigma \mu \quad \textrm{in}\quad Q_T, \qquad y_1=0 \quad \textrm{on} \quad \Sigma_T, \qquad (y_1(\cdot,0),y_{1,t}(\cdot,0))=(0,0), 
\end{equation}
and 
\begin{equation}
\nonumber
Ly_2=Ly-\sigma \mu \quad \textrm{in}\quad Q_T, \qquad y_2=0 \quad \textrm{on} \quad \Sigma_T, \qquad (y_2(\cdot,0),y_{2,t}(\cdot,0))=(0,0), 
\end{equation}
respectively and then, applying Theorem \ref{yamamoto_th} for $y_1$ and \cite[Theorem 4.1]{JLL88} for $y_2$, to write 
\begin{equation}
\nonumber
\begin{aligned}
C^{-1}\Vert \mu\Vert_{H^{-1}(\Omega)}\leq & \Vert c(x)\partial_{\nu} y_1\Vert_{L^2(\Gamma_T)} \leq \Vert c(x)\partial_{\nu} y\Vert_{L^2(\Gamma_T)}+ \Vert c(x)\partial_{\nu} y_2\Vert_{L^2(\Gamma_T)}\\
& \leq \Vert c(x)\partial_{\nu} y\Vert_{L^2(\Gamma_T)} + C(\Omega,T) \Vert Ly-\sigma \mu\Vert_{L^1(0,T;L^2(\Omega))}
\end{aligned}
\end{equation}
for some $C(\Omega,T)>0$.

Then, within this hypothesis, for any $\eta>0$  as in Section \ref{recovering_y}, we define on $Y$ the bilinear form 
\begin{equation}
\label{eq:pseta_mu}
\langle (y,\mu),(\overline{y},\overline{\mu})\rangle_Y:= \jjntGT c^2(x)\,\partial_{\nu} y\,\partial_{\nu}\overline{y}\,d\sigma dt + \eta \jjntQT (Ly-\sigma \mu)\, (L\overline{y}-\sigma \overline{\mu})  \,dxdt,
\end{equation}
for every $(y,\mu), \ (\overline y, \overline \mu) \in Y$.
We note the corresponding semi-norm $\Vert (y,\mu)\Vert_Y:=\sqrt{\langle(y,\mu),\ (y,\mu)\rangle_Y}$. Again, in view of (\ref{dependancecontinu}), $\partial_{\nu}y\in L^2(\Gamma_T)$ for any $(y,\mu)\in Y$.
We have the following result: 
\begin{lemma}
Under the hypotheses $(\mathcal{H}_2)$, the space $Y$ is a Hilbert space with the scalar product $\langle \cdot,\ \cdot \rangle_Y$ defined by (\ref{eq:pseta_mu}).
\end{lemma}
\textsc{Proof-} From the inequality $(\mathcal{H}_2)$, the semi-norm $\| \cdot \|_Y$ is indeed a norm. Let us check that $Y$ is closed with respect to this norm.  Let us consider a sequence 
$\{(y_k,\mu_k)\}_{k\geq 1}\in Y$ such that $(y_k,\mu_k)\to (y,\mu)$ for the norm $Y$. Then, there exists $f\in L^2(\Omega)$ such that $L y_k-\sigma \mu_k\to f$ in $L^2(Q_T)$ and $g\in L^2(\Gamma_T)$ such that $\partial_{\nu}y_k\to g \in L^2(\Gamma_T)$. Then, (\ref{iobs2}) implies that $\mu_k\to \mu$ in $H^{-1}(\Omega)$.    Consequently, 
$Ly_k=(Ly_k -\sigma \mu_k) +\sigma\mu_k$ converges to $f+\sigma \mu$ in $L^2(0,T; H^{-1}(\Omega))$. Consequently, $y_k$ can be considered as a sequence of solutions of the hyperbolic equation with zero initial data and second hand term converging in $L^2(0,T,H^{-1}(\Omega))$. Therefore, by the continuous dependence of the solution, $y_k\to y$  in $C([0,T];L^2(\Omega))\cap C^1([0,T];H^{-1}(\Omega))$ where $y$ is the solution of $Ly= f+\sigma \mu\in L^2(Q_T)$ with $(y(\cdot,0),y_t(\cdot,0))=(0,0) \in \boldsymbol{V}$. But, again, from Theorem \ref{yamamoto_th}, this implies that $\partial_{\nu}y\in L^2(\Gamma_T)$ (and therefore $g=\partial_{\nu}y$ on $\Gamma_T$) and $y$ enjoys the regularity $C([0,T];H^1_0(\Omega))\cap C^1([0,T];L^2(\Omega))$
as the sum of two solutions (as for $y_1$ and $y_2$ above) in this space. Therefore, $(y,\mu)\in Y$. \Fin
\
\par\noindent

Proceeding as in Section \ref{sec2_direct}, we introduce a Lagrange multiplier $\lambda\in L^2(Q_T)$ and the following mixed formulation: find $((y,\mu), \lambda)\in Y\times L^2(Q_T)$ solution of 
\begin{equation} \label{eq:mfeps}
\left\{
\begin{array}{rcll}
\noalign{\smallskip} a((y,\mu), (\overline{y},\overline{\mu})) + b((\overline{y},\overline{\mu}), \lambda) & = & l(\overline{y},\overline{\mu}), & \qquad \forall (\overline{y},\overline{\mu}) \in Y \\
\noalign{\smallskip} b((y,\mu), \overline{\lambda}) & = & 0, & \qquad \forall \overline{\lambda} \in L^2(Q_T),
\end{array}
\right.
\end{equation}
where
\begin{align}
\nonumber
 & a: Y \times Y \to \mathbb{R},  \quad a((y,\mu), (\overline{y},\overline{\mu})) := \jjntGT c^2(x) \partial_{\nu} y\partial_{\nu}\overline{y}\,d\sigma dt,
\\
\nonumber
& b: Y \times L^2(Q_T)  \to \mathbb{R},  \quad b((y,\mu),\lambda) := \jjntQT \lambda (Ly-\sigma \mu) dx\, dt,\\
\nonumber
& l: Y \to \mathbb{R},  \quad l(y,\mu) := \jjntGT c^2(x)\,\partial_{\nu} y \,y_{\nu,obs}\, d\sigma dt. 
\end{align}
\begin{theorem}\label{th:mfeps} Under the hypothesis (\ref{iobs2}), the following holds :
\begin{enumerate}
\item The mixed formulation (\ref{eq:mfeps}) is well-posed.
\item The unique solution $((y,\mu), \lambda) \in Y\times L^2(Q_T)$ is the saddle-point of the Lagrangian $\mathcal{L}:Y\times L^2(Q_T)\to \mathbb{R}$ defined by
\[
\mathcal{L}((y,\mu), \lambda):= \frac{1}{2}a((y,\mu),(y,\mu)) + b((y,\mu),\lambda)- l(y,\mu).
\]
Moreover, the pair $(y,\mu)$ solves the extremal problem (\ref{Pmu}).

\item The following estimates hold :  
\begin{equation}
\Vert (y,\mu) \Vert_Y= \Vert c(x)\partial_{\nu}y\Vert_{L^2(\Gamma_T)}\leq \Vert c(x)\,y_{\nu,obs}\Vert_{L^2(\Gamma_T)}  \label{ineqy_th1}
\end{equation}
and 
\begin{equation}
\Vert \lambda\Vert_{L^2(Q_T)} \leq 2\sqrt{C_{\Omega,T}+\eta}  \Vert c(x)\, y_{\nu,obs}\Vert_{L^2(\Gamma_T)}  \label{estimate_lambdaeps}
\end{equation}
for some constant $C_{\Omega,T}>0$.
\end{enumerate}
\end{theorem}
\par\noindent
The proof is very close to the proof of Theorem \ref{th:mf}.  
In particular, the inf-sup property can be obtained by taking, for any $\lambda\in L^2(Q_T)$, $\mu_0=0$ and $y_0$ as in (\ref{infsup_y0}) leading to $(y_0,\mu_0)\in Y$ so that the inf-sup constant
\begin{equation}
\nonumber
\delta:=\inf_{\lambda\in L^2(Q_T)} \sup_{(y,\mu)\in Y}  \frac{b((y,\mu),\lambda)}{\Vert (y,\mu) \Vert_Y \Vert \lambda\Vert_{L^2(Q_T)}}
\end{equation}
is bounded by above by $(C_{\Omega,T}+\eta)^{-1/2}$.

Then, as a consequence of the estimate (\ref{iobs2}), the source term $\mu$ is uniformly bounded in the $H^{-1}$-norm:
\begin{corollary}
Assuming (\ref{iobs2}), the solution $\mu$ of (\ref{eq:mfeps}) satisfies $\Vert \mu\Vert_{H^{-1}(\Omega)}\leq C_{obs} \Vert c(x)\,y_{\nu,obs}\Vert_{L^2(\Gamma_T)}$. 
\end{corollary}


 Again, it is very convenient to "augment" the Lagrangian (see \cite{fortinglowinski}) and consider instead the Lagrangian $\mathcal{L}_r$ defined for any $r>0$ by 
\begin{equation}
\begin{aligned}
& \mathcal{L}_r((y,\mu),\lambda):=\frac{1}{2}a_r((y,\mu),(y,\mu))+b(y,\lambda)-l(y,\mu), \\
& a_r((y,\mu),(y,\mu)):=a((y,\mu),(y,\mu))+r\Vert Ly-\sigma \mu\Vert^2_{L^2(Q_T)}. 
\end{aligned}
 \nonumber
\end{equation}
Since $a((y,\mu),(y,\mu))=a_r((y,\mu),(y,\mu))$ in $W$, the Lagrangian $\mathcal{L}$ and $\mathcal{L}_r$ share the same saddle-point. The positive number $r$ is an augmentation parameter. Moreover, proceeding as in Section \ref{sec2_dual}, we may also associate to the saddle-point problem $\sup_{\lambda\in L^2(Q_T)}\inf_{(y,\mu)\in Y} \mathcal{L}_r((y,\mu),\lambda)$ a dual problem, which again reduces the search of the couple $(y,\mu)$, solution of problem (\ref{Pmu}), to the minimization of a elliptic functional with respect to $\lambda$. 

\begin{proposition}\label{prop_equiv_dual_eps}
For any $r>0$, let $(y_0,\mu_0)\in Y$ be the unique solution of 
\[
a_r((y_0,\mu_0),(\overline{y},\overline{\mu}))= l(\overline{y},\overline{\mu}), \quad \forall (\overline{y},\overline{\mu})\in Y
\]
and let $\mathcal{P}_r$ be the operator from $L^2(Q_T)$ into $L^2(Q_T)$ defined by 
$\mathcal{P}_r\lambda:= Ly-\sigma \mu$ where $(y,\mu)\in Y$ is the unique solution to 
\begin{equation}\label{eq:imageAeps}
a_r((y,\mu), (\overline y,\overline \mu)) = b((\overline y,\overline \mu), \lambda), \quad \forall (\overline y,\overline \mu) \in Y.   
\end{equation}
The operator $\mathcal{P}_r$ is strongly elliptic and symmetric. Moreover,  the following equality holds
\begin{equation}\nonumber
\sup_{\lambda\in L^2(Q_T)}\inf_{(y,\mu)\in Y} \mathcal{L}_r((y,\mu),\lambda) = - \inf_{\lambda\in L^2(Q_T)} J_r^{\star\star}(\lambda)\quad + \mathcal{L}_r((y_0,\mu_0),0)
\end{equation}
where $J_r^{\star\star}:L^2(Q_T)\to \mathbb{R}$ is the functional defined by 
\[
J_r^{\star\star}(\lambda) = \frac{1}{2} \jjntQT (\mathcal{P}_r \lambda) \lambda dx\,dt - b((y_0,\mu_0), \lambda).
\]
\end{proposition}

\textsc{Proof-} From the definition of $a_r$, we easily get that $\Vert \mathcal{P}_r\lambda\Vert_{L^2(Q_T)}\leq r^{-1} \Vert \lambda\Vert_{L^2(Q_T)}$ and the continuity of $\mathcal{P}_r$. Next, consider any $\lambda^{\prime}\in L^2(Q_T)$ and denote by $(y^{\prime},\mu^{\prime})$ the corresponding unique solution of (\ref{eq:imageAeps}) so that $\mathcal{P}_r\lambda^{\prime}:=Ly^{\prime}-\sigma\mu^{\prime}$. Relation (\ref{eq:imageAeps}) with $(\overline{y},\overline{\mu})=(y^{\prime},\mu^{\prime})$ then implies that 
\begin{equation}
\jjntQT  (\mathcal{P}_r\lambda^{\prime})\lambda \,dx\,dt = a_r((y,\mu),(y^{\prime},\mu^{\prime}))    \label{arAlambda}
\end{equation}
and therefore the symmetry and positivity of $\mathcal{P}_r$. The last relation with $\lambda^{\prime}=\lambda$ and the observability estimate (\ref{iobs2}) imply that $\mathcal{P}_r$ is also positive definite.

Finally, let us check the strong ellipticity of $\mathcal{P}_r$, equivalently that the bilinear functional $(\lambda,\lambda^{\prime})\to \jjntQT (\mathcal{P}_r\lambda)\lambda^{\prime}\,dx\,dt$ is $L^2(Q_T)$-elliptic. Thus we want to show that 
\begin{equation}\label{ellipticity_A}
\jjntQT (\mathcal{P}_r\lambda)\lambda\,dx\,dt \geq C \Vert \lambda\Vert^2_{L^2(Q_T)}, \quad\forall \lambda\in L^2(Q_T)
\end{equation}
for some positive constant $C$. Suppose that (\ref{ellipticity_A}) does not hold; there exists then a sequence $\{\lambda_n\}_{n\geq 0}$ of $L^2(Q_T)$ such that 

\[
\Vert \lambda_n\Vert_{L^2(Q_T)}=1, \quad\forall n\geq 0, \qquad \lim_{n\to\infty} \jjntQT (\mathcal{P}_r\lambda_n)\lambda_n\,dx\,dt=0.
\]
Let us denote by $(y_n,\mu_n)$ the solution of (\ref{eq:imageAeps}) corresponding to $\lambda_n$. From (\ref{arAlambda}), we then obtain that 
\begin{equation}
\lim_{n\to\infty} (r \Vert Ly_n-\sigma \mu_n\Vert^2_{L^2(Q_T)}+ \Vert c(x)\,\partial_{\nu}y_n\Vert^2_{L^2(\Gamma_T)}) =0 .\label{limit}
\end{equation}

From (\ref{eq:imageAeps}) with $(y,\mu)=(y_n,\mu_n)$ and $\lambda=\lambda_n$, we have
\begin{equation} \label{phinlambdan}
\jjntQT  [ r (L y_n-\sigma \mu_n)-\lambda_n ] (L\overline{y}-\sigma\overline{\mu}) \, dx\,dt + \jjntGT c^2(x) \partial_{\nu}y_n \partial_{\nu}\overline{y} d\sigma\,dt =0, \quad \forall (\overline{y},\overline{\mu})\in Y.
\end{equation}
We define the sequence $\{\overline{y}_n\}_{n\geq 0}$ as follows : 
\begin{equation}\nonumber
\left\{
\begin{aligned}
& L\overline{y}_n =  r\,(Ly_n-\sigma \mu_n)-\lambda_n, & &\textrm{in}\quad Q_T, \\
& \overline{y}_n=0, & &\textrm{on}\quad \Sigma_T,\\
& \overline{y}_n(\cdot,0)=\overline{y}_{n,t}(\cdot,0)=0, & &\textrm{in}\quad \Omega,
\end{aligned}
\right.
\end{equation}
so that, for all $n$, $\overline{y}_n$ is the solution of a hyperbolic equation with zero initial data and source term $r\,(Ly_n-\sigma\mu_n)-\lambda_n$ in $L^2(Q_T)$. Using again \cite[Theorem 4.1]{JLL88}, we get $\Vert \partial_{\nu}\overline{y}_n\Vert_{L^2(\Gamma_T)} \leq  \sqrt{C_{\Omega,T}} \Vert r(Ly_n-\sigma \mu_n)-\lambda_n\Vert_{L^2(Q_T)}$, so that $(\overline{y}_n,0)\in Y$. Then, using (\ref{phinlambdan}) with $(\overline{y},\overline{\mu})=(\overline{y}_n,0) \in Y$ we get 
\[
\Vert r(Ly_n-\sigma\mu_n)-\lambda_n\Vert_{L^2(Q_T)} \leq \sqrt{C_{\Omega,T}} \Vert c(x)\partial_{\nu}y_n\Vert_{L^2(\Gamma_T)}.
\]
Then, from (\ref{limit}), we conclude that $\lim_{n\to +\infty} \Vert \lambda_n\Vert_{L^2(Q_T)}=0$ leading to a contradiction and to the strong ellipticity of the operator $\mathcal{P}_r$.
The rest of the proof is standard. \Fin

To end this section, it is worth to mention that a result similar to Theorem \ref{yamamoto_th} is given in \cite[Theorem 1]{yamamoto95} with a control of $\mu$ in $L^2(\Omega)$. However, this requires an observation of $\partial_{\nu} y$ in $H^1(0,T;L^2(\Gamma))$ (equivalently of $\partial_t\partial_\nu y$ in $L^2(\Gamma_T)$) which may be too strong in real life applications.  


\section{Numerical Analysis of the mixed formulations} \label{sec_numer}

\subsection{Numerical approximation of the mixed formulation (\ref{eq:mf})}

We now proceed to the numerical analysis of the mixed formulation (\ref{eq:mf}), assuming $r>0$. We follow \cite{NC-AM-mixedwave}, to which we refer for the details. 

Let $Z_h$ and $\Lambda_h$ be two finite dimensional spaces parametrized by the variable $h$ such that $Z_h\subset Z, \Lambda_h\subset L^2(Q_T)$ for every $h>0$.
Then, we can introduce the following approximated problems: find $(y_h,\lambda_h)\in Z_h\times \Lambda_h$ solution of 
\begin{equation} \label{eq:mfh}
\left\{
\begin{array}{rcll}
\noalign{\smallskip} a_r(y_h,\overline{y}_h) + b(\overline{y}_h, \lambda_h) & = & l(\overline{y}_h), & \qquad \forall \overline{y}_h\in Z_h \\
\noalign{\smallskip} b(y_h, \overline{\lambda}_h) & = & 0, & \qquad \forall \overline{\lambda}_h \in \Lambda_h.
\end{array}
\right.
\end{equation}

The well-posedness of this mixed formulation is again a consequence of two properties. The first one is the coercivity of the bilinear form $a_r$ on the subset 
$\mathcal{N}_h(b)=\{y_h\in Z_h; b(y_h,\lambda_h)=0\,\, \forall \lambda_h\in \Lambda_h\}$.
 Actually, from the relation $a_r(y,y)\geq (r/\eta)\Vert y\Vert_Z^2$ for all $y\in Z$, the form $a_r$ is coercive on the full space $Z$, and so a fortiori on $\mathcal{N}_h(b)\subset Z_h\subset Z$. 
 
 The second property is a discrete inf-sup condition: precisely, for any $h$, 
\begin{equation}
\delta_h:=\inf_{\lambda_h\in \Lambda_h}\sup_{y_h\in Z_h}  \frac{b(y_h,\lambda_h)}{\Vert \lambda_h\Vert_{L^2(Q_T)}\Vert y_h\Vert_Z}>0. \label{infsupdiscret}
\end{equation}

Let us assume that this condition holds, so that for any fixed $h>0$, there exists a unique couple $(y_h,\lambda_h)$ solution of (\ref{eq:mfh}). We then have the following estimate.

\begin{proposition} \label{prop_estimateh_eps0}
Let $h>0$. Let $(y,\lambda)$ and $(y_h,\lambda_h)$ be the solution of (\ref{eq:mf}) and of (\ref{eq:mfh}) respectively. Let $\delta_h$ be the discrete inf-sup constant defined by (\ref{infsupdiscret}). Then,
\begin{align}
\label{eq:est_e01}
& \Vert y-y_h \Vert_{Z} \leq  2\biggl(1+\frac{1}{\sqrt{\eta} \delta_h}\biggr)d(y,Z_h)+\frac{1}{\sqrt{\eta}} d(\lambda,\Lambda_h), \\
\label{eq:est_e02}
& \Vert \lambda-\lambda_h \Vert_{L^2(Q_T)}  \leq  \biggl(2+\frac{1}{\sqrt{\eta}\delta_h}\biggr)\frac{1}{\delta_h}d(y,Z_h)+
\frac{3}{\sqrt{\eta} \delta_h} d(\lambda,\Lambda_h)
\end{align}
where  
$ d(\lambda,\Lambda_h):=\inf_{\lambda_h\in \Lambda_h} \Vert \lambda-\lambda_h\Vert_{L^2(Q_T)}$
and
\begin{equation}
\label{def_distanceyZh}
\begin{aligned}
d(y,Z_h):=& \inf_{y_h\in Z_h}\Vert y-y_h\Vert_Z\\
=& \inf_{y_h \in Z_h}   \biggl( \Vert \partial_{\nu} y-\partial_{\nu}y_h\Vert^2_{L^2(\Gamma_T)}+  \eta \Vert L(y-y_h) \Vert^2_{L^2(Q_T)} \biggr)^{1/2}.
\end{aligned}
\end{equation}
\end{proposition}
\textsc{Proof-} From the classical theory of approximation of saddle point problems (see \cite[Theorem 5.2.2]{brezzi_new}) we have that
\begin{align}
\nonumber
\Vert y-y_h \Vert_{Z} \leq & \left( \frac{2 \| a_r \|_{(Z\times Z)^{\prime}}}{\alpha_0} + \frac{2 \| a_r \|_{(Z\times Z)^{\prime}}^\frac12 \| b \|_{(Z\times L^2(Q_T))^{\prime}}}{\alpha_0^\frac12 \delta_h} \right) d(y,Z_h) \\
& +\frac{\Vert b\Vert_{(Z\times L^2(Q_T))^{\prime}}}{\alpha_0} d(\lambda,\Lambda_h)\nonumber
\end{align}
and 
\begin{align}\nonumber
\Vert \lambda-\lambda_h \Vert_{L^2(Q_T)} \leq  &
\left( \frac{2\| a_r \|_{(Z\times Z)^{\prime}}^\frac32}{\alpha_0^\frac12 \delta_h} + \frac{\| a_r \|_{(Z\times Z)^{\prime}} \| b \|_{(Z\times L^2(Q_T))^{\prime}}}{\delta_h^2} \right) d(y,Z_h) \\
& + \frac{3 \| a_r \|^\frac12 \| b \|_{(Z\times L^2(Q_T))^{\prime}}}{\alpha_0^\frac12 \delta_h} d(\lambda,\Lambda_h). \nonumber
\end{align}
Since, $\Vert a_r\Vert_{(Z\times Z)^{\prime}}=\alpha_0=1$; $\Vert b\Vert_{(Z\times L^2(Q_T))^{\prime}}= \frac{1}{\sqrt{\eta}}$, the result follows.  \Fin 
\begin{remark}
For $r=0$, the discrete mixed formulation (\ref{eq:mfh}) is not well-posed over $Z_h\times \Lambda_h$ because the form $a_{r=0}$ is not coercive over the discrete kernel of $b$: the equality $b(y_h,\lambda_h)=0$ for all $\lambda_h\in \Lambda_h$ does not imply in general that $L y_h$ vanishes. Therefore, the term $r \Vert Ly_h\Vert^2_{L^2(Q_T)}$, which appears in the Lagrangian $\mathcal{L}_r$, may be understood as a stabilization term: for any $h>0$, it ensures the uniform coercivity of the form $a_r$ and vanishes at the limit in $h$. We also emphasize that this term is not a regularization term as it does not add any regularity on the solution $y_h$. 
\end{remark}


Let $n_h=\dim(Z_h), m_h=\dim(\Lambda_h)$ the dimension of the space $Z_h$ and $\Lambda_h$ respectively. Let the real matrices $A_{r,h}\in \mathbb{R}^{n_h,n_h}$, $B_{h}\in \mathbb{R}^{m_h,n_h}$, $J_h\in \mathbb{R}^{m_h,m_h}$ and $L_h\in \mathbb{R}^{n_h}$ be defined by   
\begin{equation}
\label{def_matrix}
\left\{
\begin{aligned}
& a_r( y_h, \overline{y_h})=   \langle A_{r,h} \{y_h\}, \{\overline{y_h}\} \rangle_{\mathbb{R}^{n_h},\mathbb{R}^{n_h}} & \forall y_h,\overline{y_h}\in Z_h,\\
& b(y_h,\lambda_h)=  \langle B_{h} \{y_h\}, \{\lambda_h\} \rangle_{\mathbb{R}^{m_h},\mathbb{R}^{m_h}} & \forall y_h\in Z_h, \lambda_h\in \Lambda_h,\\
& \jjntQT \lambda_h\overline{\lambda_h}\,dx\,dt=  \langle J_h \{\lambda_h\}, \{\overline{\lambda_h}\} \rangle_{\mathbb{R}^{m_h},\mathbb{R}^{m_h}} & \forall \lambda_h,\overline{\lambda_h}\in \Lambda_h, \\
& l(y_h)= \langle L_h,\{y_h\} \rangle_{\mathbb{R}^{n_h}} & \forall y_h\in Z_h,
\end{aligned}
\right.
\end{equation}
where $\{y_h\}\in \mathbb{R}^{n_h}$ denotes the vector associated to $y_h$ and $\langle \cdot,\cdot \rangle_{\mathbb{R}^{n_h},\mathbb{R}^{n_h}}$ the usual scalar product over $\mathbb{R}^{n_h}$. With these notations, the problem (\ref{eq:mfh}) reads as follows: find $\{y_h\}\in \mathbb{R}^{n_h}$ and $\{\lambda_h\}\in \mathbb{R}^{m_h}$ such that 
\begin{equation} \label{matrixmfh}
\left(
\begin{array}{cc}
A_{r,h} &  B_h^T \\
B_h & 0   
\end{array}
\right)_{\mathbb{R}^{n_h+m_h,n_h+m_h}}
\left(
\begin{array}{c}
\{y_h\}  \\
\{\lambda_h\}    
\end{array}
\right)_{\mathbb{R}^{n_h+m_h}}  =
\left(
\begin{array}{c}
L_h \\
0   
\end{array}
\right)_{\mathbb{R}^{n_h+m_h}}.
\end{equation}
The matrix $A_{r,h}$ as well as the mass matrix $J_h$ are symmetric and positive definite for any $h>0$ and any $r>0$. On the other hand, the full matrix of order $m_h+n_h$ in (\ref{matrixmfh}) is symmetric but not necessarily positive definite. 

We recall  (see \cite[Theorem 3.2.1]{brezzi_new})  that the inf-sup property (\ref{infsupdiscret}) is equivalent to the injective character of the matrice $B_h^T$ of size $n_h\times m_h$, that is $Ker(B_h^T)=0$. If a necessary condition is given by $m_h\leq n_h$, this property strongly depends on the structure of the spaces $Z_h$ and $\Lambda_h$. We will discuss numerically this property in Remark \ref{better_estimate} for a specific choice of approximation.

\subsubsection{$C^1$-finite elements and order of convergence for $N=1$}

The finite dimensional and conformal space $Z_h$ must be chosen such that $Ly_h$ belongs to $L^2(Q_T)$ for any $y_h\in Z_h$. This is guaranteed, for instance, as soon as $y_h$ possesses second-order derivatives in $L^2(Q_T)$. As in \cite{NC-AM-mixedwave}, we consider a conformal approximation based on functions continuously differentiable with respect to both variables $x$ and $t$. 

We introduce a triangulation $\mathcal{T}_h$ such that $\overline{Q_T}=\cup_{K\in \mathcal{T}_h} K$ and we assume that $\{\mathcal{T}_h\}_{h>0}$ is a regular family. We note 
$
h:=\max\{\textrm{diam}(K), K\in \mathcal{T}_h\}
$
where $\textrm{diam}(K)$ denotes the diameter of $K$. Then, we introduce the space $Z_h$ as follows : 
\begin{equation}
Z_h:=\{z_h\in C^1(\overline{Q_T}):  z_h\vert_K\in \mathbb{P}(K) \quad \forall K\in \mathcal{T}_h, \,\, z_h=0 \,\,\textrm{on}\,\,\Sigma_T\},   \label{defZh}
\end{equation}
where $\mathbb{P}(K)$ denotes an appropriate space of functions in $x$ and $t$. In this work, we consider two choices, in the one-dimensional setting (for which $\Omega\subset \mathbb{R}$ and $Q_T\subset \mathbb{R}^2$): 
\begin{enumerate}
\item The \textit{Bogner-Fox-Schmit} (BFS for short) $C^1$-element defined for rectangles. It involves $16$ degrees of freedom, namely the values of $y_h, y_{h,x}, y_{h,t}, y_{h,xt}$ on the four vertices of each rectangle $K$. 
Therefore $\mathbb{P}(K)= \mathbb{P}_{3,x}\otimes \mathbb{P}_{3,t}$ where $\mathbb{P}_{r,\xi}$ is by definition the space of polynomial functions of order $r$ in the variable $\xi$. We refer to \cite[ch. II, sec. 9, p. 94]{ciarletfem}.
\item The reduced \textit{Hsieh-Clough-Tocher} (HCT for short) $C^1$-element defined for triangles. This is a so-called composite finite element and involves $9$ degrees of freedom, namely, the values of $y_h,y_{h,x}, y_{h,t}$ on the three vertices of each triangle $K$. We refer to \cite[ch. VII, sec. 46, p. 285]{ciarletfem} and to \cite{bernadouHCT,meyer} where the implementation is discussed.
\end{enumerate}
We also define the finite dimensional space
\[
\Lambda_h:=\{\lambda_h\in C^0(\overline{Q_T}),  \lambda_h\vert_K\in \mathbb{Q}(K) \quad \forall K\in \mathcal{T}_h \}
\]
where $\mathbb{Q}(K)$ denotes the space of affine functions both in $x$ and $t$ on the element $K$.

For any $h>0$, we have $Z_h\subset Z$ and $\Lambda_h\subset L^2(Q_T)$.

We then have the following result: 
\begin{proposition}[BFS element for $N=1$ - Rate of convergence for the norm $Z$]   \label{propeps0}
Let $h>0$, let $k\leq 2 $ be a nonnegative integer. Let $(y,\lambda)$ and $(y_h,\lambda_h)$ be the solution of (\ref{eq:mf}) and (\ref{eq:mfh}) respectively. If the solution $(y,\lambda)$ belongs to $H^{k+2}(Q_T)\times H^k(Q_T)$, then there exist two positives constants 
\[ K_i=K_i(\Vert y\Vert_{H^{k+2}(Q_T)}, \Vert c \Vert_{C^1(\overline{Q_T})},\Vert d \Vert_{L^{\infty}(Q_T)}), \qquad i \in \{ 1,2\},\]
independent of $h$, such that 
\begin{align}
\Vert y-y_h \Vert_Z & \leq  K_1 \biggl(1+\frac{1}{\sqrt{\eta} \delta_h}+ \frac{1}{\sqrt{\eta}} \biggr) h^k, \label{eq:yyh0}\\
\Vert \lambda-\lambda_h \Vert_{L^2(Q_T)}&\leq  K_2\biggl(\biggl(1+\frac{1}{\sqrt{\eta} \delta_h}\biggr)\frac{1}{\delta_h} +\frac{1}{\sqrt{\eta}\delta_h} \biggr) h^k\label{eq:llh0}.
\end{align}
\end{proposition}

\par\noindent
\textsc{Proof -}
 From \cite[ch. III, sec. 17]{ciarletfem}, for any $\lambda\in H^k(Q_T)$, $k\leq 2$, there exists $C_1=C_1(\Vert \lambda\Vert_{H^k(Q_T)})$
such that 
\begin{equation} \label{eq:est0h}
\Vert \lambda-\Pi_{\Lambda_h,\mathcal{T}_h}(\lambda)\Vert_{L^2(Q_T)} \leq  C_1 h^k, \quad \forall h>0
\end{equation}
where $\Pi_{\Lambda_h,\mathcal{T}_h}$ designates the interpolant operator from $L^2(Q_T)$ to $\Lambda_h$ associated to the regular mesh $\mathcal{T}_h$. 
Similarly, for any $y\in H^{k+2}(Q_T)$, there exists $C_2=C_2(\Vert y\Vert_{H^{k+2}(Q_T)})$ such that for every $h > 0$ we have
\begin{equation} \label{eq:est1h}
\Vert y-\Pi_{Z_h,\mathcal{T}_h}(y)\Vert_{L^2(Q_T)} \leq  C_2  h^{k+2}, \quad  \Vert y-\Pi_{Z_h,\mathcal{T}_h}(y)\Vert_{H^2(Q_T)} \leq  C_2  h^k. 
\end{equation}
Then, using that the linear operator $P:H^{3/2}(Q_T)\to L^2(\Sigma_T)$ defined by $Py:=\partial_\nu y$ is continuous, there exists a positive constant $C_P$ such that 
\begin{equation}
\Vert c(x)(\partial_\nu y- \partial_\nu  (\Pi_{Z_h,\mathcal{T}_h}(y)))\Vert_{L^2(\Gamma_T)} \leq \Vert c\Vert_{L^\infty(\Gamma_T)}C_P \Vert y-\Pi_{Z_h,\mathcal{T}_h}(y)\Vert_{H^{3/2}(Q_T)}.  \label{operator_neumann}
\end{equation}
We then observe that 
\begin{equation} \label{eq:est2h}
\Vert Ly-L(\Pi_{Z_h,\mathcal{T}_h}(y))\Vert_{L^2(Q_T)} \leq K(\Vert c\Vert_{C^1(\overline{Q_T})}, \Vert d\Vert_{L^{\infty}(Q_T)})\Vert y-\Pi_{Z_h,\mathcal{T}_h}(y)\Vert_{H^2(Q_T)},
\end{equation}
for some positive constant $K$;  (\ref{def_distanceyZh}) then leads to 
\begin{equation}
\begin{aligned}
d(y,Z_h)&\leq  \left( \|\partial_\nu y-\partial_\nu(\Pi_{Z_h,\mathcal{T}_h}(y))\|^2_{L^2(\Gamma_T)} + \eta \|Ly - L(\Pi_{Z_h,\mathcal{T}_h}(y))\|^2_{L^2(Q_T)}\right)^2 \\
&\leq  \sqrt{\Vert c\Vert_{L^\infty(\Gamma_T)}^2 C_p^2+\eta K^2} C_2 h^{k}
\end{aligned}
\label{eq:ch}
\end{equation}
and then from Proposition \ref{prop_estimateh_eps0}, we get that 
\begin{equation}\nonumber
 \Vert y-y_h \Vert_{Z} \leq  2\biggl(1+\frac{1}{\sqrt{\eta} \delta_h}\biggr)\sqrt{\Vert c\Vert_{L^\infty(\Gamma_T)}^2 C_p^2+\eta K^2} C_2 h^k+\frac{1}{\sqrt{\eta}} C_1 h^k.
\end{equation}
Similarly, 
 \begin{equation*}
\Vert \lambda-\lambda_h \Vert_{L^2(Q_T)}\leq  \biggl(2+\frac{1}{\sqrt{\eta} \delta_h}\biggr)\frac{1}{\delta_h} \sqrt{\Vert c\Vert_{L^\infty(\Gamma_T)}^2 C_p^2+\eta K^2} C_2 h^k+\frac{3}{\sqrt{\eta}\delta_h} C_1h^k.
\end{equation*}
From the last two estimates, we obtain the conclusion of the proposition. 
\Fin

It remains now to deduce the convergence of the approximated solution $y_h$ for the $L^2(Q_T)$ norm. This is done using the observability estimate (\ref{iobs}). Precisely, we write that  $(y-y_h)$ solves the hyperbolic equation
\begin{equation}
\nonumber
\left\{
\begin{aligned}
& L(y-y_h)=-L y_h &&\quad \text{in } Q_T \\
& ((y-y_h), (y-y_h)_t)(\cdot, \ 0) \in \boldsymbol{V} && \\
& y-y_h=0 && \quad \text{on } \Sigma_T.
\end{aligned}
\right.
\end{equation}
The continuous dependence of the solution with respect to the right hand side and the initial data leads 
\begin{equation}\nonumber 
\Vert y-y_h\Vert_{L^2(Q_T)}^2 \leq C_{\Omega,T} (\Vert ((y-y_h)(\cdot,0),(y-y_h)_t(\cdot,0))\Vert^2_{\boldsymbol{V}} + \Vert Ly_h\Vert^2_{L^2(Q_T)} ).
\end{equation}
Combined with (\ref{iobs}) applied for $y-y_h$, this leads to 
\begin{equation} \nonumber
\Vert y-y_{h} \Vert^2_{L^2(Q_T)} \leq C_{\Omega,T}(1+C_{obs})( \Vert y-y_h \Vert^2_{L^2(\Gamma_T)} +  \Vert L y_h\Vert^2_{L^2(Q_T)})  \nonumber
\end{equation}
from which we deduce, in view of the definition of the norm $Y$, that 
\begin{equation}
\Vert y-y_h \Vert_{L^2(Q_T)} \leq C_{\Omega,T}(1+C_{obs})  \max(1,\frac{2}{\sqrt{\eta}}) \Vert y-y_h\Vert_Z.   \label{eq:estyL2}
\end{equation}

Eventually, by coupling (\ref{eq:estyL2}) and Proposition \ref{propeps0}, we obtain the following result :

\begin{theorem}[BFS element for $N=1$ - Rate of convergence for the norm $L^2(Q_T)$]   \label{theps0}
Assume that the hypothesis (\ref{iobs}) holds. Let $h>0$, let $k\leq 2$ be a positive integer. Let $(y,\lambda)$ and $(y_h,\lambda_h)$ be the solution of (\ref{eq:mf}) and (\ref{eq:mfh}) respectively. If the solution $(y,\lambda)$ belongs to $H^{k+2}(Q_T)\times H^k(Q_T)$, then there exists two positives constant $K=K(\Vert y\Vert_{H^{k+2}(Q_T)}, \Vert c \Vert_{C^1(\overline{Q_T})},\Vert d \Vert_{L^{\infty}(Q_T)},C_{\Omega,T},C_{obs})$, independent of $h$, such that 
\begin{equation}
\Vert y-y_h \Vert_{L^2(Q_T)}  \leq  K \max(1,\frac{2}{\sqrt{\eta}})\biggl(1+\frac{1}{\sqrt{\eta} \delta_h}+ \frac{1}{\sqrt{\eta}} \biggr) h^k .   \label{estimate_L2_eps0}
\end{equation}
\end{theorem}

\begin{remark}[About the choice of the augmentation parameter $r$]\label{better_estimate}
Estimate (\ref{estimate_L2_eps0}) is not fully satisfactory as it depends on the constant $\delta_h$. 
In view of the complexity of both the constraint $Ly=0$ and of the structure of the space $Z_h$, the theoretical estimation of the constant $\delta_h$ with respect to $h$ is a difficult issue. However, as discussed at length in  \cite[Section 2.1]{NC-AM-mixedwave}, $\delta_h$ can be evaluated numerically for any $h$,  as the solution of the following generalized eigenvalue problem (taking $\eta=r$, so that $a_r(y,y)$ is exactly $\Vert y \Vert^2_{Z}$):   
\begin{equation}\label{eigenvalue}
\delta_h = \inf\biggl\{\sqrt{\delta}:   B_h A_{r,h}^{-1} B_h^T  \{\lambda_h\} = \delta \,J_h \{\lambda_h\}, \quad \forall\, \{\lambda_h\}\in \mathbb{R}^{m_h}\setminus\{0\}\biggr\}
\end{equation}
where the matrices $A_{r,h}$, $B_h$ and $J_h$ are defined in (\ref{def_matrix}).

The computation of $\delta_h$ given by the problem (\ref{eigenvalue}) with respect to $h$ and to $r=\eta$ has been performed in \cite[Section 4.2]{NC-AM-mixedwave} for the following data : $T=2$, $Q_T=(0,1)\times (0,T)$, $T=2$, $\Gamma_T=\{1\}\times (0,T)$. There, it is observed, that for both the BFS and the HCT finite element on regular meshes (non necessarily uniform), the constant $\delta_h$ behaves like  
\begin{equation}
\delta_h\approx C_r \frac{h}{\sqrt{r}} \quad \textrm{as} \quad h \to 0^+   \label{behavior_deltah_eps0}
\end{equation}
with $C_r>0$, a uniformly bounded constant w.r.t. $h$.  For the BFS finite element,  Table \ref{tab:infsup} reports some numerical values of $\delta_h$ while Figure \ref{fig:infsupconstant} depicts the evolution of $\sqrt{r}\delta_h$ w.r.t $h$ for $r=1, 10^{-2}, h$ and $r=h^2$.

\begin{table}[http]
\centering
\begin{tabular}{|c|ccccc|}
\hline
$h$  & $1.41\times 10^{-1}$ & $7.01\times 10^{-2}$ & $3.53\times 10^{-2}$ & $1.76\times 10^{-2}$ & $8.83\times 10^{-3}$  \tabularnewline

\hline $r=10^2$  & $4.71\times 10^{-2}$ & $2.53\times 10^{-2}$ & $1.23\times 10^{-2}$ & $6.06\times 10^{-3}$ & $3.01\times 10^{-3}$ \tabularnewline

$r=1$  & $4.17\times 10^{-1}$ & $2.53\times 10^{-1}$ & $1.23\times 10^{-1}$ & $6.05\times 10^{-3}$ & $3.01\times 10^{-2}$ \tabularnewline

$r=10^{-2}$  & $1.474$ & $1.427$ & $1.178$ & $5.99\times 10^{-1}$ & $3.01\times 10^{-1}$ \tabularnewline

$r=h$  & $1.301$ & $1.086$ & $7.65\times 10^{-1}$ & $5.37\times 10^{-1}$ & $3.79\times 10^{-1}$ \tabularnewline

$r=h^2$  & $1.665$ & $1.483$ & $1.485$ & $1.489$ & $1.497$ \tabularnewline \hline
\end{tabular}
\caption{$\delta_h$  w.r.t. $r$  and $h$ - $T=2.$  for the BFS element.} 
\label{tab:infsup}
\end{table}







\begin{figure}[http]
\begin{center}
\includegraphics[scale=0.5]{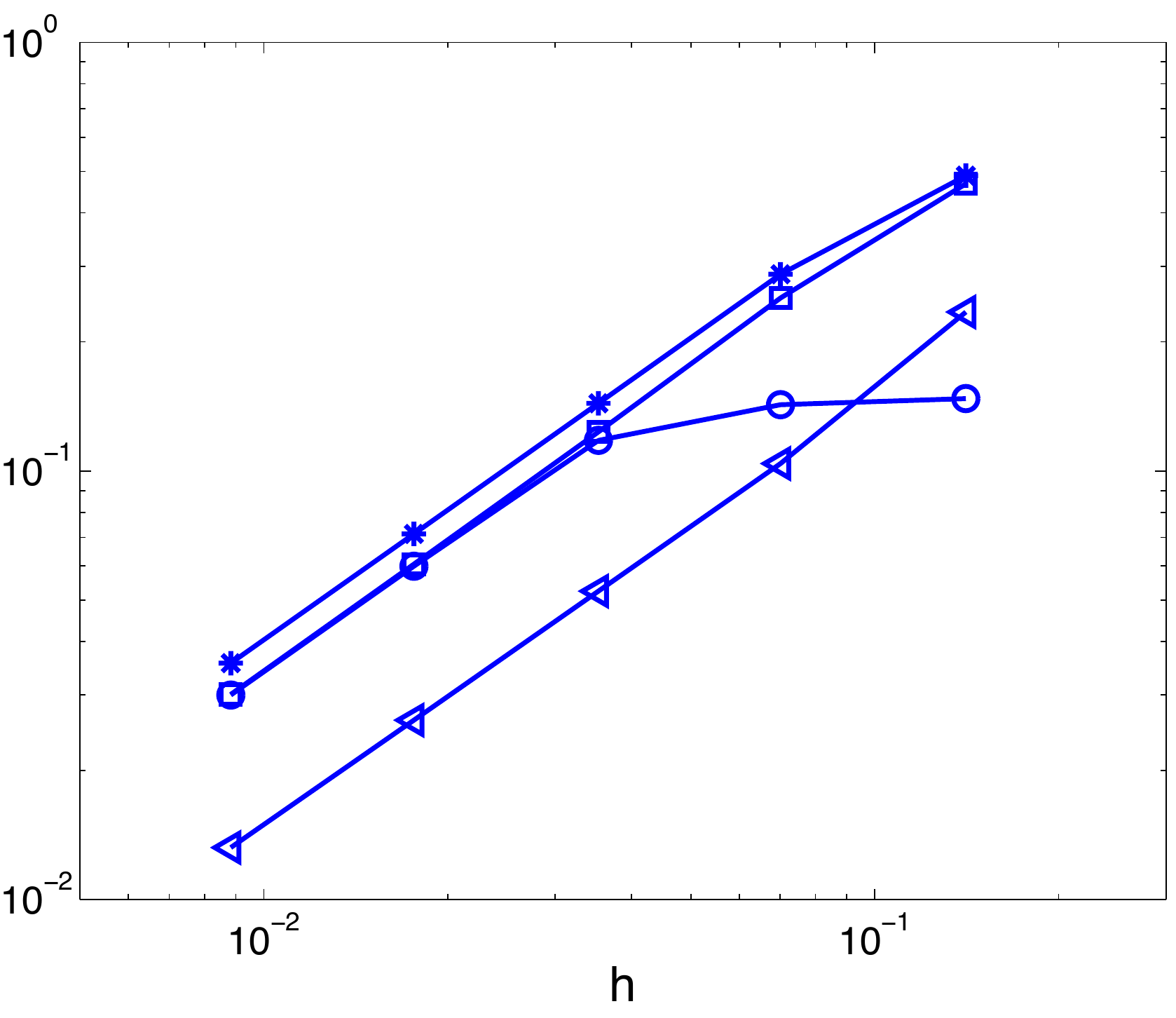}
\caption{BFS finite element - Evolution of $ \sqrt{r} \delta_{h,r}$ with respect to $h$ (see Table \ref{tab:infsup}) for $r=1$ ($\square$), $r=10^{-2}$ ($\circ$), $r=h$ ($\star$) and $r=h^2$ ($<$).}\label{fig:infsupconstant}
\end{center}
\end{figure}

Reporting this behavior in the estimate (\ref{estimate_L2_eps0}), the error $\Vert y-y_h\Vert_{L^2(Q_T)}$ behaves like (taking $\eta=r$)
$$
\Vert y-y_h\Vert_{L^2(Q_T)} \leq  K\max(1,\frac{2}{\sqrt{r}})\biggl(1+\frac{1}{h}+ \frac{1}{\sqrt{r}} \biggr) h^k.
$$
The right hand side is minimal for $r$ of the order one leading to $\Vert y-y_h\Vert_{L^2(Q_T)} \leq  K h^{k-1}$.
This estimate is very likely not optimal but shows the strong convergence of the approximation for regular enough solution. In particular, the estimate (\ref{operator_neumann}) of the boundary term by the $H^2(Q_T)$-norm may result in a loss of precision. 

On the other hand, the same argument for the variable $\lambda_h$ indicates, using (\ref{eq:llh0}) that the estimate of $\Vert \lambda-\lambda_h\Vert_{L^2(Q_T)}$ is 
\begin{equation} \nonumber
\Vert \lambda-\lambda_h \Vert_{L^2(Q_T)}\leq  K_2  \frac{\sqrt{r}}{h}\biggl(1+\frac{1}{h}+\frac{1}{\sqrt{r}}\biggr) h^k.
\end{equation}
The optimal value for the augmentation parameter is now $r=h^2$ leading to $\Vert \lambda-\lambda_h \Vert_{L^2(Q_T)}\leq  K_2 h^{k-1}$. Remark that for $r=h^2$, the discrete-inf sup constant $\delta_h$ remains uniformly bounded by below w.r.t. $h$. 
As observed in \cite[Section 4.]{NC-AM-mixedwave} (see also the numerical experiments of the present paper), if the influence of the parameter $r>0$ is not really important (for $h$ small enough),
the choice $r=h^2$ offers the best rate of convergence. 

As discussed and used in \cite[Section 4.3]{CC-NC-AM}, when $y_h$ is approximated with the BFS and HCT finite element, the quantity $h\Vert Ly_h\Vert_{L^2(Q_T)}$ is asymptotically equivalent w.r.t. $h$ to $\Vert L y_h\Vert_{L^2(0,T; H^{-1}(\Omega))}$. Therefore, taking $r=h^2$ in the augmented mixed formulation amounts to relax the constraint $Ly=0$ in $L^2(Q_T)$ by the weaker one $Ly=0$ in $L^2(0,T, H^{-1}(\Omega))$ enough in practice to approximate weak solutions of (\ref{eq:wave}).
\end{remark}

\subsection{Numerical approximation of the \textit{stabilized} mixed formulation (\ref{mf_stab})}

We address the numerical approximation of the \textit{stabilized} mixed formulation (\ref{mf_stab}) with $\alpha \in (0,1)$, assuming again that $r>0$. Let $h$ be a real parameter. 
Let $Z_h$ and $\widetilde\Lambda_h$ be two closed finite dimensional spaces such that 
$$
Z_h\subset Z, \quad \widetilde\Lambda_h\subset \Lambda, \qquad \forall h>0.
$$
The problem (\ref{mf_stab}) becomes : find $(y_{h}, \lambda_{h}) \in Z_h\times \widetilde\Lambda_h$ solution of 
\begin{equation} \label{mf_stab_h}
\left\{
\begin{array}{rcll}
\noalign{\smallskip} a_{r,\alpha}(y_h, \overline{y}_h) + b_{\alpha}( \lambda_h, \overline{y}_h) & = & l_{1,\alpha}(\overline{y}_h), & \qquad \forall \overline{y}_h \in Z_h \\
\noalign{\smallskip} b_{\alpha}( \overline{\lambda}_h, y_h) - c_{\alpha}(\lambda_h,\overline{\lambda}_h)& = & l_{2,\alpha}(\overline{\lambda}_h), & \qquad \forall \overline{\lambda}_h \in \widetilde\Lambda_h.
\end{array}
\right.
\end{equation}
In view of the properties of the forms $a_{r,\alpha}$, $c_{\alpha}$, $l_{1,\alpha}$ and $l_{2,\alpha}$, this formulation is well-posed. 

Proceeding as in the proof of \cite[Theorem 5.5.2]{brezzi_new}, we first check that the following estimate holds.
\begin{lemma} \label{lemmastab} Let $(y,\lambda)\in Y\times \Lambda$ be the solution of (\ref{mf_stab}) and $(y_h,\lambda_h)\in Z_h\times \widetilde{\Lambda}_h$ be the solution of (\ref{mf_stab_h}). Then we have, 
\begin{align}
\frac{1}{4}\theta \Vert y-y_h\Vert_{Z}^2 +   \frac{1}{4}\alpha \Vert \lambda-\lambda_h\Vert_{\Lambda}^2  \leq  & \biggl(\frac{\Vert a_{r,\alpha}\Vert_{(Z\times Z)^{\prime}}^2}{\theta}+\frac{\Vert b_\alpha\Vert_{(Z\times \Lambda)^{\prime}}^2}{\alpha}+\frac{\theta}{2}\biggr) \inf_{\overline{y}_h\in Z_h} \Vert\overline{y}_h-y\Vert_Z^2 \nonumber \\
& + \biggl(\frac{\Vert b_{\alpha}\Vert_{(Z\times \Lambda)^{\prime}}^2}{\theta}+\alpha+\frac{\alpha}{2}\biggr) \inf_{\overline{\lambda}_h\in \widetilde{\Lambda_h}}\Vert \overline{\lambda}_h-\lambda\Vert_{\Lambda}^2
\end{align}
with  $\Vert a_{r,\alpha}\Vert_{(Z\times Z)^{\prime}}\leq \max(1-\alpha, \eta^{-1}r)$,  $\Vert b_{\alpha}\Vert_{(Z\times \Lambda)^{\prime}}\leq (C_{\Omega,T}+\alpha)\,\eta^{-1}$. Parameter $\theta$ is defined in Proposition \ref{pr:mfs}.
\end{lemma}
Concerning the space $\widetilde\Lambda_h$, since $L\lambda_h$ should belong to $L^2(0,T,H^{-1}(\Omega))$, a natural choice is 
\begin{equation}
\widetilde\Lambda_h=\{\lambda\in Z_h; \lambda(\cdot,0)=\lambda_t(\cdot,0)=0\}.  \label{choice_lambdah_alpha}
\end{equation}
where $Z_h\subset Z$ is defined by (\ref{defZh}).
Then, using Lemma \ref{lemmastab} and the estimate (\ref{eq:ch}), we obtain the following result. 
\begin{proposition}[BFS element for $N=1$ - Rate of convergence for the norm $Z\times \Lambda$]\label{propZLambda}
Let $h>0$, let $k \leq 2$ be a positive integer and $\alpha\in (0,1)$. Let $(y,\lambda)$ and $(y_h,\lambda_h)$ be the solution of (\ref{mf_stab}) and (\ref{mf_stab_h}) respectively. If $(y,\lambda)$ belongs to $H^{k+2}(Q_T)\times H^{k+2}(Q_T)$, then there exists a positive constant $K=K(\Vert y\Vert_{H^{k+2}(Q_T)}, \Vert c \Vert_{C^1(\overline{Q_T})},\Vert d \Vert_{L^{\infty}(Q_T)},\alpha,r,\eta)$  independent of $h$, such that 
\begin{equation}
\begin{aligned}
\Vert y-y_h \Vert_Z + \Vert \lambda-\lambda_h \Vert_{\Lambda} \leq K h^k.
\end{aligned}
\end{equation}
\end{proposition}
\textsc{Sketch of the proof-} The estimate of $\Vert y-\Pi_{Z_h,\mathcal{T}_h}(y)\Vert_Z$ for any $y\in Z$ in term of $\mathcal{O}(h^{k})$,   is detailed in the proof of Proposition \ref{propeps0}: precisely, we refer to (\ref{eq:ch}). Similarly, we write that, for any $\lambda\in \Lambda$, 
\begin{equation}
\begin{aligned}
\inf_{\lambda_h\in \widetilde\Lambda_h}\Vert \lambda-\lambda_h\Vert_{\Lambda} & \leq \Vert \lambda-\Pi_{\widetilde\Lambda_h,\mathcal{T}_h}(\lambda)\Vert_{\Lambda}\\
& \leq   \biggl(  \Vert L(\lambda-\Pi_{\widetilde\Lambda_h,\mathcal{T}_h}(\lambda))\Vert^2_{L^2(H^{-1}(\Omega))} + \Vert c(x) (\lambda-\Pi_{\widetilde\Lambda_h,\mathcal{T}_h}(\lambda))\Vert^2_{L^2(\Gamma_T)}\biggr)^{1/2}\\
& \leq  \biggl(  \Vert L(\lambda)-L(\Pi_{\widetilde\Lambda_h,\mathcal{T}_h}(\lambda))\Vert^2_{L^2(Q_T)} + \Vert c\Vert^2_{L^{\infty}(\Gamma_T)} \Vert \lambda-\Pi_{\widetilde\Lambda_h,\mathcal{T}_h}(\lambda)\Vert^2_{H^{1/2}(Q_T)}\biggr)^{1/2}\\
& \leq K (\Vert c \Vert_{C^1(\overline{Q_T})},\Vert d \Vert_{L^{\infty}(Q_T)})   \Vert \lambda-\Pi_{\widetilde\Lambda_h,\mathcal{T}_h}(\lambda)\Vert_{H^{2}(Q_T)}
\end{aligned}
\end{equation}
for which the $\mathcal{O}(h^k)$ estimate follows, in view of the definition of $\widetilde\Lambda_h$.

In particular, arguing as in the previous section, using Proposition \ref{propZLambda}, the observability estimate (\ref{iobs}) for the variable $y$ and the estimate (\ref{global_lambda}) for the variable $\lambda$, we get the following desired global estimate:

\begin{theorem}[BFS element for $N=1$ - Rate of convergence for the $L^2(Q_T)$ norm]
Assume that the hypothesis (\ref{iobs}) holds. Let $h>0$, let an integer $k\leq 2$. Let $(y,\lambda)$ and $(y_{h},\lambda_{h})$ be the solution of (\ref{mf_stab}) and (\ref{mf_stab_h}) respectively. If the solution $(y, \lambda)$ belongs to $H^{k+2}(Q_T)\times H^{k+2}(Q_T)$, then there exist two positive constants $K_i=K_i(\Vert y\Vert_{H^{k+2}(Q_T)}, \Vert \lambda\Vert_{H^{k+2}(Q_T)}, \Vert c \Vert_{C^1(\overline{Q_T})},\Vert d \Vert_{L^{\infty}(Q_T)},\alpha,r,\eta)$, $i=1,2$ independent of $h$ such that
\begin{equation}
\Vert y-y_h\Vert_{L^2(Q_T)} \leq  K_1 \frac{h^k}{\sqrt{\eta}},  \quad \Vert \lambda-\lambda_h\Vert_{L^2(Q_T)} \leq  K_2 h^k. 
\end{equation}
\end{theorem}
We emphasize that, within the stabilized formulation, these estimates do not depend of any discrete inf-sup constant. In particular, the positive augmentation parameter $r$ can be chosen arbitrarily.

\subsection{Numerical approximation of the mixed formulation (\ref{eq:mfeps})}

We now consider the numerical analysis of the mixed formulation (\ref{eq:mfeps}) where both the solution $y$ and the spatial part $\mu$ of the source are unknown. 
We take a strictly positive augmentation parameter $r$.

Let $Y_h$ and $\Lambda_h$ be two finite dimensional spaces parametrized by the variable $h$ such that $Y_h\subset Y, \Lambda_h\subset L^2(Q_T)$, and $\dim(Y_h)\geq \dim(\Lambda_h)$ for every $h>0$.
Then, we can introduce the following approximated problems: find $((y_h,\mu_h),\lambda_h)\in Y_h\times \Lambda_h$ solution of 
\begin{equation} \label{eq:mfepsh}
\left\{
\begin{array}{rcll}
\noalign{\smallskip} a_r((y_h,\mu_h),(\overline{y}_h,\overline{\mu}_h)) + b((\overline{y}_h,\overline{\mu}_h), \lambda_h) & = & l(\overline{y}_h), & \qquad \forall (\overline{y}_h,\overline{\mu}_h)\in Y_h \\
\noalign{\smallskip} b((y_h,\mu_h), \overline{\lambda}_h) & = & 0, & \qquad \forall \overline{\lambda}_h \in \Lambda_h.
\end{array}
\right.
\end{equation}

Again, for any $r>0$, the coercivity of the form $a_r$ holds true on the whole space $Y$: 
\begin{equation}\nonumber
a_r((y,\mu),(y,\mu)) \geq  \frac{r}{\eta} \Vert (y,\mu)\Vert^2_Y, \quad \forall (y,\mu)\in Y
\end{equation}
so a fortiori, uniformly on the subspace $Y_h$. Therefore, assuming that the discrete inf-sup constant defined by
\begin{equation}
\delta_h:=\inf_{\lambda_h\in \Lambda_h}\sup_{(y_h,\mu_h)\in Y_h}  \frac{b((y_h,\mu_h),\lambda_h)}{\Vert \lambda_h\Vert_{L^2(Q_T)}\Vert (y_h,\mu_h)\Vert_Y} \label{infsupepsdiscret}
\end{equation}
is, for any $h>0$, strictly positive, there exists a unique solution $(y_h,\mu_h)$ of (\ref{eq:mfepsh}) and the following holds true: 

\begin{proposition} \label{prop_estimateh_eps}
Let $h>0$. Let $((y,\mu),\lambda)$ and $((y_h,\mu_h),\lambda_h)$ be the solution of (\ref{eq:mfeps}) and of (\ref{eq:mfepsh}) respectively. Let $\delta_h$ be the discrete inf-sup constant defined by (\ref{infsupepsdiscret}). Then,
\begin{align}
\label{eq:est_e01eps}
& \Vert (y,\mu)-(y_h,\mu_h) \Vert_{Y} \leq  2\biggl(1+\frac{1}{\sqrt{\eta} \delta_h}\biggr)d((y,\mu),Y_h)+\frac{1}{\sqrt{\eta}} d(\lambda,\Lambda_h), \\
\label{eq:est_e02eps}
& \Vert \lambda-\lambda_h \Vert_{L^2(Q_T)}  \leq  \biggl(2+\frac{1}{\sqrt{\eta}\delta_h}\biggr)\frac{1}{\delta_h}d(y,Z_h)+
\frac{3}{\sqrt{\eta} \delta_h} d(\lambda,\Lambda_h)
\end{align}
where  $d(\lambda,\Lambda_h)$ is as in Proposition \ref{prop_estimateh_eps0}
and
\begin{equation}
\label{def_distanceyYh}
\begin{aligned}
d((y,\mu),Y_h):=& \inf_{(y_h,\mu_h)\in Y_h}\Vert (y,\mu)-(y_h,\mu_h)\Vert_Y\\
=& \inf_{y_h \in Y_h}   \biggl( \Vert c(x)\partial_{\nu} (y-y_h)\Vert^2_{L^2(\Gamma_T)}+  \eta \Vert L(y-y_h)-\sigma(\mu-\mu_h) \Vert^2_{L^2(Q_T)}\biggr)^{1/2}.
\end{aligned}
\end{equation}
\end{proposition}
\textsc{Proof-} The proof is similar to the proof of Proposition \ref{prop_estimateh_eps0} using again that $\Vert a_r\Vert_{(Y\times Y)^{\prime}}\leq 1$ and $\Vert b\Vert_{(Y\times \Lambda)^{\prime}}\leq  \frac{1}{\sqrt{\eta}}$.  \Fin

The finite dimensional problem (\ref{eq:mfepsh}) reads as follows: find $\{y_h,\mu_h\}\in \mathbb{R}^{n_h}$ and $\{\lambda_h\}\in \mathbb{R}^{m_h}$ such that 
\begin{equation} \label{matrixmfepsh}
\left(
\begin{array}{cc}
A_{r,h} &  B_h^T \\
B_h & 0   
\end{array}
\right)_{\mathbb{R}^{n_h+m_h,n_h+m_h}}
\left(
\begin{array}{c}
\{y_h,\mu_h\}  \\
\{\lambda_h\}    
\end{array}
\right)_{\mathbb{R}^{n_h+m_h}}  =
\left(
\begin{array}{c}
L_h \\
0   
\end{array}
\right)_{\mathbb{R}^{n_h+m_h}}
\end{equation}
with matrices similarly defined as in the previous section. 

Since the variable $\mu$ is a function of $x\in \Omega$ only, we first introduce a triangulation $\mathcal{T}_{\Delta x}$ such that $\overline{\Omega}=\cup_{K\in \mathcal{T}_{\Delta x}} K$ and we assume that $\{\mathcal{T}_{\Delta x}\}_{\Delta x>0}$ is a regular family. We denote 
$
\Delta x:=\max\{\textrm{diam}(K), K\in \mathcal{T}_{\Delta x}\}.
$
We then introduce the subspace $M_{\Delta x}$ of $H^{-1}(\Omega)$ defined by
\begin{equation}
M_{\Delta x}=\{\mu_{\Delta x}\subset C^0(\overline{\Omega}):  \mu_{\Delta x}\vert_K\in \mathbb{Q}_x(K), \quad \forall K\in \mathcal{T}_{\Delta x}\},   \label{defMDeltax}
\end{equation}
where $\mathbb{Q}_x(K)$ denotes the space of affine functions in the variable $x$ in the element $K$.

Similarly, we define a subdivision $\mathcal{T}_{\Delta t}$ of $[0,T]$ such that $[0,T]=\cup_{Q\in \mathcal{T}_{\Delta t}} Q$ and denote $\Delta t:=\max\{\textrm{diam}(Q), Q\in \mathcal{T}_{\Delta t}\}$. Then, we consider the triangulation defined by $\{\mathcal{T}_h\}:= \{\mathcal{T}_{\Delta x}\}_{\Delta x>0}\otimes \{\mathcal{T}_{\Delta t}\}_{\Delta t>0}$ such that  
$$
\overline{Q_T}=\bigcup_{K\in \mathcal{T}_{\Delta x}, Q\in \mathcal{T}_{\Delta t}}  \{ K \times Q\}.
$$
Again, we denote by $h:=\max\{\textrm{diam}(K\times Q), K\in \mathcal{T}_{\Delta x}, Q\in \mathcal{T}_{\Delta t}\}$.
The triangulation $\{\mathcal{T}_h\}_{h}$ is a regular family for $Q_T$ as soon as $\{\mathcal{T}_{\Delta x}\}_{\Delta x>0}$ is a regular family for $\Omega$.
 
As in the previous section, we now go on in the one dimensional case in space. 

\subsubsection{$C^1$-finite elements and order of convergence for $N=1$}

We only discuss the BFS finite element and introduce the space $Z_h$ as follows : 
\begin{equation}
Z_h=\{z_h\subset C^1(\overline{Q_T}):  z_h\vert_M\in \mathbb{P}(M), \quad \forall M\in \mathcal{T}_h, \,\, z_h=0 \,\,\textrm{on}\,\,\Sigma_T, \,\, z_h=z_{h,t}=0 \,\, \textrm{on}\,\, \Omega\times \{0\}\},   \label{defZhb}
\end{equation}
with $\mathbb{P}(M):=\mathbb{P}_{3,x}(Q)\times \mathbb{P}_{3,t}(K)$,   for all $(Q,K)\in (\mathcal{T}_{\Delta x},T_{\Delta t})$. Finally,  we define the space $Y_h$ by
\begin{equation}\nonumber
Y_h:=\{y_h=(z_h,\mu_{\Delta x}):   z_h\in Z_h, \mu_{\Delta x}\in M_{\Delta x}\}
\end{equation}
so that, for each value of $t\in [0,T]$, the variables $z_h(\cdot,t)$ and $\mu_{\Delta x}$ share the same triangulation with respect to the variable $x$. We check that $Y_h\subset Y$ for each $h>0$. 
In the sequel, for simplicity, we use the notation $\mu_h$ for $\mu_{\Delta x}$.

We then have the following result: 
\begin{proposition}[BFS element for $N=1$ - Rate of convergence for the norm $Y$]   \label{propeps}
Let $h>0$, let $k,q\leq 2 $ be two nonnegative integers. Let $((y,\mu),\lambda)$ and $((y_h,\mu_h),\lambda_h)$ be the solution of (\ref{eq:mfeps}) and (\ref{eq:mfepsh}) respectively. If the solution $((y,\mu),\lambda)$ belongs to $(H^{k+2}(Q_T)\times H^q(\Omega))\times H^k(Q_T)$, then there exist two positive constants 
\[ K_i=K_i(\Vert y\Vert_{H^{k+2}(Q_T)}, \Vert \mu\Vert_{H^q(\Omega)}, \Vert \lambda\Vert_{H^k(Q_T)} ,\Vert c \Vert_{C^1(\overline{Q_T})},\Vert d \Vert_{L^{\infty}(Q_T)},\Vert \sigma\Vert_{L^{\infty}([0,T])}), \qquad i \in \{ 1,2\},\]
independent of $h$, such that 
\begin{align}
\Vert (y,\mu)-(y_h,\mu_h) \Vert_Y & \leq K_1 \biggl(1+\frac{1}{\sqrt{\eta} \delta_h}+ \frac{1}{\sqrt{\eta}} \biggr) h^k + K_1 \biggl(1+\frac{1}{\sqrt{\eta} \delta_h}\biggr) (\Delta x)^q , \label{eq:yyheps}\\
\Vert \lambda-\lambda_h \Vert_{L^2(Q_T)}&\leq  K_2\biggl(\biggl(1+\frac{1}{\sqrt{\eta} \delta_h}\biggr)\frac{1}{\delta_h} +\frac{1}{\sqrt{\eta}\delta_h} \biggr) h^k+ K_2\biggl(\biggl(1+\frac{1}{\sqrt{\eta} \delta_h}\biggr)\frac{1}{\delta_h} \biggr) (\Delta x)^q\label{eq:llheps}.
\end{align}
\end{proposition}
\par\noindent
\textsc{Proof -} We proceed as in the proof of Proposition \ref{propeps0}. We write that 
\begin{equation}
\nonumber
\begin{aligned}
 \Vert L(y-\Pi_{Y_h,\mathcal{T}_h}(y))&-\sigma(\mu-\Pi_{Z_h,\mathcal{T}_h}(\mu))\Vert_{L^2(Q_T)}  \\
&\leq \Vert L(y-\Pi_{Z_h,\mathcal{T}_h}(y))\Vert_{L^2(Q_T)}+ \Vert \sigma\Vert_{L^{\infty}([0,T])}
\Vert \mu-\Pi_{M_{\Delta x},\mathcal{T}_{\Delta x}}(\mu)\Vert_{L^2(\Omega)}\\
& \leq K \Vert y-\Pi_{Z_h,\mathcal{T}_h}(y)\Vert_{H^2(Q_T)}+ \Vert \sigma\Vert_{L^{\infty}([0,T])}
\Vert \mu-\Pi_{M_{\Delta x},\mathcal{T}_{\Delta_x}}(\mu)\Vert_{L^2(\Omega)}
\end{aligned}
\end{equation}
where $K$ is from (\ref{eq:est2h}). 

We use that, if $\mu\in H^q(\Omega)$, $q\leq 2$, there exists a positive constant $C_3=C_3(\Vert \mu\Vert_{H^q(\Omega)})$ such that 
$\Vert \mu-\Pi_{M_{\Delta x},\mathcal{T}_{\Delta_x}}(\mu)\Vert_{L^2(\Omega)}\leq C_3 (\Delta x)^q$ for any $\Delta x>0$. Consequently, from (\ref{eq:est1h}), we have the estimate
\begin{equation}\nonumber
 \Vert L(y-\Pi_{Y_h,\mathcal{T}_h}(y))-\sigma(\mu-\Pi_{Z_h,\mathcal{T}_h}(\mu))\Vert_{L^2(Q_T)} \leq KC_2 h^k + C_3 (\Delta x)^q \leq (KC_2 +C_3)(h^k+(\Delta x)^q).
 \end{equation}
 This leads to the following estimates : 
 
 \begin{equation}
 \nonumber
 \begin{aligned}
 d((y,\mu),Y_h) & \leq \biggl((\Vert c\Vert^2_{L^{\infty}(\Gamma_T)} C_P^2 +\eta K^2) C^2_2 h^{2k}+\eta \Vert \sigma\Vert^2_{L^\infty([0,T])}C_3^2 (\Delta x)^{2q}\biggr)^{1/2} \\
 & \leq \sqrt{2} \sqrt{(\Vert c\Vert^2_{L^{\infty}(\Gamma_T)} C_P^2 +\eta K^2)} C_2 h^k + \sqrt{2}\sqrt{\eta} \Vert \sigma\Vert_{L^\infty([0,T])} C_3 (\Delta x)^q
 \end{aligned}
 \end{equation}
 Using (\ref{eq:est_e01eps}), we then get
\begin{equation*}
\begin{aligned}
 \Vert (y,\mu)-&(y_h,\mu_h) \Vert_{Y} \\
 \leq & 2\sqrt{2}\biggl(1+\frac{1}{\sqrt{\eta} \delta_h}\biggr)\biggl[\sqrt{(\Vert c\Vert^2_{L^{\infty}(\Gamma_T)} C_P^2 +\eta K^2)} C_2 h^k + \sqrt{\eta} \Vert \sigma\Vert_{L^\infty([0,T])}) C_3 (\Delta x)^q\biggr]\\
 &+\frac{1}{\sqrt{\eta}} C_1h^k.
 \end{aligned}
\end{equation*}
Similarly, using (\ref{eq:est_e02eps}), we get 

 \begin{equation*}
\begin{aligned}
 \Vert \lambda-\lambda_h& \Vert_{L^2(Q_T)} \\
 \leq & \sqrt{2}\biggl(2+\frac{1}{\sqrt{\eta} \delta_h}\biggr)\frac{1}{\delta_h}\biggl[\sqrt{(\Vert c\Vert^2_{L^{\infty}(\Gamma_T)} C_P^2 +\eta K^2)} C_2 h^k + \sqrt{\eta} \Vert \sigma\Vert_{L^\infty([0,T])} C_3 (\Delta x)^q\biggr]\\
 &+\frac{3}{\sqrt{\eta}\delta_h} C_1h^k.
 \end{aligned}
\end{equation*}
The proposition then follows from the last two estimates. 
\Fin

It remains now to deduce the convergence of the approximation $y_h$ for a global norm, typically $L^2(Q_T)$. We write that $y-y_h$ solves 
the problem
\begin{equation}
\nonumber
\left\{
\begin{aligned}
& L(y-y_h)=\sigma(\mu-\mu_h)+  L(y-y_h)-\sigma(\mu-\mu_h)  && \quad \text{in } Q_T \\
& ((y-y_h), (y-y_h)_t)(\cdot, 0)=(0,0) && \quad \text{in } \Omega\\
& y-y_h=0  && \quad \text{on } \Sigma_T
\end{aligned}
\right.
\end{equation}
leading to 
$$
\begin{aligned}
\nonumber
\Vert y-y_h\Vert_{L^2(Q_T)} & \leq C_{\Omega,T}\biggl( \Vert \sigma\Vert_{L^{\infty}(0,T)}\Vert (\mu-\mu_h)  \Vert_{H^{-1}(\Omega)} +  \Vert L(y-y_h)-\sigma(\mu-\mu_h)\Vert_{L^2(Q_T)} \biggr).
\end{aligned}
$$
On the other hand, proceeding as before, assuming (\ref{iobs2}), we have 
\begin{equation}
\nonumber
\Vert \mu-\mu_h\Vert^2_{H^{-1}(\Omega)} \leq C_{obs} \biggl(\Vert c(x)(\partial_\nu y-\partial_\nu y_h)\Vert^2_{L^2(\Gamma_T)} + \Vert L(y-y_h)-\sigma(\mu-\mu_h)\Vert^2_{L^2(Q_T)}\biggr)
\end{equation}
leading to 
\begin{equation}
\nonumber
\begin{aligned}
\Vert y-y_h\Vert_{L^2(Q_T)}&\leq C_{\Omega,T} \biggl( \Vert\sigma\Vert_{L^{\infty}([0,T])} \sqrt{C_{obs}} \Vert c(x)\partial_{\nu}(y-y_h)\Vert_{L^2(\Gamma_T)}\\
&\hspace{2cm}+ (1+\Vert\sigma\Vert_{L^{\infty}([0,T])} \sqrt{C_{obs}})\Vert L(y-y_h)-\sigma(\mu-\mu_h)\Vert_{L^2(Q_T)} \biggr)\\
& \leq  C_{\Omega,T} (1+\Vert\sigma\Vert_{L^{\infty}([0,T])} \sqrt{C_{obs}})\biggl(\Vert c(x)\partial_{\nu}(y-y_h)\Vert_{L^2(\Gamma_T)}\\
&\hspace{3cm}+ \Vert L(y-y_h)-\sigma(\mu-\mu_h)\Vert_{L^2(Q_T)} \biggr)\\
& \leq \sqrt{3}C_{\Omega,T} (1+\Vert\sigma\Vert_{L^{\infty}([0,T])} \sqrt{C_{obs}})\max(1,\frac{1}{\sqrt{\eta}}) \Vert (y-y_h,\mu-\mu_h)\Vert_Y.
\end{aligned}
\end{equation}

Therefore, in view of Proposition \ref{propeps}, we get the following a priori estimate for the $L^2$-norm :
\begin{theorem}[BFS element for $N=1$ - Rate of convergence for the $L^2(Q_T)$-norm]   \label{theps}
Let $h>0$, let $k,q\leq 2 $ be two nonnegative integers. Let $((y,\mu),\lambda)$ and $((y_h,\mu_h),\lambda_h)$ be the solution of (\ref{eq:mfeps}) and (\ref{eq:mfepsh}) respectively. If  $((y,\mu),\lambda)$ belongs to $(H^{k+2}(Q_T)\times H^q(\Omega))\times H^k(Q_T)$, then there exists a positive constant 
\[ K=K(\Vert y\Vert_{H^{k+2}(Q_T)}, \Vert \mu\Vert_{H^k(\Omega)}, \Vert c \Vert_{C^1(\overline{Q_T})},\Vert d \Vert_{L^{\infty}(Q_T)}),\]
independent of $h$, such that 
\begin{equation}
\nonumber
\begin{aligned}
\Vert y-y_h\Vert_{L^2(Q_T)} & \leq    K C_{\Omega,T} (1+\Vert\sigma\Vert_{L^{\infty}([0,T])} \sqrt{C_{obs}})\max(1,\frac{1}{\sqrt{\eta}}) \\
&\hspace{2cm}\biggl[\biggl(1+\frac{1}{\sqrt{\eta} \delta_h}+ \frac{1}{\sqrt{\eta}} \biggr) h^k + \biggl(1+\frac{1}{\sqrt{\eta} \delta_h}\biggr) (\Delta x)^q\biggr].
\end{aligned}
\end{equation}
\end{theorem}
\par\noindent

\begin{remark}
We have presented along Section \ref{sec_numer} error estimates in the one-dimensional case in the case where the space $Z_h$ (and $\Lambda_h$ for the present section) is based on the BFS finite element. Similar results may be obtained within the HCT finite element. Precisely, we refer to \cite[ch. VII, sec. 48, p. 295]{ciarletfem} where the following interpolation estimates for the HCT element are developed: 
\begin{equation}\nonumber
\Vert y-\Pi_{Z_h,\mathcal{T}_h}(y)\Vert_{L^2(Q_T)} \leq  C h^{k+2}, \quad  \Vert y-\Pi_{Z_h,\mathcal{T}_h}(y)\Vert_{H^2(Q_T)} \leq  C  h^k, \quad \forall y\in H^2(Q_T)
\end{equation}
for $k\in \{0,1\}$ and some constant $C>0$.
\end{remark}

\section{Numerical experiments}\label{sec_experiment}

We now report and discuss some numerical experiments corresponding to mixed formulations (\ref{eq:mfh}), (\ref{mf_stab_h}) and (\ref{eq:mfepsh}) for $N=1$ and $N=2$ respectively.

\subsection{Reconstruction of the solution - One dimensional case ($N = 1$)}
\label{sec:1d}
\begin{figure}[!ht]
\centering
\begin{tabular}{cc}
\includegraphics[width=0.52\textwidth]{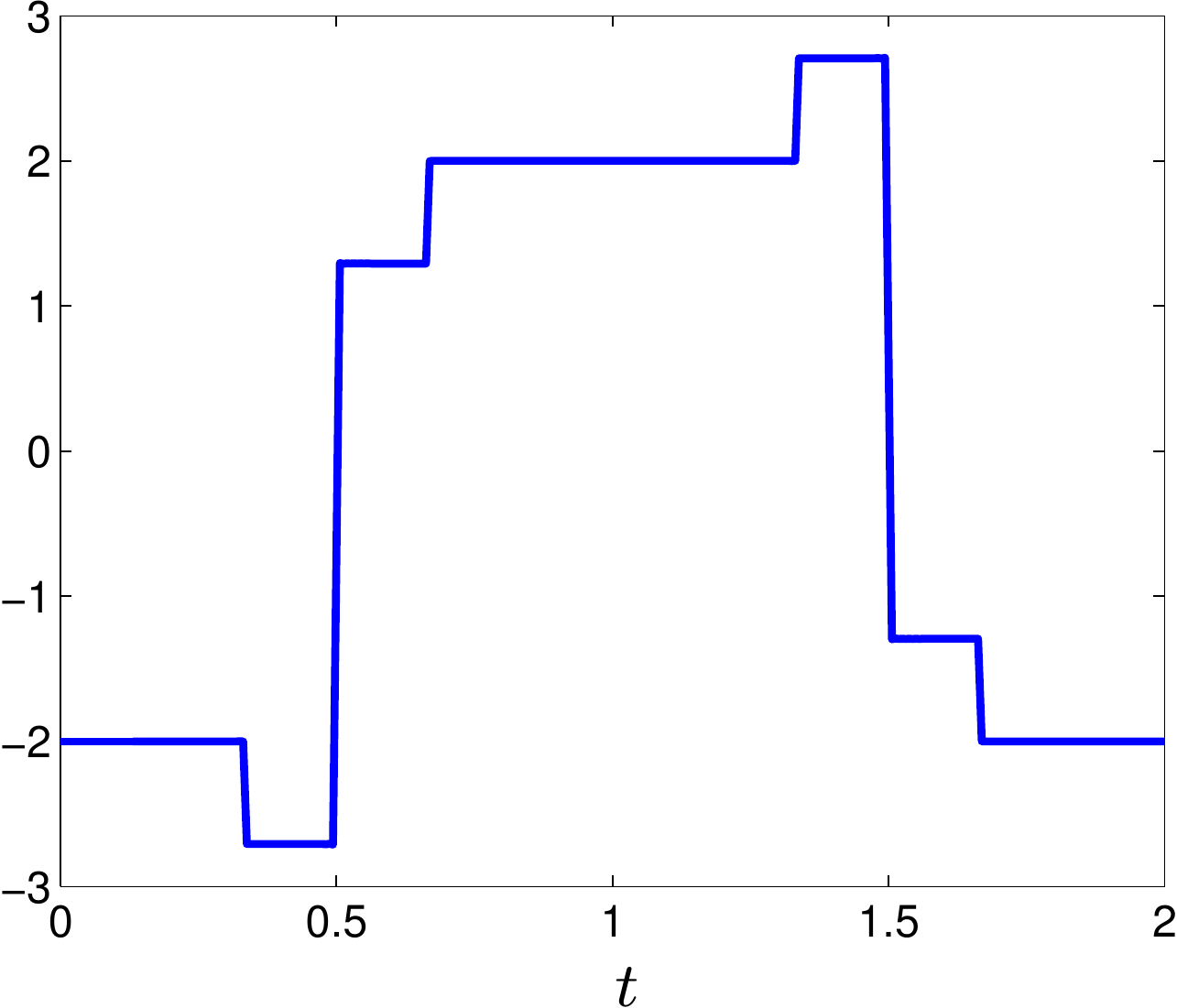} & \hspace{1cm}\includegraphics[width=0.28\textwidth]{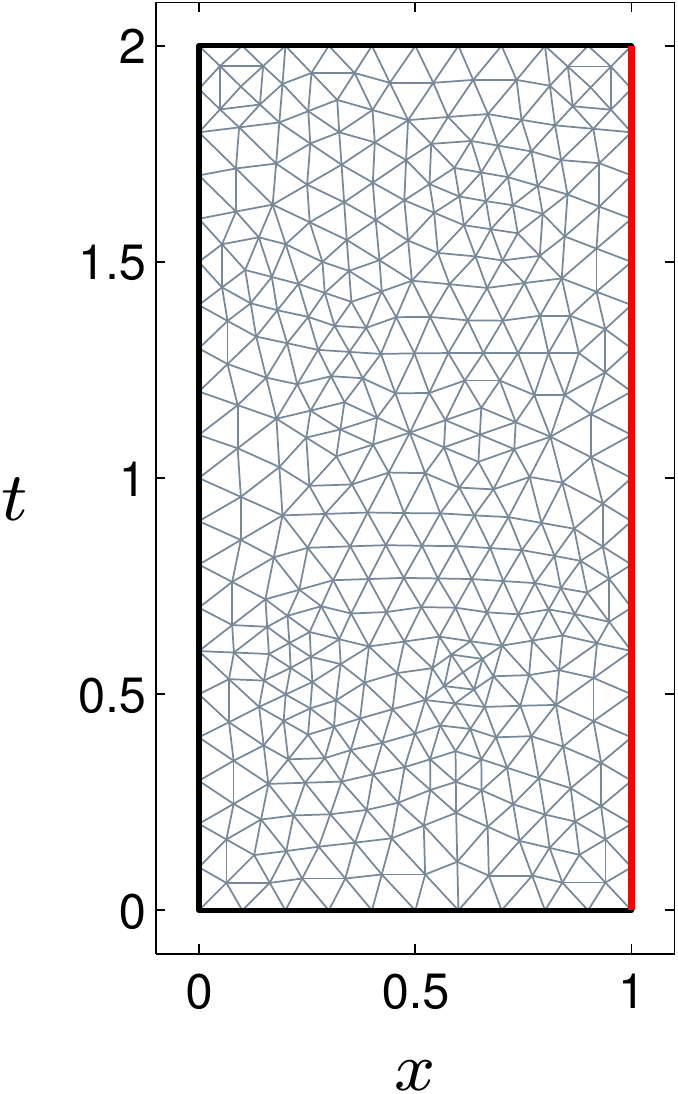} \tabularnewline
\end{tabular}
\caption{The observation $y_{\nu, obs}$ on $\{1\}\times (0,T)$ associated to initial data \textbf{EX1}. Example of mesh of the domain $Q_T$.}
\label{fig:yobs_EX}
\end{figure}

We take $\Omega=(0,1)$ and $\Gamma=\{1\}\subset \partial\Omega$, $T=2$, $c\equiv 1$ and $d\equiv 0$. We check that for these data, the inequality (\ref{iobs}) holds true. 
Moreover, in order to check the convergence of the numerical approximations, we consider explicit solutions of (\ref{eq:wave}). We first define the following initial condition in $H_0^1(\Omega)\times L^2(\Omega)$ (see \cite{cindea_moireau}): 
\begin{equation}
\nonumber
(\textbf{EX1}) \quad
y_0(x)=1-\vert 2x-1\vert, \quad y_1(x)=  1_{(1/3,2/3)}(x),  \qquad x\in (0,1)
\end{equation}
and $f=0$. The corresponding solution of (\ref{eq:wave}) is given by 
\begin{equation}
y(x,t)=\sum_{k>0}  \biggl(a_k \cos(k\pi  t) + \frac{b_k}{k\pi}\sin(k\pi t)\biggr) \sqrt{2} \sin(k\pi x)  \label{yobs_EX}
\end{equation}
with 
$$
a_k=\frac{4\sqrt{2}}{\pi^2 k^2} \sin(\pi k/2), \quad b_k=  \frac{1}{\pi k}(\cos(\pi k/3)-\cos(2\pi k/3)), \quad k>0.
$$
In particular $\Vert y\Vert_{L^2(Q_T)}=\sqrt{\sum_{k>0}  (a_k^2+\frac{b_k^2}{(k\pi)^2})}\approx 0.58663$.   
The corresponding normal derivative $\partial_{\nu} y \vert_{\Gamma_T}=y_x(1,t)$ is depicted in Figure \ref{fig:yobs_EX}. This example is rather stiff: the normal derivative is in $L^2(0,T)$ but is discontinuous (in view of the regularity of the initial condition). We compute $\Vert \partial_{\nu}y\Vert_{L^2(\Gamma_T)}=\sqrt{  2\sum_{k>0} ((k\pi)^2 a_k^2+b_k^2)}\approx 8.33298$.

We recall that the direct method amounts to solve, for any $h$, the linear system (\ref{matrixmfh}). We use exact integration methods developed in \cite{dunavant} for the evaluation of the coefficients of the matrices. Moreover, the linear system (\ref{matrixmfh}) is solved using the LU decomposition method.

We first consider the BFS finite element with uniform triangulation (each element $K$ of the triangulation $\mathcal{T}_h$ is a rectangle of lengths $\Delta x$ and $\Delta t$ so that $h=\sqrt{(\Delta x)^2+(\Delta t)^2}$).
 Table \ref{tab:ex1_rh2_T2} collects some numerical values with respect to $h$ for $r=h^2$ and for $\Delta x=\Delta t$. We observe the following behavior with respect to $h$: 
\begin{equation}
\nonumber
r=h^2:\qquad
\begin{aligned}
&\frac{\Vert y-y_h\Vert_{L^2(Q_T)}}{\Vert y\Vert_{L^2(Q_T)}}=\mathcal{O}(h^{1.20}), \quad \frac{\Vert \partial_\nu(y-y_h)\Vert_{L^2(\Gamma_T)}}{\Vert \partial_{\nu} y\Vert_{L^2(\Gamma_T)}}=\mathcal{O}(h^{0.59}), \\
& \Vert \lambda_h\Vert_{L^2(Q_T)}=\mathcal{O}(h^{1.11}), \quad \Vert L y_h\Vert_{L^2(Q_T)}=\mathcal{O}(h^{-0.29}).
\end{aligned}
\end{equation}
The evolution of the norm  $\Vert y-y_h\Vert_{L^2(Q_T)}$ suggests that the unique solution $y$ is correctly reconstructed from the observation. Moreover, in agreement with Remark \ref{rk_lambda_sys}, since $y_{\nu,obs}$ is by construction the restriction to $\{1\}\times (0,T)$ of the normal derivative of a solution of (\ref{eq:wave}), we check, that the sequence $\lambda_h$, approximation of $\lambda$, vanishes as $h\to 0$. Eventually, in view of the Remark \ref{better_estimate}, we have $r\Vert Ly_h\Vert_{L^2(Q_T)}=h\Vert Ly_h\Vert_{L^2(Q_T)}\approx \Vert Ly_h\Vert_{L^2(0,T, H^{-1}(0,1))}\approx\mathcal{O}(h^{1-0.29})$, so that $y_h$ approximates correctly the unique corresponding weak solution of (\ref{eq:wave}). 

We also check that the minimization of the functional $J_r^{\star\star}$ introduced in Proposition \ref{prop_equiv_dual} leads exactly to the same result approximation: we recall that the minimization of the functional $J_r^{\star\star}$ corresponds to the resolution of the associated mixed formulation by an iterative Uzawa type method. The minimization is done using a conjugate gradient algorithm (we refer to \cite[Section 2.2]{NC-AM-mixedwave} for the algorithm). Each iteration amounts to solve a linear system involving the matrix $A_{r,h}$ which is sparse, symmetric and positive definite. The Cholesky method is used. The performance of the algorithm depends on the condition number of the operator $\mathcal{P}_r$: precisely, it is known that (see for instance \cite{Daniel1971}), 
\[
\Vert \lambda^n-\lambda\Vert_{L^2(Q_T)} \leq 2\sqrt{\nu(\mathcal{P}_r)}   \biggl(\frac{\sqrt{\nu(\mathcal{P}_r)}-1}{\sqrt{\nu(\mathcal{P}_r)}+1}  \biggr)^n \Vert \lambda^0-\lambda\Vert_{L^2(Q_T)}, \quad \forall n\geq 1
\]
where $\lambda$ minimizes $J_r^{\star\star}$. $\nu(\mathcal{P}_r)=\Vert \mathcal{P}_r\Vert_{\mathcal{L}(L^2(Q_T))} \Vert \mathcal{P}_r^{-1}\Vert_{\mathcal{L}(L^2(Q_T))}$ denotes the condition number of the operator $\mathcal{P}_r$. As discussed in \cite[Section 4.4]{NC-AM-mixedwave}, the condition number of $\mathcal{P}_r$ restricted to $\Lambda_h\subset L^2(Q_T)$ (which coincides with the condition number of the matrix $B_h A_{r,h}B_h^T$) behaves asymptotically as $C_r^{-2}h^{-2}$, where $C_r$ is the constant appearing in (\ref{behavior_deltah_eps0}). This quadratic behavior is the typical one for well-posed elliptic problems. Table \ref{tab:ex1_rh2_T2} reports the number of iterations of the algorithm, initialized with $\lambda^0=0$ in $Q_T$. We take $\epsilon=10^{-10}$ as a stopping threshold for the algorithm (the algorithm is stopped as soon as the norm of the residue $g^n$ given here by $Ly^n$ satisfies $\Vert g^n\Vert_{L^2(Q_T)}\leq 10^{-10} \Vert g^0\Vert_{L^2(Q_T)}$). We observe that the number of iterates is sub-linear with respect to $h$, precisely, with respect to the dimension $m_h=card(\{\lambda_h\})$ of the approximated problems. This renders this method very attractive from a numerical point of view. 

Table \ref{tab:ex1_r1_T2} reports the results for $r=1$. We get a slightly better estimate for the norm $\Vert Ly_h\Vert_{L^2(Q_T)}$ but slightly worst estimate for the norm $\Vert y-y_h\Vert_{L^2(Q_T)}$. On the other hand, since $r$ acts as an augmentation parameter, the convergence of the conjugate gradient algorithm is faster. We also check in view of Remark \ref{better_estimate} (specifically, (\ref{eq:llh0}) and (\ref{behavior_deltah_eps0})) that  the convergence of $\lambda_h$ toward zero is much lower for $r=1$ than for $r=h^2$:
\begin{equation}
\nonumber
r=1 : \qquad 
\begin{aligned}
 &\frac{\Vert y-y_h\Vert_{L^2(Q_T)}}{\Vert y\Vert_{L^2(Q_T)}}=\mathcal{O}(h^{1.02}), \quad \frac{\Vert \partial_\nu(y-y_h)\Vert_{L^2(\Gamma_T)}}{\Vert \partial_{\nu} y\Vert_{L^2(\Gamma_T)}}=\mathcal{O}(h^{0.16}), \\
& \Vert \lambda_h\Vert_{L^2(Q_T)}=\mathcal{O}(h^{0.17}), \quad \Vert L y_h\Vert_{L^2(Q_T)}=\mathcal{O}(h^{0.27}).
\end{aligned}
\end{equation}
%
\begin{table}[http]
\centering
\begin{tabular}{|c|ccccc|}
\hline
$h$  &  $7.07\times 10^{-2}$  & $3.53\times 10^{-2}$ & $1.76\times 10^{-2}$ & $8.83\times 10^{-3}$ & $4.42\times 10^{-3}$  \tabularnewline

\hline 
$\frac{\Vert y-y_h\Vert_{L^2(Q_T)}}{\Vert y\Vert_{L^2(Q_T)}}$ & $1.63\times 10^{-2}$ & $6.63\times 10^{-3}$ & $2.78\times 10^{-3}$ & $1.29\times 10^{-3}$ & $5.72\times 10^{-4}$\tabularnewline

$\frac{\Vert \partial_{\nu}(y-y_h)\Vert_{L^2(\Gamma_T)}}{\Vert \partial_{\nu} y\Vert_{L^2(\Gamma_T)}}$ & $7.67\times 10^{-3}$ & $4.95\times 10^{-3}$ & $3.24\times 10^{-3}$ & $2.16\times 10^{-3}$ & $1.48\times 10^{-3}$\tabularnewline

$\Vert L y_h\Vert_{L^2(Q_T)}$ & $0.937$ & $1.204$ & $1.496$ & $1.798$ & $2.135$\tabularnewline

$\Vert \lambda_h\Vert_{L^2(Q_T)}$ & $7.74\times 10^{-3}$ & $3.74\times 10^{-3}$ & $1.72\times 10^{-3}$ & $7.90\times 10^{-4}$  & $3.60\times 10^{-4}$\tabularnewline


card($\{\lambda_h\}$) & $861$ & $3\ 321$ & $13\ 041$ & $51\ 681$ & $205\ 761$\tabularnewline

$\sharp$ \textrm{CG iterates} & $57$ & $103$ & $172$ & $337$ & $591$\tabularnewline

\hline
\end{tabular}
\caption{Example \textbf{EX1} - BFS element - $r=h^2$ - $T=2$.}
\label{tab:ex1_rh2_T2}
\end{table}

\begin{table}[http]
\centering
\begin{tabular}{|c|ccccc|}
\hline
$h$  &  $7.07\times 10^{-2}$  & $3.53\times 10^{-2}$ & $1.76\times 10^{-2}$ & $8.83\times 10^{-3}$ & $4.42\times 10^{-3}$  \tabularnewline

\hline 
$\frac{\Vert y-y_h\Vert_{L^2(Q_T)}}{\Vert y\Vert_{L^2(Q_T)}}$ & $2.25\times 10^{-2}$ & $1.07\times 10^{-2}$ & $5.23\times 10^{-3}$ & $2.62\times 10^{-3}$ & $1.29\times 10^{-3}$\tabularnewline

$\frac{\Vert \partial_{\nu}(y-y_h)\Vert_{L^2(\Gamma_T)}}{\Vert \partial_{\nu} y\Vert_{L^2(\Gamma_T)}}$ & $3.7\times 10^{-2}$ & $3.46\times 10^{-2}$ & $3.14\times 10^{-2}$ & $2.76\times 10^{-2}$ & $2.37\times 10^{-2}$\tabularnewline

$\Vert L y_h\Vert_{L^2(Q_T)}$ & $0.242$ & $0.207$ & $0.171$ & $0.139$ & $0.112$\tabularnewline

$\Vert \lambda_h\Vert_{L^2(Q_T)}$ & $0.147$ & $0.142$ & $0.131$ & $0.116$  & $0.101$\tabularnewline





$\sharp$ \textrm{CG iterates} & $35$ & $60$ & $106$ & $179$ & $312$\tabularnewline

\hline
\end{tabular}
\caption{Example \textbf{EX1} - BFS element - $r=1.$ - $T=2$.}
\label{tab:ex1_r1_T2}
\end{table}

We now discuss the results obtained with the reduced HCT finite element on regular but non-uniform triangulations of the rectangle $\Omega \times (0, T)$. Precisely, we consider $5$ levels of meshes of $Q_T$ described in Table \ref{tab:mesh1D}. For each of these meshes, we compute $h$ as the maximum of the diameters of the triangles composing the triangulation. The coarsest of this meshes is displayed in Figure \ref{fig:yobs_EX}.

\begin{table}[http]
\centering
\begin{tabular}{|c|ccccc|}
\hline
Mesh number  &  1 & 2 & 3 & 4 & 5 \tabularnewline
\hline 
$\sharp$ elements & $688$ & $2\ 752$       & $11\ 008$                     & $44\ 032$                    & $176\ 128$ \tabularnewline
$\sharp$ points      & $375$ & $1\ 437$       & $ 5\  625$                     & $22\ 257$                    & $  88\ 545$ \tabularnewline
$h$                         & $7.62\times10^{-2}$ & $3.81\times10^{-2}$ & $1.91\times10^{-2}$ & $9.53\times 10^{-3}$ &  $4.77\times 10^{-3}$ \tabularnewline
\hline
\end{tabular}
\caption{Example \textbf{EX1} - HCT element  - Information concerning the meshes of the domain $Q_T$.}
\label{tab:mesh1D}
\end{table}

Table \ref{tab:1D} collects the numerical results on the reconstruction of the solution $y$ from the observation $y_{\nu, obs}$ obtained with again $r=h^2$. We observe a slightly super linear convergence for the variable $y_h$ and $\lambda_h$:
\begin{equation}
\nonumber
\frac{\|y - y_h\|_{L^2(Q_T)}}{\|y\|_{L^2(Q_T)}} = \mathcal{O}(h^{1.22}), \qquad \|\lambda_h\|_{L^2(Q_T)} = \mathcal{O}(h^{1.04}).
\end{equation}

\begin{table}[http]
\centering
\begin{tabular}{|c|ccccc|}
\hline
$h$ &
      $7.62\times10^{-2}$ & $3.81\times10^{-2}$ & $1.91\times10^{-2}$ & $9.53\times 10^{-3}$ & $4.77\times 10^{-3}$ \tabularnewline
\hline 
$\frac{\|y- y_h\|_{L^2(Q_T)}}{\|y\|_{L^2(Q_T)}}$ & $3.67 \times 10^{-2}$ & $1.35 \times 10^{-2}$ & $5.99 \times 10^{-3}$ & $2.63 \times 10^{-3}$ & $1.22 \times 10^{-3}$ \tabularnewline
$\|\lambda_h\|_{L^2(Q_T)}$  & $2.12 \times 10^{-2}$ & $1.08 \times 10^{-2}$ & $5.45 \times 10^{-3}$ & $2.53 \times 10^{-3}$ & $1.18 \times 10^{-3}$ \tabularnewline
$\|L y_h\|_{L^2(Q_T)}$ & $1.88 $ & $2.51$ & $3.26$ & $4.13$ & $5.13$ \tabularnewline
$\kappa_h$ & $2.15\times 10^6$ & $1.11 \times 10^7$ & $1.03 \times 10^{8}$ & $8.67 \times 10^{8}$ & $6.94\times 10^{9}$ \tabularnewline
\hline
\end{tabular}
\caption{Example
      \textbf{EX1} - HCT element - $r=h^2$ - $T=2$.}
\label{tab:1D}
\end{table}

\begin{figure}[ht!]
\centering
\begin{tabular}{cc}
\includegraphics[width=0.46\textwidth]{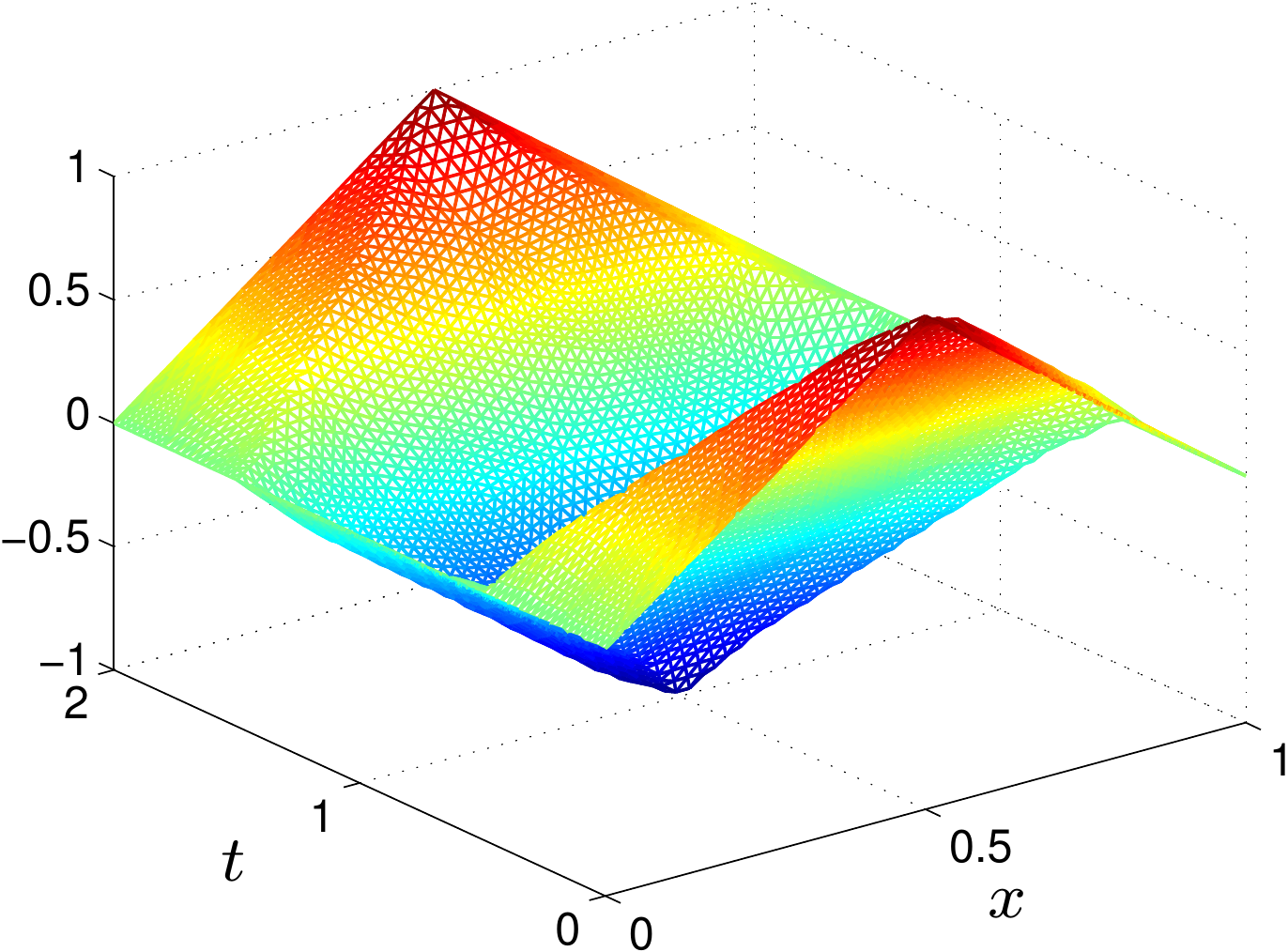} & \includegraphics[width=0.46\textwidth]{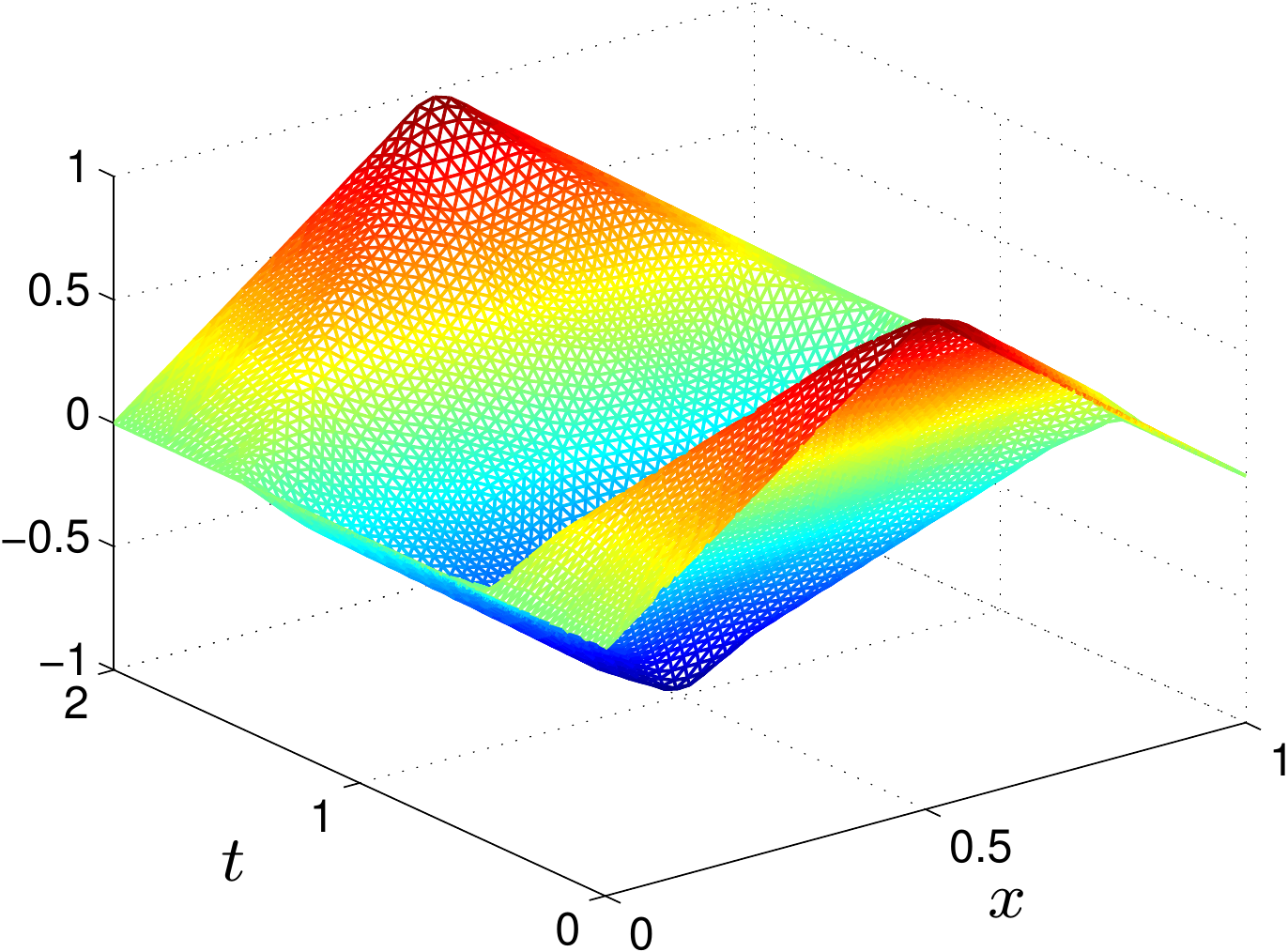} 
\end{tabular}
\caption{Example \textbf{EX1}. Exact solution $y$ and approximated solution $y_h$ on the mesh $\sharp \ 3$.}
\label{fig:EX}
\end{figure}

\begin{figure}[ht!]
\centering
\includegraphics[width=0.55\textwidth]{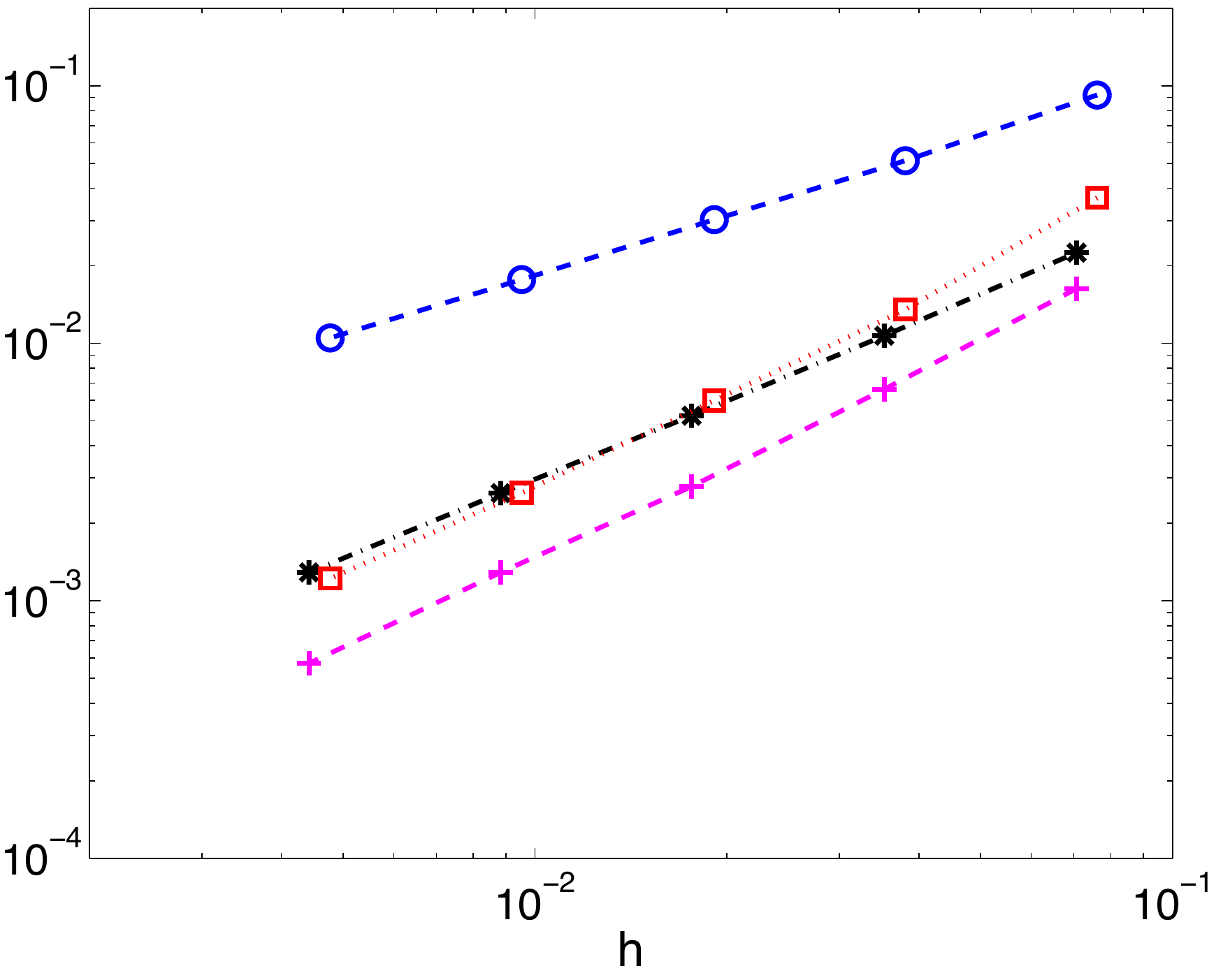}
\caption{Example \textbf{EX1} -$T=2$ - Relative error $\|y - y_h\|_{L^2(Q_T)}/\|y\|_{L^2(Q_T)}$ w.r.t. $h$ for the BFS element with $r=h^2$ ($+$) and $r=1$ ($\star$), the HCT element with $r=h^2$ ($\square$) and $r=1$ ($\circ$)}
\label{fig:EX_Err}
\end{figure}

Figure \ref{fig:EX} depicts the exact solution $y$ computed by (\ref{yobs_EX}) and its approximation $y_h$ (computed with the mesh $\sharp \ 3$). Figure \ref{fig:EX_Err} represents the relative error $\|y - y_h\|_{L^2(Q_T)}/\|y\|_{L^2(Q_T)}$ reported also in Table \ref{tab:1D}. 
Tables \ref{tab:1D}, \ref{tab:EX_adapt} and \ref{tab:1D_stab} also report the value of the condition number $\kappa_h$ of the matrix $A_{r,h}$. In both situations, as it is usual for elliptic problems, $\kappa_h$ behaves quadratically with respect to $h^{-1}$. 

\begin{figure}[ht!]
\hspace{-1.76cm}
\centering
\begin{tabular}{ccc}
\includegraphics[width=0.29\textwidth]{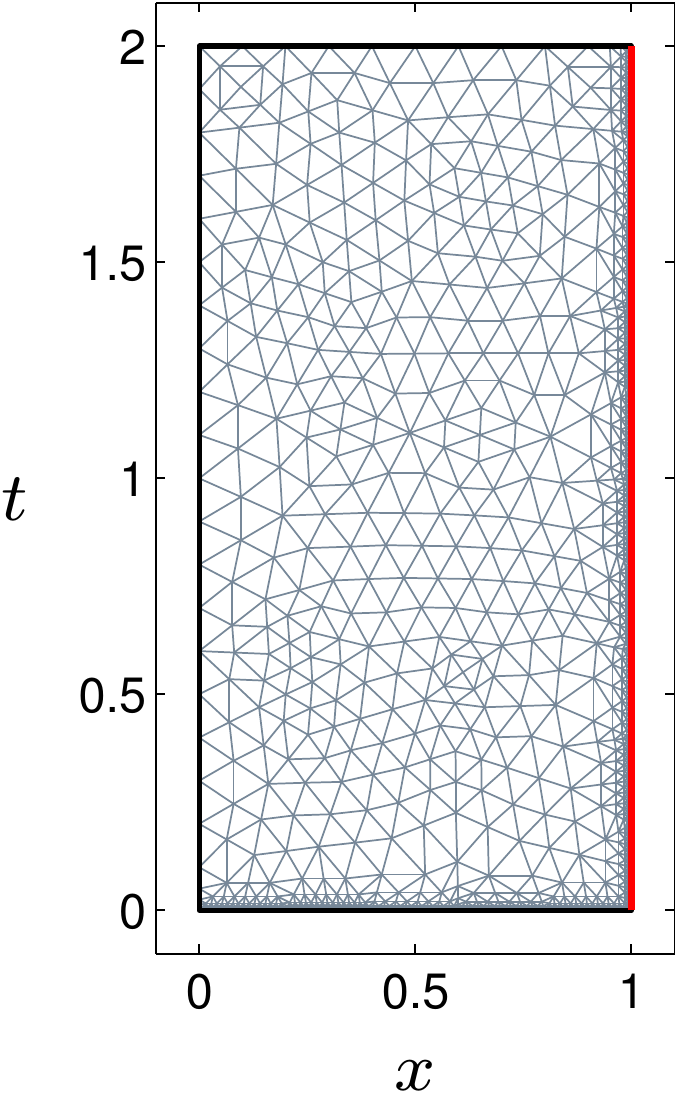} & \qquad&
\includegraphics[width=0.29\textwidth]{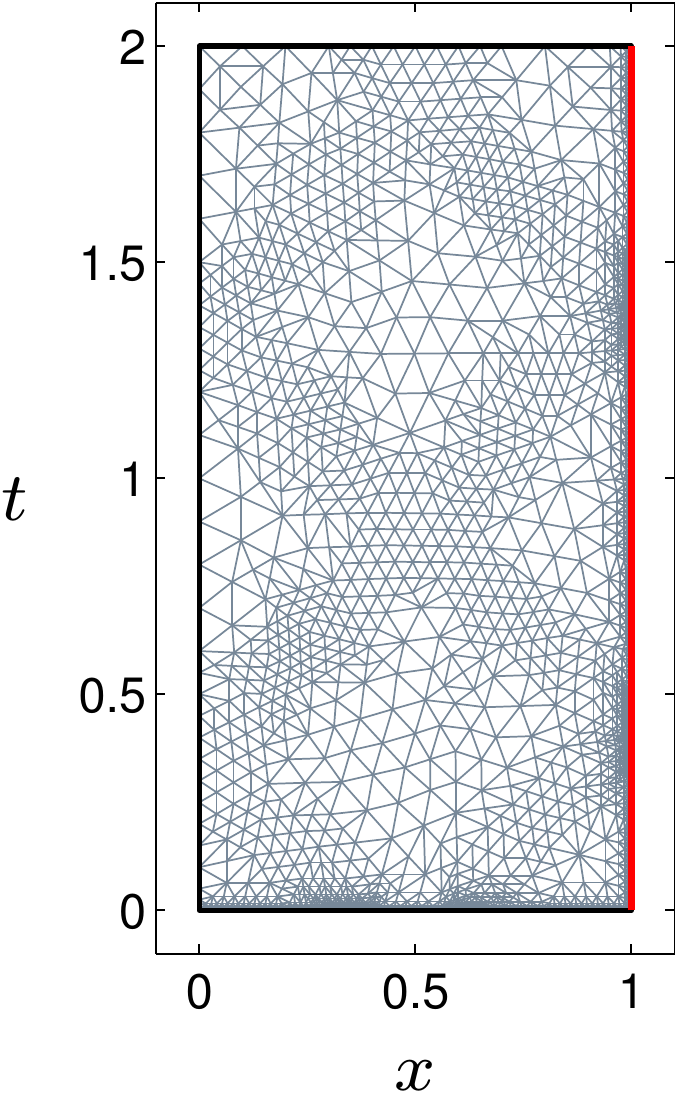} \\
\includegraphics[width=0.29\textwidth]{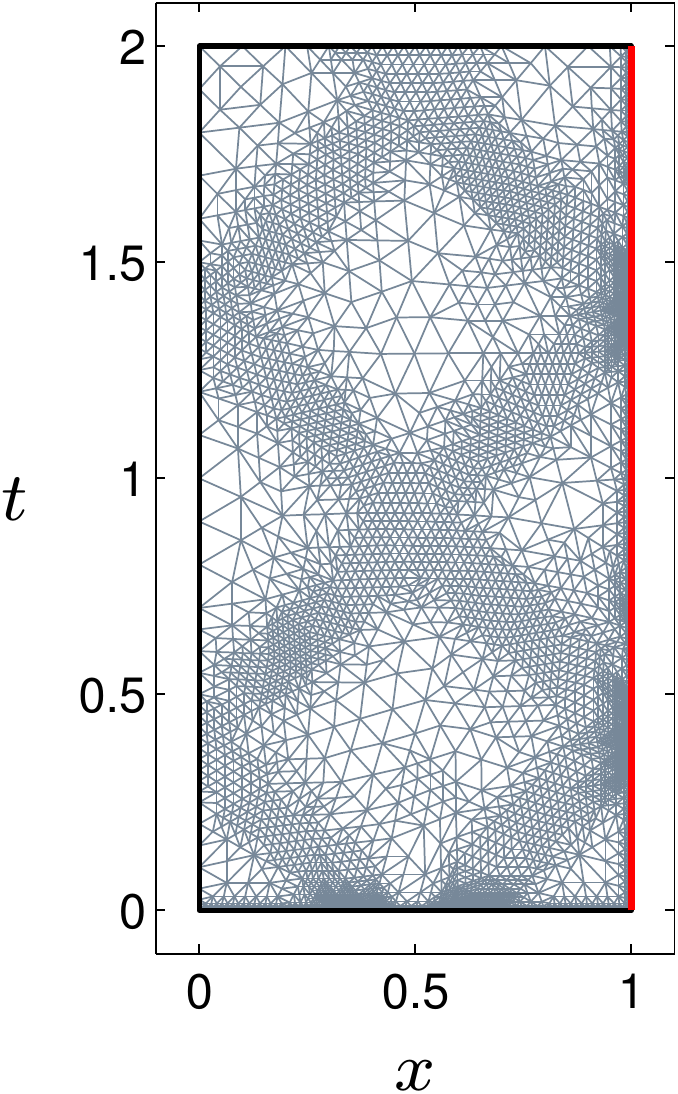} & \qquad &
\includegraphics[width=0.29\textwidth]{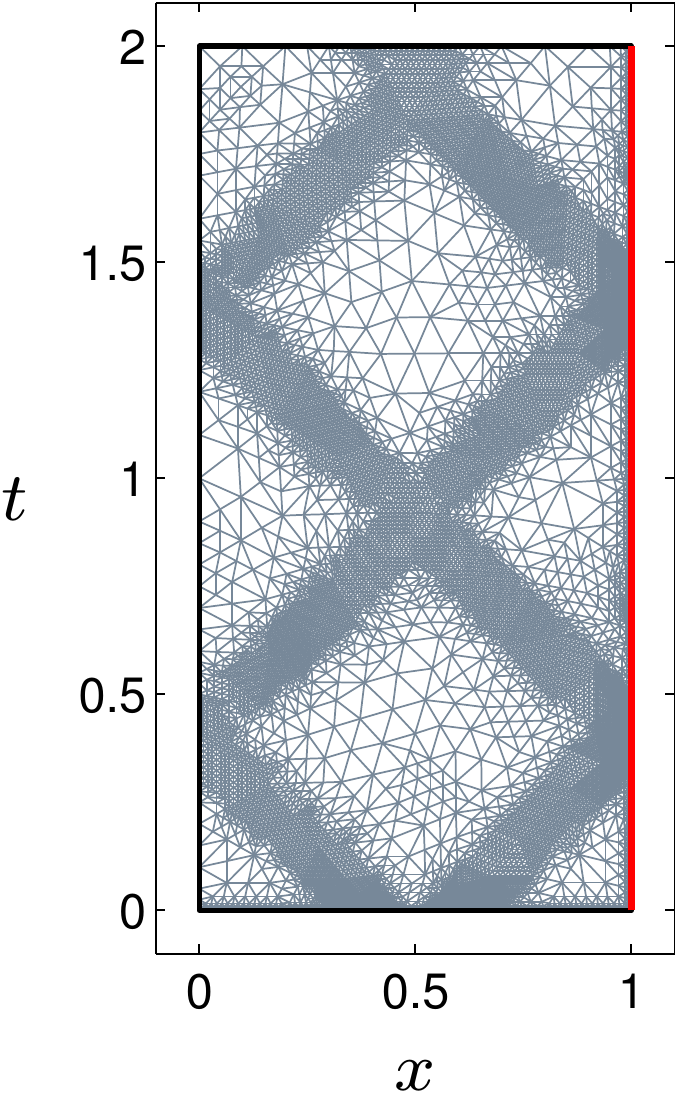}
\end{tabular}
\caption{Example \textbf{EX1} - reduced HCT finite element: Iterative refinement of the triangular mesh over $Q_T$ with respect to the variable $y_h$}
\label{fig:EX_adapt}
\end{figure}

We also emphasize that this variational method which requires a finite element discretization of the time-space $Q_T$ is particularly well-adapted to mesh optimization. Still for the example \textbf{EX1}, Figure \ref{fig:EX_adapt} depicts a sequence of four meshes of $Q_T=(0,1)\times (0,T)$: the sequence is initialized with the coarsest mesh described in Table \ref{tab:mesh1D} which is locally refined near boundary $\Gamma_T$ (where the observation $y_{\nu, obs}$ is localized) and near $\Omega \times \{0\}$ (for a better representation of the initial data). The three other meshes are successively obtained by local refinement based on the norm of the gradient of $y_h$ on each triangle of $\mathcal{T}_h$. As expected, the refinement is concentrated around the lines of singularity of $y_h$ traveling in $Q_T$, generated by the singularity of the initial position $y_0$. Some information concerning these meshes and the approximation errors obtained where this mesh adaptation strategy is employed are reported in Table \ref{tab:EX_adapt}.

\begin{table}[ht!]
\centering
\begin{tabular}{|c|cccc|}
\hline
Mesh
      number & 1 & 2 & 3 & 4 \tabularnewline
\hline
$\sharp$
      elements & $3\ 452$ & $6\ 262$ & $16\ 440$ & $49\ 499$ \tabularnewline
$\sharp$
      points & $1\ 986$ & $3 \ 442$ & $8\ 629$ & $25\ 318$\tabularnewline
$\frac{\|y
      - y_h\|_{L^2(Q_T)}}{\|y\|_{L^2(Q_T)}}$ & $1.35 \times 10^{-2}$ & $8.02 \times 10^{-3}$ & $5.5 \times 10^{-3}$ & $4.51
      \times 10^{-3}$ \tabularnewline
$\|\lambda_h\|_{L^2(Q_T)}$
      & $2.11 \times 10^{-2}$ & $1.38 \times 10^{-2}$ & $8.57 \times 10^{-3}$ & $4.56 \times 10^{-3}$\tabularnewline
$\|L
      y_h\|_{L^2(Q_T)}$& $7.03$& $6.98$& $7.82$ & $9.01$
      \tabularnewline
$\kappa_h$ & $7.9 \times 10^8$ & $3.05 \times 10^9$ & $4.38 \times 10^{10}$ & $6.19 \times 10^{11}$\tabularnewline
\hline
\end{tabular}
\caption{Example \textbf{EX1}- HCT element -  $r=h^2$ - $T=2$ - Iterative refinement.}
\label{tab:EX_adapt}
\end{table}

We end this section with some numerical results for the stabilized mixed formulation (\ref{mf_stab_h}). The main difference is that the multiplier $\lambda$ is now approximated in a richer space  (see (\ref{choice_lambdah_alpha})) leading to larger linear systems. Table \ref{tab:1D_stab} considers the case of the example \textbf{EX1} with $T=2$ and $\alpha=1/2$. In order to compare with the formulation (\ref{eq:mfh}), we take again $r=h^2$. We observe the convergence w.r.t. $h$ and obtain similar rates and constants to the ones in Table \ref{tab:1D}: in particular, we have
\begin{equation}
\frac{\|y - y_h\|_{L^2(Q_T)}}{\|y\|_{L^2(Q_T)}} = \mathcal{O}(h^{1.21}), \qquad \|\lambda_h\|_{L^2(Q_T)} = \mathcal{O}(h^{1.04}). \nonumber
\end{equation}
Finally, we also check - in contrast with the mixed formulation (\ref{eq:mf}) - that the positive parameter $r$ does not affect the numerical results. 

\begin{table}[!ht]
\centering
\begin{tabular}{|c|ccccc|}
\hline
$h$ & $7.62\times10^{-2}$ & $3.81\times10^{-2}$ &  $1.91\times10^{-2}$ & $9.53\times 10^{-3}$ & $4.77\times 10^{-3}$ \tabularnewline
\hline 
$\frac{\|y- y_h\|_{L^2(Q_T)}}{\|y\|_{L^2(Q_T)}}$ & $4.1 \times 10^{-2}$ & $1.55 \times 10^{-2}$ & $6.88 \times 10^{-3}$ & $3.03 \times 10^{-3}$ & $1.39 \times 10^{-3}$ \tabularnewline
$\|\lambda_h\|_{L^2(Q_T)}$ & $6.96 \times 10^{-3}$ & $3.67 \times 10^{-3}$ & $1.71 \times 10^{-3}$ & $8.42 \times 10^{-4}$ & $3.98 \times 10^{-4}$ \tabularnewline
$\|Ly_h\|_{L^2(Q_T)}$ & $1.47 $ & $1.86$ & $2.28$ & $2.76$ & $3.23$ \tabularnewline
$\kappa_h$ & $5.05 \times 10^8$ & $3.01 \times 10^8$ & $2.59 \times 10^{9}$ & $1.82 \times 10^{10}$ & $1.31 \times 10^{11}$ \tabularnewline
\hline
\end{tabular}
\caption{Example
      \textbf{EX1} - HCT element - $r=h^2$ - $T=2$ - Stabilized mixed formulation (\ref{mf_stab_h}).}
\label{tab:1D_stab}
\end{table}

\subsection{Reconstruction of the solution - Two dimensional case ($N = 2$)}

In this section we illustrate the method introduced in Section \ref{recovering_y} on a two-dimensional example.  The procedure is similar but a bit more involved on a computational point of view, since $Q_T$ is now a subset of $\mathbb{R}^3$. We take again $c:=1$ in $\Omega$ and $d:=0$ in $Q_T$.

\begin{figure}[ht]
\hspace{-0.6cm}
\centering
\begin{tabular}{m{.4\textwidth}cm{.4\textwidth}}
\includegraphics[width=0.4\textwidth]{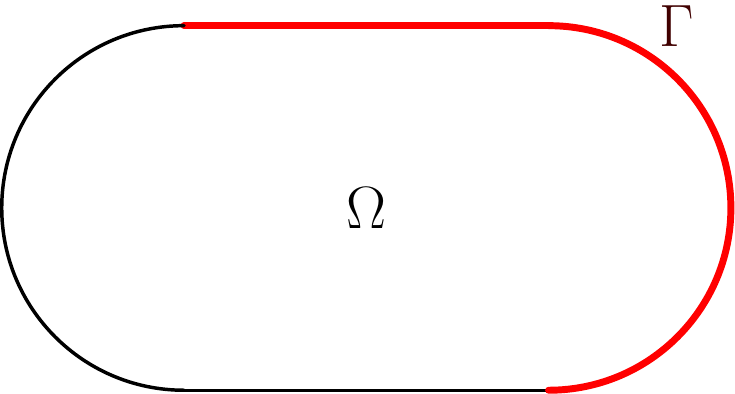} & \qquad\qquad &
\includegraphics[width=0.35\textwidth]{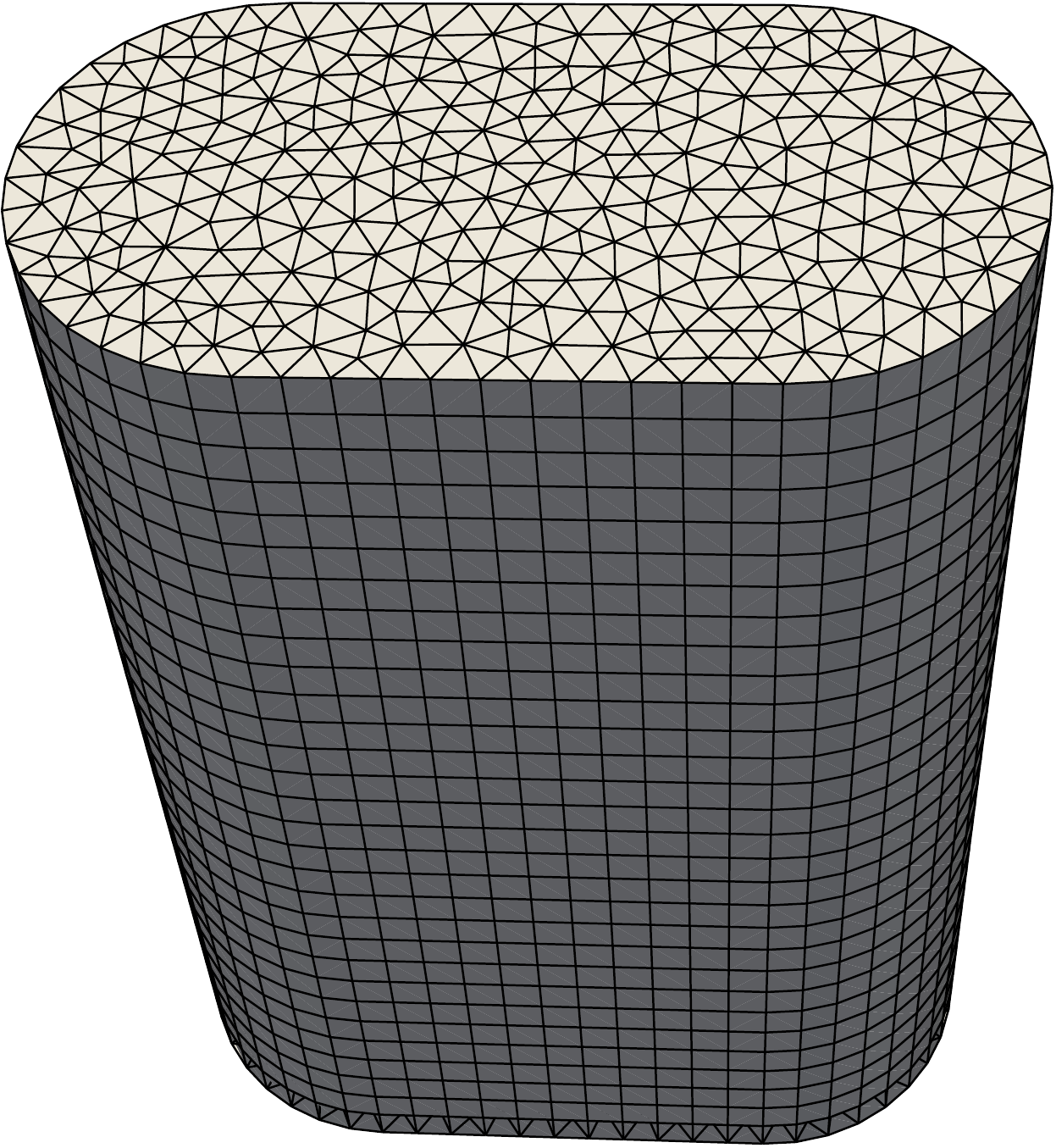} 
\end{tabular}
\caption{Bunimovich's stadium and the subset $\Gamma$ of $\partial \Omega$ on which the observation is available (Left). Example of mesh of the domain $Q_T$ (Right).}
\label{fig:stadium}
\end{figure}

In order to approximate the mixed-formulation (\ref{eq:mf}), we consider a mesh $\mathcal{T}_h$ of the domain $Q_T = \Omega \times (0,T)$ formed by triangular prisms. This mesh is obtained by extrapolating along the time axis a triangulation of the spatial domain $\Omega$. For all the simulations considered in this section, $\Omega$ is the so called Bunimovich's stadium (see \cite{Bunimovich}) and $T=3$. Figure \ref{fig:stadium}-Left displays the domain $\Omega$ and the part $\Gamma$ of the boundary on which the observation is available for $t \in (0, T)$ while Figure \ref{fig:stadium}-right displays an example of mesh of domain $Q_T$.

Let $Z_h$ be the finite dimensional space defined as follows
\begin{equation}
\label{eq:Yh}
Z_h = \left\{ 
\begin{array}{l|l}
\ph_h = \psi(x_1, x_2) \theta(t) \in C^1(Q_T) & \psi|_{K_{x_1x_2}} \in \mathbb{P}(K_{xy}), \theta|_{K_t} \in \mathbb{Q}(K_t) \\ 
 \ph_h = 0 \text{ on } \Sigma_T & \text {for every } K = K_{x_1x_2} \times K_t \in \mathcal{T}_h.
\end{array}
\right\},
\end{equation}
$\mathbb{P}(K_{x_1x_2})$ is the space of functions corresponding to the reduced {\em Hsieh-Clough-Tocher} (HCT) $C^1$-element recalled in Section 4.1.1; $\mathbb{Q}(K_t)$ is a space of degree three polynomials on the interval $K_t$ of the form $[t_j,t_{j+1}]$ defined uniquely by their value and the value of their first derivative at the point $t_j$ and $t_{j+1}$. In other words, $Y_h$ is the finite element space obtained as a tensorial product between the reduced HCT finite element and the cubic Hermite finite element. 
We check that on each element $K=K_{x_1x_2}\times K_t$, the function $\ph_h$ is
determined uniquely in term of the values of $\Sigma_K:=\{\ph(a_i),\ph_{x_1}(a_i),\ph_{x_2}(a_i),\ph_t(a_i),\ph_{x_1,t}(a_i),\ph_{x_2,t}(a_i), i=1,\cdots,6\}$ at the six nodes $a_i$ of $K$. Therefore, $\dim \Sigma_K=36$.

We consider meshes formed by triangular prisms of the domain $Q_T = \Omega \times (0, T)$. An example of such a mesh associated to the domain $Q_T$ is displayed in Figure \ref{fig:stadium} right. This mesh is composed by $10 \ 261$ nodes distributed in $18 \ 060$ prismatic elements (this mesh corresponds to the mesh number 2 described in Table \ref{tab:mesh_h}).

We consider three levels of meshes of the domain $Q_T$ formed by the number of prisms and containing the number of nodes reported in Table \ref{tab:mesh_h}.

\begin{table}[ht]
\centering
\begin{tabular}{|c|ccc|}
\hline
Mesh number & 1 & 2 & 3\\
\hline 
Number of elements & 1  860 & 18 060 & 158 280 \\
Number of nodes & 1 216 & 10 261 & 84 241 \\
$\Delta x$ & $1.82\times 10^{-1}$ & $8.2\times 10^{-2}$ & $3.95\times 10^{-2}$ \\
$\Delta t$ (Height of elements) & $0.2$ & $0.1$ & $0.05$ \\
$h$ & $2.7\times 10^{-1}$ & $1.29\times 10^{-1}$ &  $6.37\times 10^{-2}$\\
\hline
\end{tabular}
\caption{Characteristics of the three meshes associated with $Q_T$.}
\label{tab:mesh_h}
\end{table}

Comparing to the one dimensional situation described in Section \ref{sec:1d}, the eigenfunctions and eigenvectors of the Dirichlet Laplace operator defined on $\Omega$ are not explicitly available. Consequently, from a given set of initial data, we define as the "exact" solution and note $\overline{y}$ the solution obtained numerically with a very fine discretization, from which we can extract an observation on $\Gamma_T$. Precisely, we solve the hyperbolic equation (\ref{eq:wave}) using a standard time-marching method: we employ a {\em HCT} finite elements method in space coupled with a Newmark unconditionally stable scheme for the time discretization. 

Here,  we solve the hyperbolic equation on the spatial mesh which was extrapolated in time in order to obtain the mesh number 3 of $Q_T$. This two-dimensional mesh contains $1 \ 381$ nodes and $2\ 638$ triangles and corresponds to the value $\Delta x\approx 3.95\times 10^{-2}$. As for the time discretization, we use the value $\Delta t = 10^{-2}$. We denote $\overline{y}$ the solution obtained in this way for the initial data $(y_0, y_1) \in H_0^1(\Omega) \times L^2(\Omega)$ given by

\begin{equation}
(\textbf{EX2})\hspace{1cm}
\label{eq:idheart}
\left\{ 
\begin{array}{ll}
-\Delta y_0 = 10,  &\quad \text{in } \Omega \\
y_0 = 0, & \quad \text{on } \partial \Omega,
\end{array}
\right.
\qquad y_1 =0 \quad \text{in } \Omega.
\end{equation}

From $\overline{y}$, we then generate the observation $y_{\nu, obs}$ as the restriction of $\partial_{\nu}\overline{y}$ on $\Gamma_T$. The geometry used here allows to compute easily 
the normal derivative of $\overline{y}$ on the boundary. Finally, from this observation we reconstruct $y_h$ as the solution of the mixed formulation (\ref{eq:mfeps}) on each of the three meshes described in Table \ref{tab:mesh_h}. For this simulations we take the augmentation parameter $r = h^3$.  Table \ref{tab:h} displays some norms of $y_h$ and $\lambda_h$ obtained for the three meshes and illustrates again the convergence of the method. 

\begin{table}[ht!]
\centering
\begin{tabular}{|c|ccc|}
\hline
Mesh number & 1 & 2 & 3\\
\hline
$\frac{\|\overline{y} - y_h\|_{L^2(Q_T)}}{\|\overline{y}\|_{L^2(Q_T)}}$ & $3.75 \times 10^{-2}$ & $1.53 \times 10^{-2}$ & $1.39\times 10^{-2}$ \\
$\|Ly_h\|_{L^2(Q_T)}$ & $2.07$ & $1.39$ & $1.09$\\
$\|\lambda_h\|_{L^2(Q_T)}$ & $9.3 \times 10^{-6}$ & $4.89 \times 10^{-6}$ & $4.57 \times 10^{-6}$ \\
$\| \overline y(\cdot, 0) - y_h(\cdot, 0) \|_{L^2(\Omega)}$ & $5.9 \times 10^{-2}$ & $2.23 \times 10^{-2}$ & $1.08 \times 10^{-2}$ \\
\hline
\end{tabular}
\caption{Initial data $(y_0,y_1)$ given by (\ref{eq:idheart}).}
\label{tab:h}
\end{table}

\begin{figure}[ht]
\centering
\begin{tabular}{cc}
\includegraphics[width=0.45\textwidth]{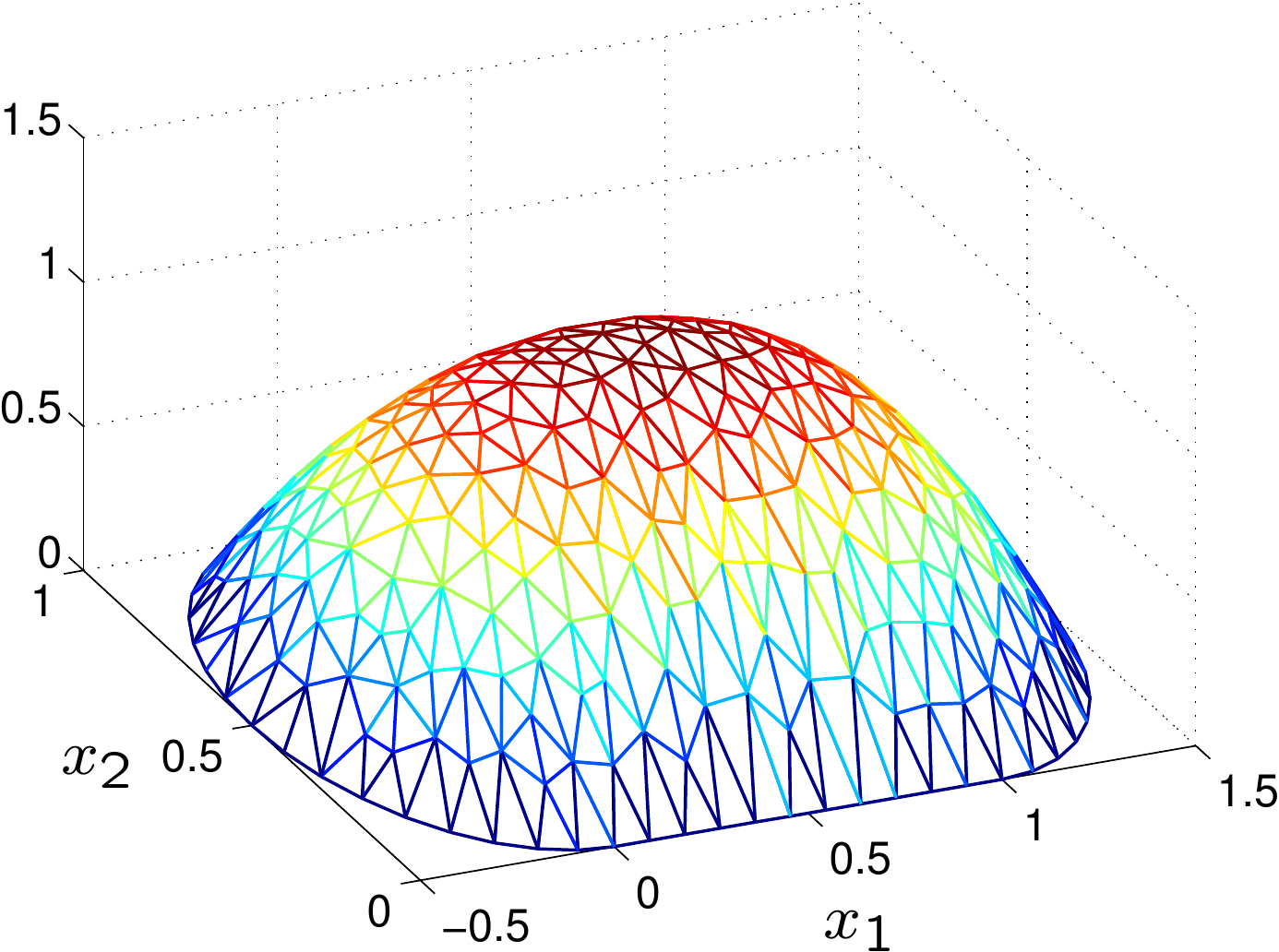} & \includegraphics[width=0.45\textwidth]{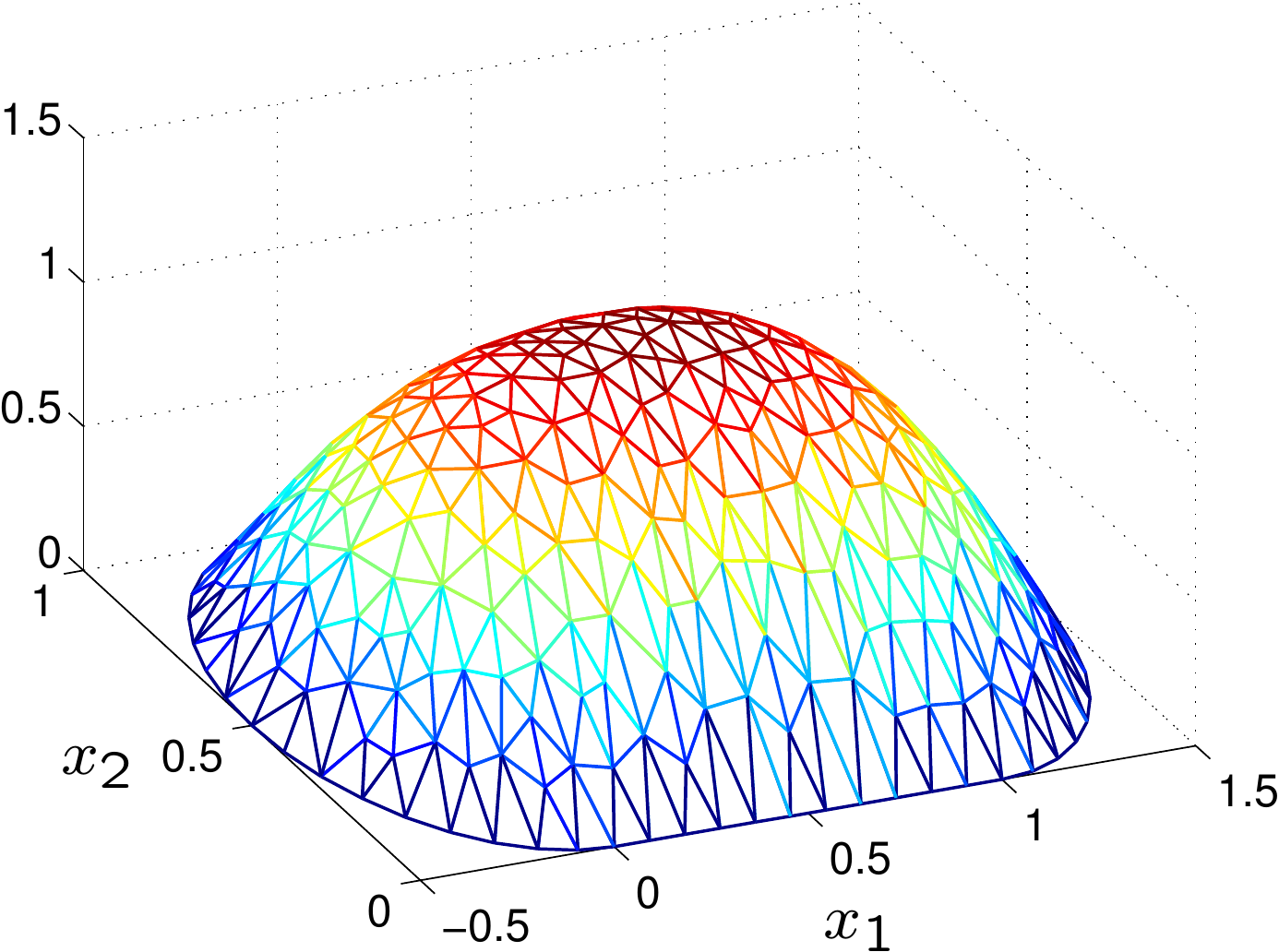}  
\end{tabular}
\caption{Initial data $y_0$ given by (\ref{eq:idheart}) (Left). Reconstructed initial data $y_h(\cdot, 0)$ (Right).}
\label{fig:idheart}
\end{figure}

Figure \ref{fig:idheart}-Left displays the initial position $y_0$ of (\ref{eq:idheart}) while Figure \ref{fig:idheart}-Right displays the initial position $y_h(\cdot, 0)$ corresponding to restriction at time $t=0$ of the solution $y_h$ of the inverse problem. Figures correspond to the mesh number 2.. The errors between these two functions are given in the last row of Table \ref{tab:h}.

\subsection{Reconstruction of the solution and the source $\mu\in H^{-1}(\Omega)$ - One dimensional case $N=1$}
We now consider the reconstruction of both the state and the source from a partial observation, as discussed in Section \ref{recovering_y_mu}.
In order to construct explicit solution, we first recall that the solution of (\ref{eq:wave}) with zero initial condition can be expanded as follows : 
$$
\nonumber
\left\{
\begin{aligned}
& y(x,t)=\sum_{p>0} b_p(t)\sin(p\pi x) \\
& b_p(t):=\frac{1}{p\pi}\int_0^t \sin(p\pi (t-s))f_p(s)ds, \quad f_p(s):=2\sigma(s)\int_{\Omega} \sin(p\pi x) \mu(x) dx.
\end{aligned}
\right.
$$
In the following examples, we take $\sigma(t)=1+t$ defined on $[0,T]$, $T=2$ and  $\Gamma_T=\{1\}\times (0,T)$.

\begin{figure}[ht]
\centering
\begin{tabular}{ccc}
\includegraphics[width=0.46\textwidth]{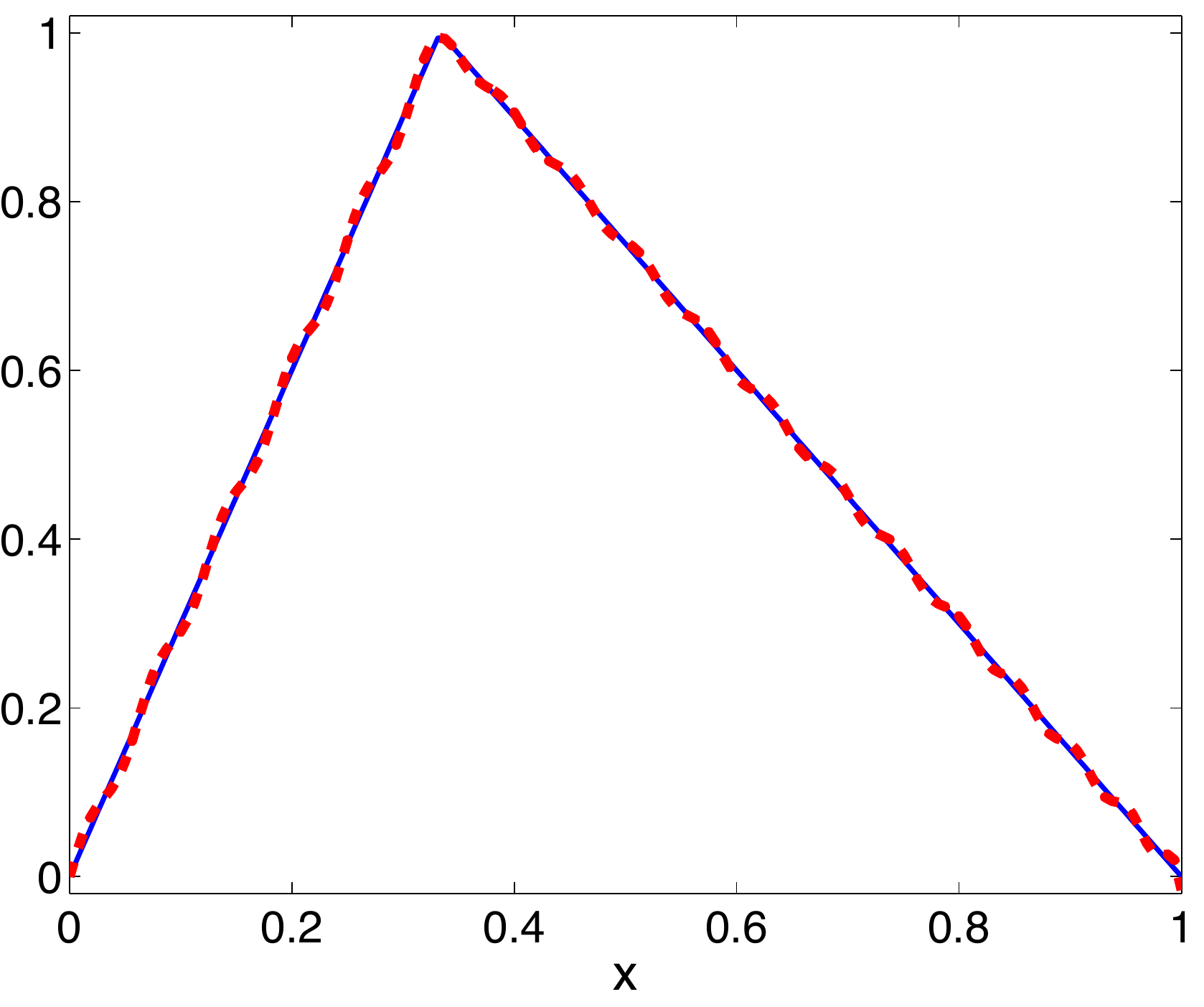} & \qquad & 
\includegraphics[width=0.46\textwidth]{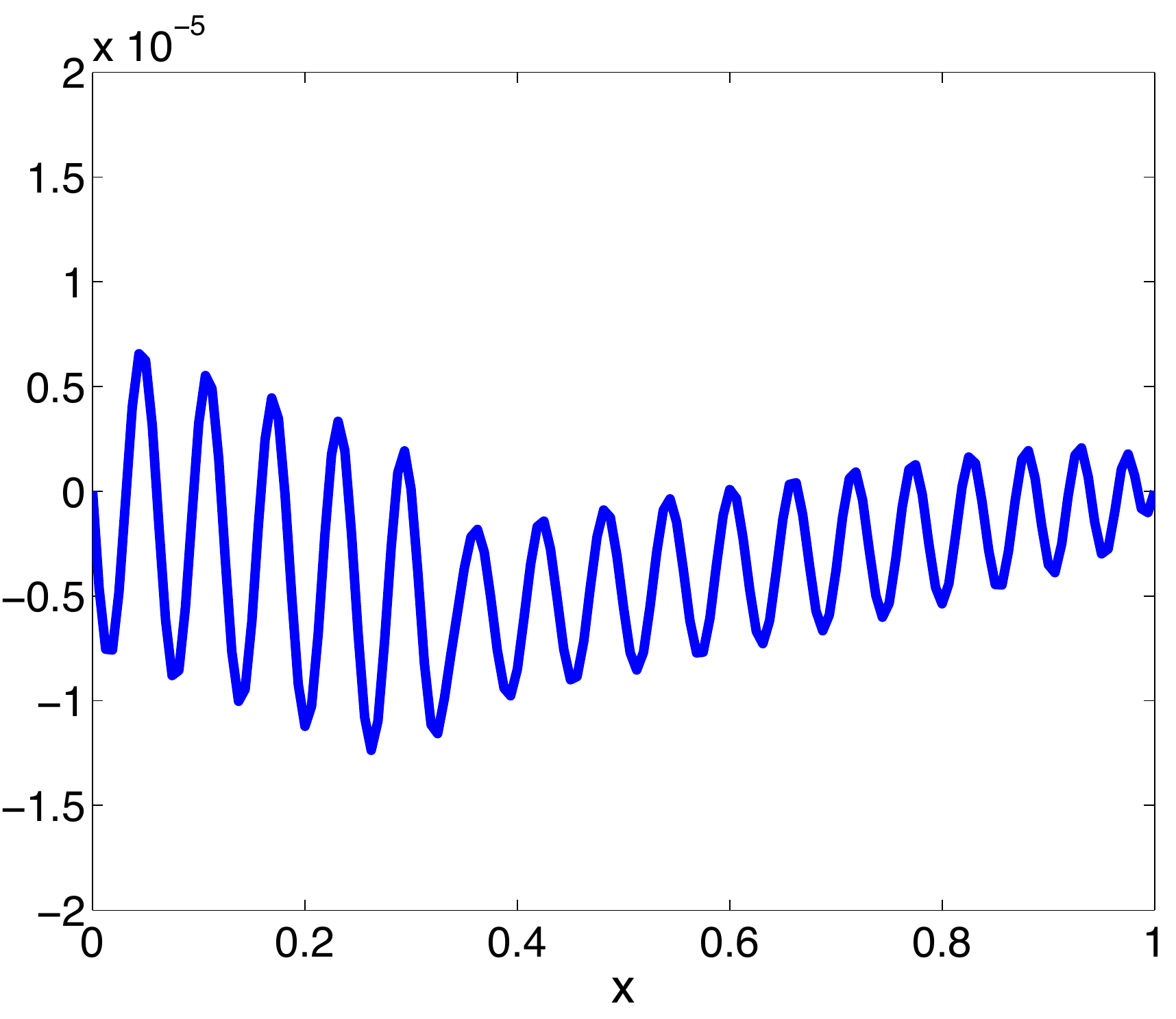} 
\end{tabular}
\caption{\textbf{EX3} -\textbf{Left}: Function $\mu$ (full line) and its approximation $\mu_h$ (dotted line) along $\Omega$;   \textbf{Right}:   $\frac{(-\Delta)^{-1}(\mu-\mu_h)}{\Vert(-\Delta)^{-1}\mu\Vert_{H^1_0(\Omega)}}$ along $\Omega$.}
\label{fig:ex3_mu}
\end{figure}

We now report the resolution of the discrete mixed formulation (\ref{eq:mfepsh}).
 We first consider a rather smooth case, with $\mu\in H^{1}(\Omega)$ given by 
\begin{equation}
\nonumber
(\textbf{EX3})\qquad   \mu(x)=  \frac{x}{\theta}\, 1_{[0,\theta]}(x) + \frac{(1-x)}{1-\theta}\, 1_{[\theta,1]}(x), \quad \theta=1/3.
\end{equation}
Table  \ref{tab:ex3_rh4_T2} reports the main norms with respect to $h$. Concerning the augmentation parameter $r$, we use $r=h^4$ which leads to slightly better approximation of the function $\mu$ than $r=h^2$. We check the convergence of the approximations $(y_h,\mu_h)$ as $h$ tends to $0$.  In particular, we get 
\begin{equation}\nonumber
\frac{\Vert y-y_h\Vert_{L^2(Q_T)}}{\Vert y\Vert_{L^2(Q_T)}}=\mathcal{O}(h^{1.9}), \qquad \frac{\Vert \mu-\mu_h \Vert_{H^{-1}(\Omega)} }{ \Vert \mu\Vert_{H^{-1}(\Omega)}}=\mathcal{O}(h^{1.4}) .
\end{equation}
For instance,  for $h=8.83\times 10^{-3}$ (precisely, $\Delta x=\Delta t=1/160$) leading to a linear system with $258\ 566$ unknowns, we get a relative error for $\mu_h$ of the order of $10^{-4}$. This allows a very good reconstruction of the corresponding solution $y$. Figure \ref{fig:ex3_mu}-Left depicts the function $\mu$ and its corresponding approximation $\mu_h$. Figure \ref{fig:ex3_mu}-Right  depicts along $\Omega$ the function $(-\Delta)^{-1}(\mu-\mu_h)/\Vert (-\Delta)^{-1} \mu\Vert_{H^1_0(\Omega)}$ of magnitude $10^{-5}$. Here, $\Delta$ denotes the Dirichlet Laplacian.  








\begin{table}[http]
\centering
\begin{tabular}{|c|ccccc|}
\hline
$h$  &  $7.07\times 10^{-2}$  & $3.53\times 10^{-2}$ & $1.72\times 10^{-2}$ & $8.83\times 10^{-3}$ & $7.07\times 10^{-3}$  \tabularnewline

\hline 

$\frac{\Vert y-y_h\Vert_{L^2(Q_T)}}{\Vert y\Vert_{L^2(Q_T)}}$  & $1.72\times 10^{-3}$ & $5.06\times 10^{-4}$  & $1.28\times 10^{-4}$ &  $3.45\times 10^{-4}$ &  $2.14\times 10^{-5}$\tabularnewline

$\frac{\Vert \mu-\mu_h\Vert_{H^{-1}(\Omega)}}{\Vert \mu\Vert_{H^{-1}(\Omega)}}$ & $5.9\times 10^{-3}$ & $1.63\times 10^{-3}$ & $8.3\times 10^{-4}$  & $3.79\times 10^{-4}$ &   $1.68\times 10^{-4}$\tabularnewline

$\frac{\Vert \partial_{\nu}(y-y_h)\Vert_{L^2(\Gamma_T)}}{\Vert \partial_{\nu} y\Vert_{L^2(\Gamma_T)}}$ & $6.5\times 10^{-4}$  & $2.03\times 10^{-4}$ & $5.17\times 10^{-5}$ & $1.38\times 10^{-5}$ &  $8.89\times 10^{-6}$\tabularnewline

$\Vert \lambda_h\Vert_{L^2(Q_T)}$ & $7.44\times 10^{-2}$ & $5.26\times 10^{-2}$ & $3.68\times 10^{-2}$  & $2.63\times 10^{-2}$  & $2.35\times 10^{-2}$ \tabularnewline

\hline
\end{tabular}
\caption{Example \textbf{EX3} - BFS element - $r=h^4$ - $T=2$.}
\label{tab:ex3_rh4_T2}
\end{table}

\begin{figure}[ht]
\centering
\begin{tabular}{ccc}
\includegraphics[width=0.46\textwidth]{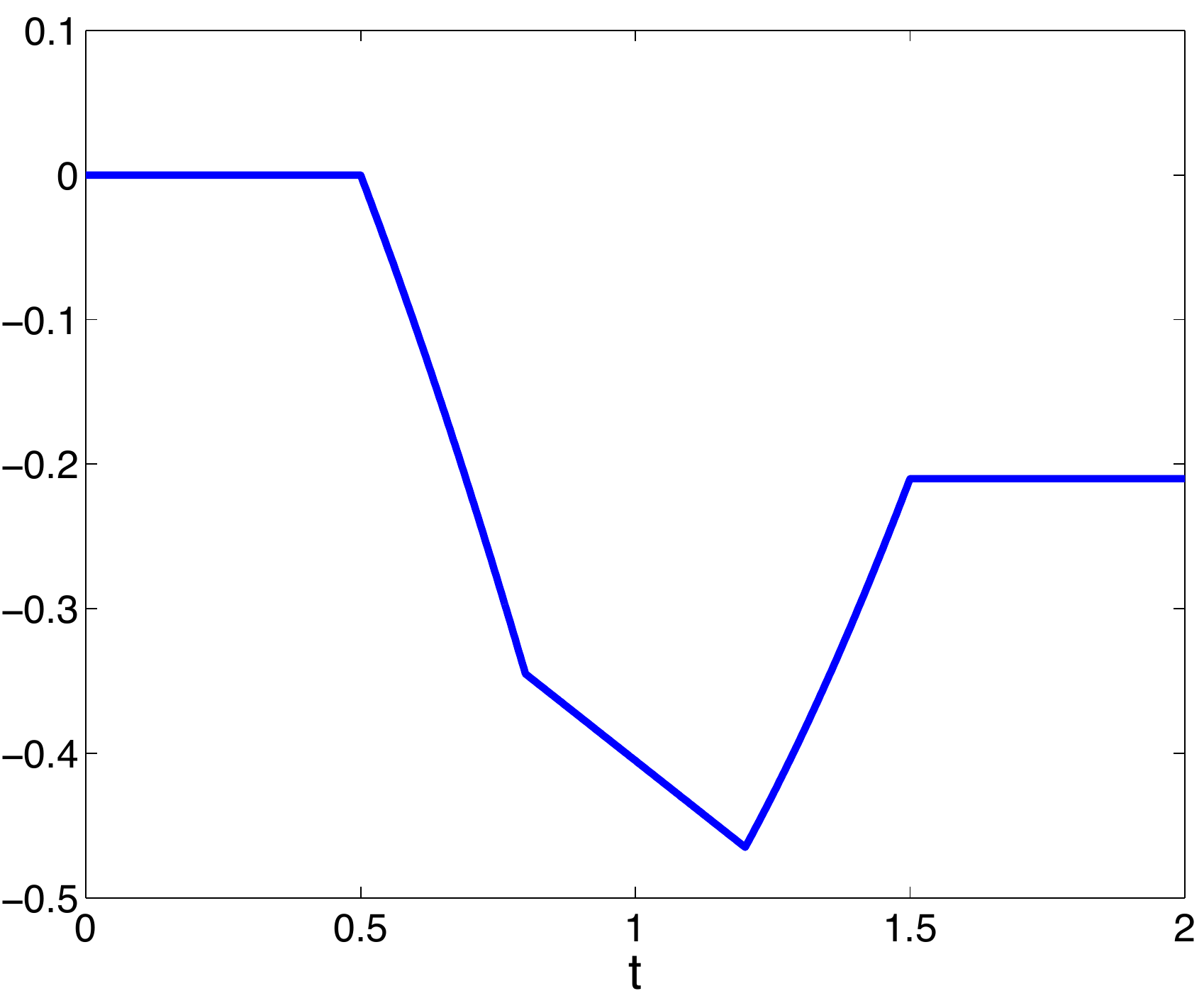} & \qquad & 
\includegraphics[width=0.46\textwidth]{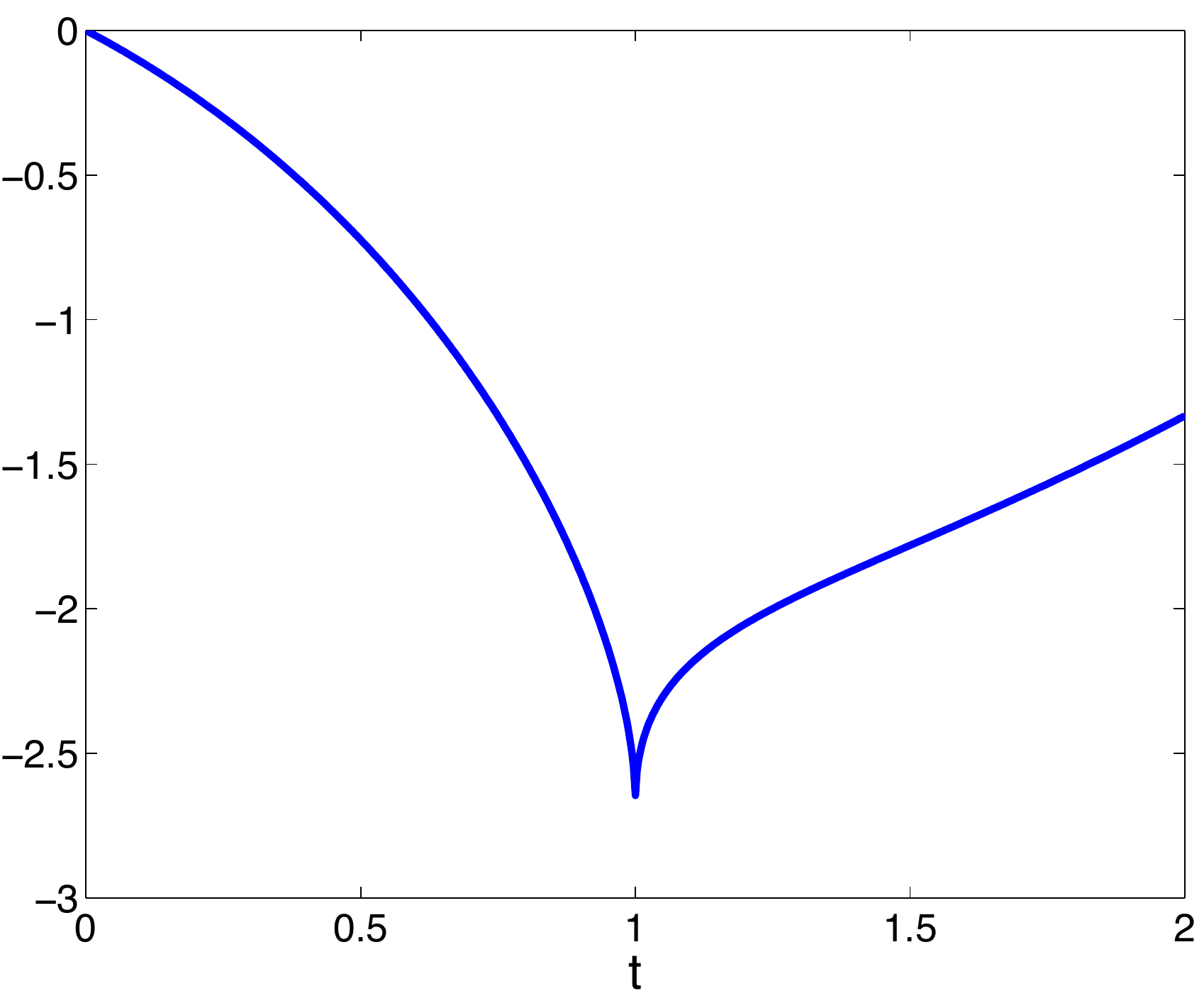} 
\end{tabular}
\caption{Normal derivative $\partial_{\nu}y\vert_{q_T}=y_x(1,t)$ on $(0,T)$ associated to $\mu(x)=1_{(0.2,0.5)}(x)$ (left) and $\mu(x)=1/\sqrt{x}$ (right).}
\label{fig:ynux}
\end{figure}

Remarkably, the approach also provides good reconstruction of the solution $y$ when the observation $y_{\nu,obs}$ is obtained from less regular $\mu$ function. We consider the following rather stiff examples, respectively in $L^2(\Omega)$ and $H^{-1}(\Omega)$:
  \begin{equation}
  \nonumber
\begin{aligned}
&(\textbf{EX4})\qquad   \mu(x)=1_{[a,b]}(x)\in L^2(0,1),  \qquad a=0.2, b=0.5\\
&(\textbf{EX5})\qquad   \mu(x)=1/\sqrt{x}\in H^{-1}(0,1).
\end{aligned}
\end{equation}
The corresponding normal derivatives $y_x(1,\cdot)$ in $H^1(0,T; L^2(\partial \Omega))$ (see \cite{yamamoto}) and in $L^2(0,T)$ respectively, are depicted on Figure \ref{fig:ynux}.  

\begin{figure}[ht]
\centering
\begin{tabular}{ccc}
\includegraphics[width=0.46\textwidth]{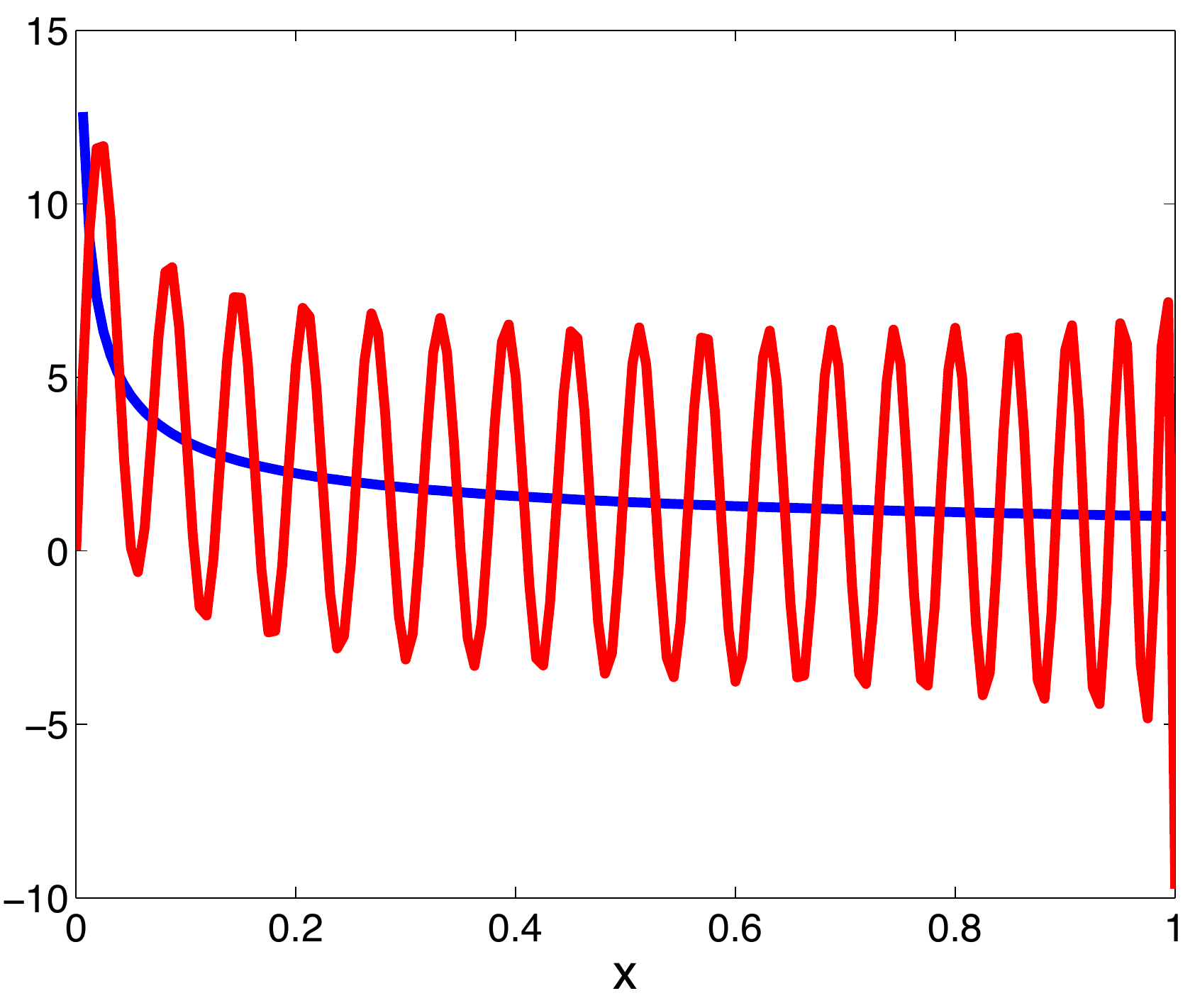} & \qquad & 
\includegraphics[width=0.46\textwidth]{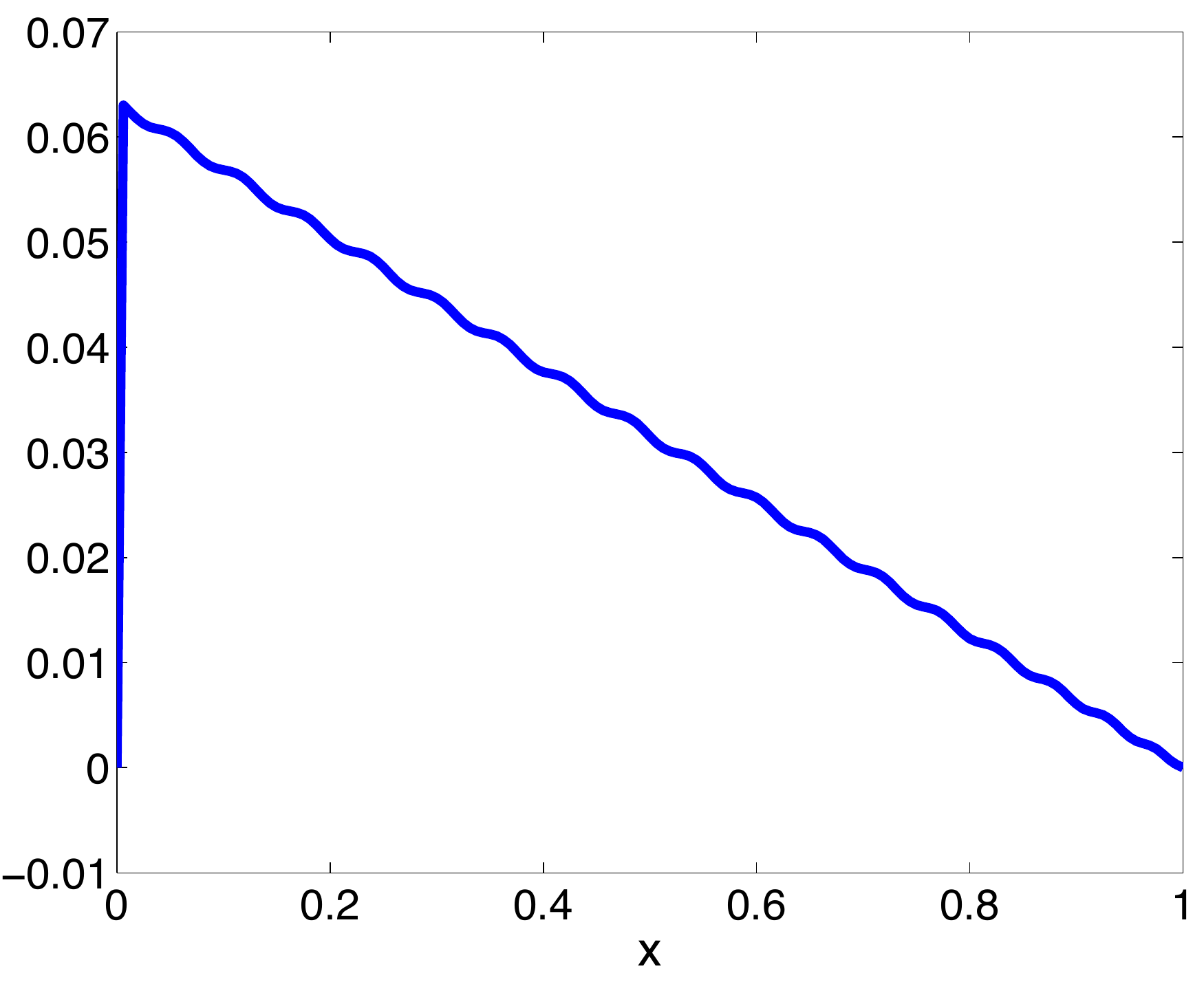} 
\end{tabular}
\caption{\textbf{EX5} -\textbf{Left}: Function $\mu$ (blue) and its approximation $\mu_h$ (red) along $\Omega$;   \textbf{Right}:   $\frac{(-\Delta)^{-1}(\mu-\mu_h)}{\Vert(-\Delta)^{-1}\mu\Vert_{H^1_0(\Omega)}}$ along $\Omega$.}
\label{fig:ex5_mu}
\end{figure}

Table \ref{tab:ex4_rh4_T2} reports some norms with respect to $h$ for the example \textbf{(EX4)}. We observe the convergence as $h\to 0$ with as expected a lower rate : precisely, we compute 
\begin{equation}\nonumber
\frac{\Vert y-y_h\Vert_{L^2(Q_T)}}{\Vert y\Vert_{L^2(Q_T)}}=\mathcal{O}(h^{1.37}), \qquad \frac{\Vert \mu-\mu_h \Vert_{H^{-1}(\Omega)} }{ \Vert \mu\Vert_{H^{-1}(\Omega)}}=\mathcal{O}(h^{0.95}).
\end{equation}
For  $h=8.83\times 10^{-3}$, Figure \ref{fig:ex5_mu} depicts the functions $\mu$ and $\mu_h$ leading to a relative error equal to $2.16\times 10^{-2}$. The approximation $\mu_h$ oscillates around $\mu$ and suggests that the convergence (in agreement with Theorem 3.1) may not hold pointwise but in a weaker (average) sense. Again, this weak convergence of the source term is enough to reconstruct with robustness the solution $y$.
The function $-\Delta^{-1}(\mu-\mu_h)/\Vert -\Delta^{-1} \mu\Vert_{H^1_0(\Omega)}$ of magnitude $10^{-3}$ suggests the efficiency of the method to reconstruct the spatial term $\mu$ of the source $f$.








\begin{table}[http]
\centering
\begin{tabular}{|c|ccccc|}
\hline
$h$  &  $7.07\times 10^{-2}$  & $3.53\times 10^{-2}$ & $1.72\times 10^{-2}$ & $8.83\times 10^{-3}$ & $7.07\times 10^{-3}$  \tabularnewline

\hline 

$\frac{\Vert y-y_h\Vert_{L^2(Q_T)}}{\Vert y\Vert_{L^2(Q_T)}}$ & $4.72\times 10^{-3}$ & $2.34\times 10^{-3}$ &  $5.58\times 10^{-4}$ &  $2.96\times 10^{-4}$ & $2.18\times 10^{-4}$\tabularnewline

$\frac{\Vert \mu-\mu_h\Vert_{H^{-1}(\Omega)}}{\Vert \mu\Vert_{H^{-1}(\Omega)}}$ & $1.53\times 10^{-1}$ & $7.88\times 10^{-2}$ & $2.5\times 10^{-2}$   & $2.16\times 10^{-2}$ & $1.76\times 10^{-2}$ \tabularnewline

$\frac{\Vert \partial_{\nu}(y-y_h)\Vert_{L^2(\Gamma_T)}}{\Vert \partial_{\nu} y\Vert_{L^2(\Gamma_T)}}$  & $2.76\times 10^{-3}$  & $1.08\times 10^{-3}$ & $3.57\times 10^{-4}$ & $1.4\times 10^{-4}$ &  $1.01\times 10^{-4}$ \tabularnewline

$\Vert \lambda_h\Vert_{L^2(Q_T)}$  & $6.81\times 10^{-2}$ & $5.07\times 10^{-2}$ & $3.4\times 10^{-2}$  & $2.56\times 10^{-2}$  & $2.31\times 10^{-2}$\tabularnewline

\hline
\end{tabular}
\caption{Example \textbf{EX4} - BFS element - $r=h^4$ - $T=2$.}
\label{tab:ex4_rh4_T2}
\end{table}

\begin{figure}[ht]
\centering
\begin{tabular}{ccc}
\includegraphics[width=0.46\textwidth]{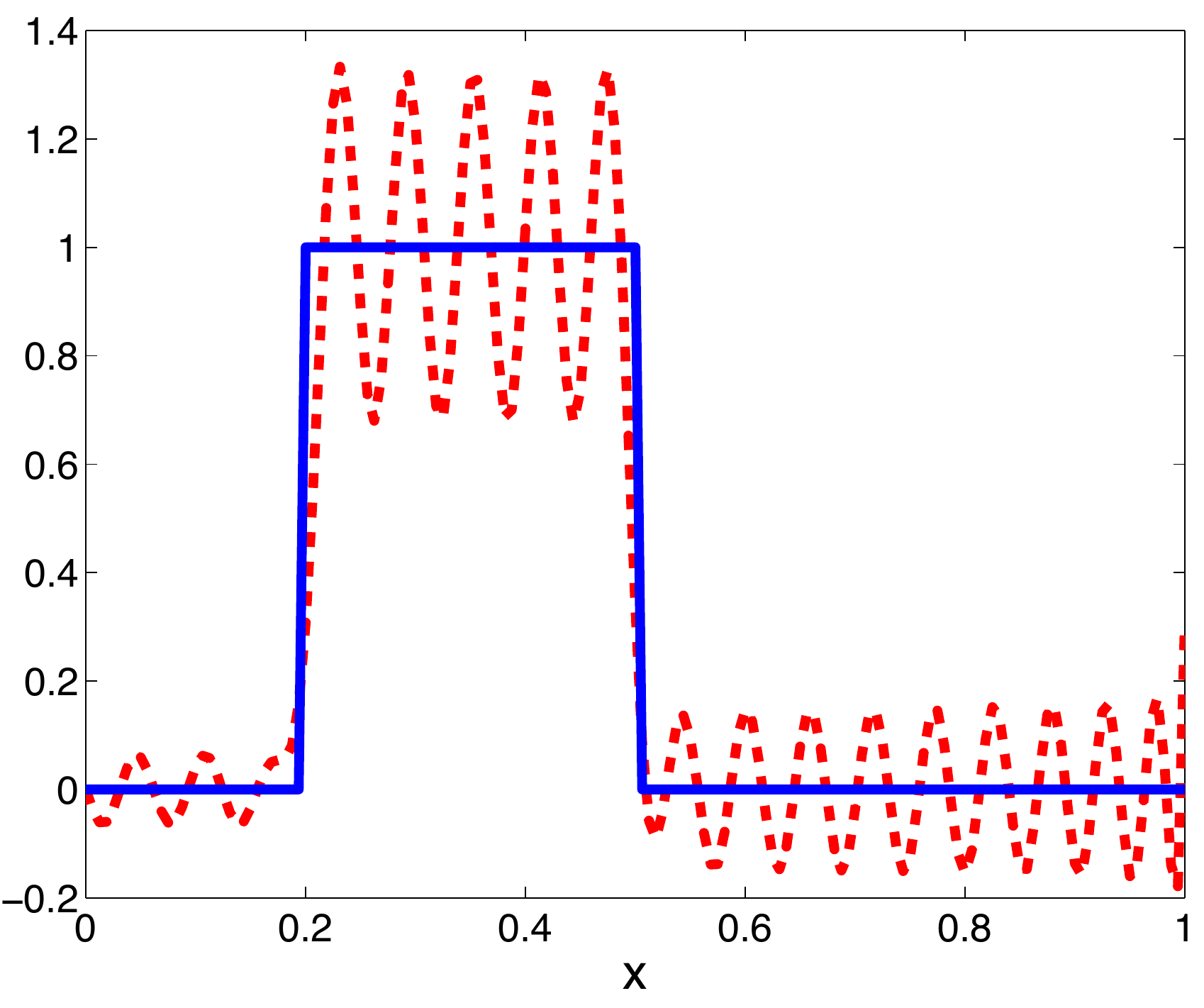} & \qquad & 
\includegraphics[width=0.46\textwidth]{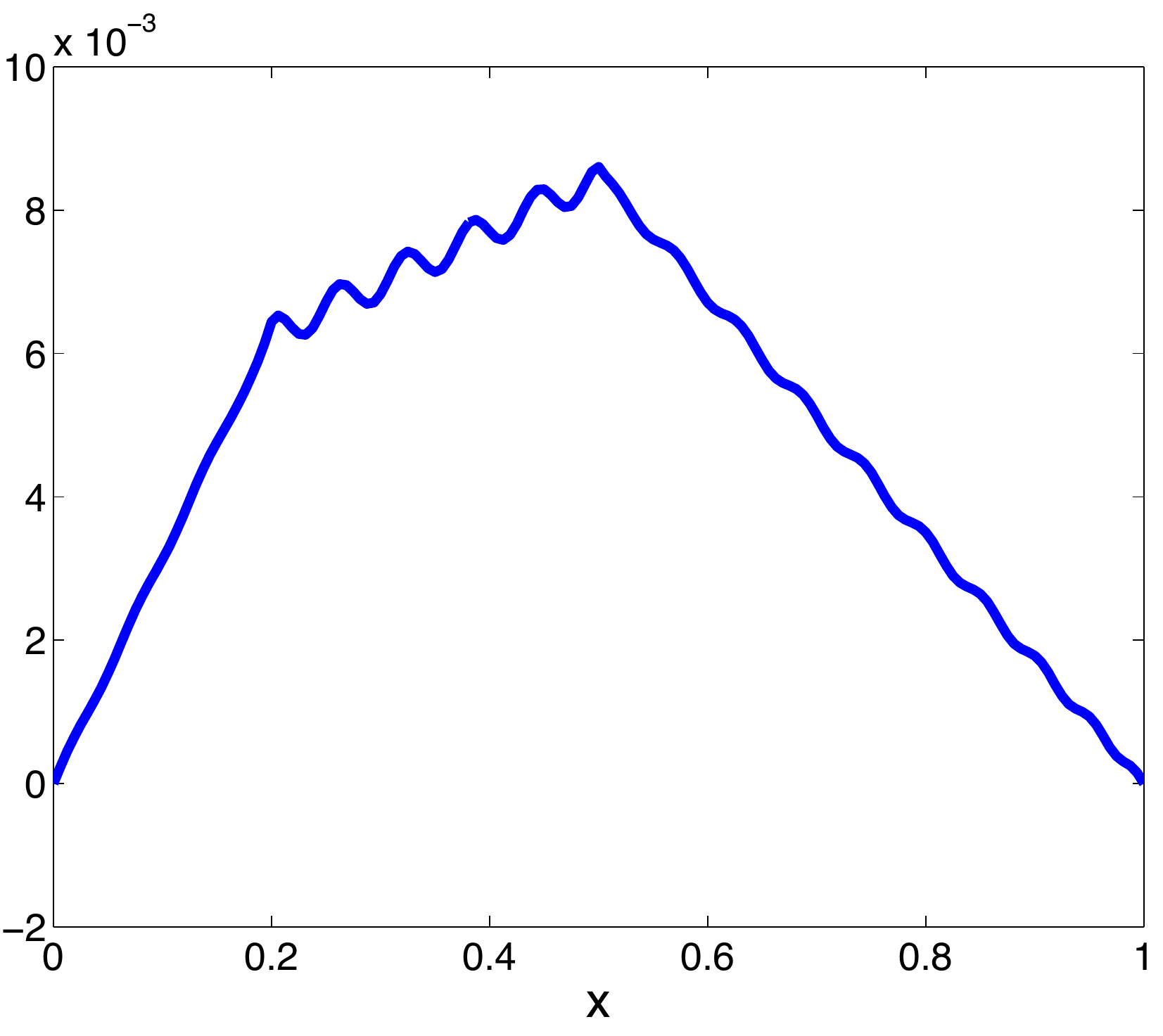} 
\end{tabular}
\caption{\textbf{EX4} -\textbf{Left}: Function $\mu$ (blue) and its approximation $\mu_h$ (red) along $\Omega$;   \textbf{Right}:   $\frac{(-\Delta)^{-1}(\mu-\mu_h)}{\Vert(-\Delta)^{-1}\mu\Vert_{H^1_0(\Omega)}}$ in $\Omega$.}
\label{fig:ex4_mu}
\end{figure}








Similar remarks can be made for the example \textbf{EX5}: we refer to Table \ref{tab:ex5_rh4_T2} and Figure \ref{fig:ex5_mu}.

\begin{table}[http]
\centering
\begin{tabular}{|c|ccccc|}
\hline
$h$  &  $7.07\times 10^{-2}$  & $3.53\times 10^{-2}$ & $1.72\times 10^{-2}$ & $8.83\times 10^{-3}$ & $7.07\times 10^{-3}$  \tabularnewline

\hline 

$\frac{\Vert y-y_h\Vert_{L^2(Q_T)}}{\Vert y\Vert_{L^2(Q_T)}}$ & $1.82\times 10^{-2}$ & $7.74\times 10^{-2}$ & $3.18\times 10^{-3}$ &  $1.87\times 10^{-3}$ & $1.17\times 10^{-3}$\tabularnewline

$\frac{\Vert \mu-\mu_h\Vert_{H^{-1}(\Omega)}}{\Vert \mu\Vert_{H^{-1}(\Omega)}}$ & $31.44$ & $11.27$ & $3.96$ & $1.42$ & $1.02$ \tabularnewline

$\frac{\Vert \partial_{\nu}(y-y_h)\Vert_{L^2(\Gamma_T)}}{\Vert \partial_{\nu} y\Vert_{L^2(\Gamma_T)}}$  & $2.49\times 10^{-1}$  & $2.88\times 10^{-1}$ & $1.01\times 10^{-1}$ & $5.68\times 10^{-2}$ & $4.72\times 10^{-2}$ \tabularnewline

$\Vert \lambda_h\Vert_{L^2(Q_T)}$   & $2.99\times 10^{-1}$ & $2.35\times 10^{-1}$ & $1.91\times 10^{-1}$  & $1.62\times 10^{-1}$  & $1.52\times 10^{-1}$\tabularnewline

\hline
\end{tabular}
\caption{Example \textbf{EX5} - BFS element - $r=h^4$ - $T=2$.}
\label{tab:ex5_rh4_T2}
\end{table}

\section{Concluding remarks and perspectives}  \label{sec_conclusion}

The mixed formulations we have introduced here in order to address inverse problems for hyperbolic equations seems original. These formulations are nothing else than the Euler systems associated to least-squares type functionals and depend on both the state to be reconstruct and a Lagrange multiplier. This Lagrange multiplier is introduced to take into account the state constraint $Ly-f=0$ and turns out to be the controlled solution of a hyperbolic equation with the source term $(\partial_{\nu}y-y_{\nu,obs})\,1_{\Gamma_T}$. This approach, recently used in a controllability context in \cite{NC-AM-mixedwave}, leads to a variational problem defined over time-space functional Hilbert spaces, without distinction between the time and the space variable. 
The main ingredients are, first a unique continuation type property for the hyperbolic equation (assuming some geometric conditions on the measurement zone) allowing to prove the well-posedness of the mixed formulation, and second, a (strong) generalized observability inequality, allowing to quantify the global reconstruction of the solution.  

At the practical level, the discrete mixed time-space formulation is solved in a systematic way in the framework of the finite element theory. The approximation is conformal allowing to obtain the strong convergence of the approximation as the discretization parameters tends to zero. In particular, we emphasize that there is no need, contrarily to the classical approach, to prove some uniform discrete observability inequality: we simply use the observability inequality on the finite dimensional discrete space. The resolution amounts to solve a sparse symmetric linear system : the corresponding matrix can be preconditioned if necessary, and may be computed once for all as it does not depend on the observation $y_{\nu,obs}$. Eventually, the space-time discretization of the domain allows an adaptation of the mesh so as to reduce the computational cost and capture the main features of the solutions.  Similarly, this space-time formulation is very appropriate to the non-cylindrical situation.

In agreement with the theoretical convergence, the numerical experiments reported here display a very good behavior and robustness of the approach: the reconstructed approximate solution converges strongly to the solution of the hyperbolic equation associated to the available observation. Remark that from the continuous dependence of the solution with respect to the observation, the method is robust with respect to the possible noise on the data. 


Eventually, since the mixed formulations rely essentially on a generalized observability inequality, it may be employed to any other observable systems for which such property is available :
we mention notably the parabolic case usually badly conditioned -- in view of regularization property -- and for which direct and robust methods are certainly very advantageous. This is in contrast with the Luenberger approach (mentioned in the introduction) which assume the reversibility in time of the equation  We also mention that this kind of approach may be used to reconstruct potential and coefficient. 


\bibliographystyle{siam}




\end{document}